\documentstyle[12pt,amstex,amssymb,psfig,epic,eepic]{book}
\let\egtrm\relax
\let\tenrm\relax
\let\twlrm\relax

\def\finpreuve
{\hskip 3pt \vrule height6pt width6pt depth 0pt}

\newcommand{\CC}{{\Bbb C}}
\newcommand{\ZZ}{{\Bbb Z}}
\newcommand{\RR}{{\Bbb R}}
\newcommand{\QQ}{{\Bbb Q}}
\newcommand{\PP}{{\Bbb P}}
\newcommand{\codim}{\operatorname{codim}}
\newcommand{\Pic}{\operatorname{Pic}}
\newcommand{\NE}{\operatorname{NE}}

\newcommand{\Hilb}{\operatorname{Hilb}}
\newcommand{\mult}{\operatorname{mult}}
\newcommand{\vol}{\operatorname{vol}}
\newcommand{\divi}{\operatorname{div}}
\newcommand{\pr}{\operatorname{pr}}
\newcommand{\con}{\operatorname{cont}}
\newcommand{\ima}{\operatorname{Im}}

\newcounter{subsub}[subsection]
\def\thesubsub{\thesubsection .\arabic{subsub}}
\def\subsub#1{\addtocounter{subsub}{1}\par\vspace{3mm}
\noindent{\bf \thesubsub ~ #1 }\par\vspace{2mm}}

\def\pr{\mathop{\rm pr}\nolimits}

\def\hfl#1#2{\smash{\mathop{\hbox to 12mm{\rightarrowfill}}
\limits^{\scriptstyle#1}_{\scriptstyle#2}}}

\def\limind{\mathop{\oalign{lim\cr
\hidewidth$\longrightarrow$\hidewidth\cr}}}

\begin{document}

\tableofcontents

\newpage

\chapter*{Introduction}

Le sujet central de cette th\`ese est l'\'etude
de certaines propri\'et\'es d'une classe
de vari\'et\'es analytiques complexes compactes~: les
vari\'et\'es de Moishezon. Ces derni\`eres sont particuli\`erement
int\'eressantes car elles forment la plus petite classe de vari\'et\'es
complexes stable par application bim\'eromorphe et contenant les
vari\'et\'es projectives. Il est bien connu depuis K.\ Kodaira
que les vari\'et\'es projectives sont caract\'eris\'ees par
l'existence d'un fibr\'e en droites ample ou de
fa\c con \'equivalente d'un fibr\'e en droites
muni d'une m\'etrique hermitienne \`a courbure
strictement positive. Le fil conducteur de cette th\`ese
est l'\'etude de l'existence ou de
l'inexistence de fibr\'es en droites v\'erifiant
des propri\'et\'es de positivit\'e faible sur les vari\'et\'es de
Moishezon.

L'\'etude que nous
avons faite est divis\'ee en deux parties~: un point de vue
analytique suivant et g\'en\'eralisant une d\'emarche pr\'esente dans
certains travaux de J.-P.\ Demailly et Y.-T.\ Siu, et
un point de vue plus alg\'ebrique reposant sur l'utilisation de
la th\'eorie de Mori, d\'emarche pr\'esente dans certains travaux
de J.\ Koll\'ar et T.\ Peternell.

\bigskip

\bigskip

\noindent{\bf \'Etude analytique}

\medskip

Cette \'etude a d\'emarr\'e avec les travaux de J.-P.\ Demailly et
Y.-T.\ Siu qui,
r\'epondant \`a une conjecture de H.\ Grauert
et O.\ Riemenschneider, ont donn\'e ind\'ependamment des conditions analytiques
suffisantes (existence
de fibr\'es en droites \`a courbure
semi-positive et g\'en\'eriquement positive)
pour qu'une vari\'et\'e complexe compacte soit de Moishezon.

Une de nos motivations vient du fait qu'aucune de ces conditions n'est
n\'ecessaire, comme le montre l'\'etude de constructions r\'ecentes.

L'un des premiers r\'esultats de cette th\`ese est de donner une
caract\'erisation analytique des vari\'et\'es de Moishezon.
Pour cela, nous montrons, et c'est le th\'eor\`eme principal
de la premi\`ere partie de notre travail, que les
in\'egalit\'es de Morse holomorphes de J.-P.\ Demailly
se g\'en\'eralisent au cas d'un fibr\'e en droites $E$ muni
d'une m\'etrique singuli\`ere $h$ au dessus d'une
vari\'et\'e complexe compacte $X$.
Nos in\'egalit\'es donnent une estimation
asymptotique de la dimension des groupes de cohomologie
\`a valeurs dans les puissances tensorielles $E^{\otimes k}$,
tordues
par une suite de faisceaux d'id\'eaux $ {\cal I}_{k} (h)$
naturellement associ\'ee
aux singularit\'es de la m\'etrique $h$~: la suite
des faisceaux d'id\'eaux multiplicateurs de Nadel.
La pr\'esence de ces faisceaux d'id\'eaux constitue le ph\'enom\`ene
nouveau par rapport au cas o\`u la m\'etrique est lisse.
Comme dans ce dernier cas, l'estimation fait intervenir des
int\'egrales de
la courbure $\Theta(E)$.
Notre r\'esultat est le suivant~:

\bigskip

\noindent {\bf Th\'eor\`eme } {\em Si la m\'etrique $h$
a des singularit\'es analytiques, alors
pour tout fibr\'e $F$ de rang $r$ et pour
tout $q$ compris entre
$0$ et $n = \dim (X)$, on a~:
$$ \sum _{j=0}^{q}(-1)^{q-j} \dim H^{j}(X,{\cal O}(E^{k}\otimes F) \otimes
{\cal I}_{k}(h)) \leq r\frac{k^{n}}{n!} \int _{X(\leq q,E)} (-1)^{q} \Theta
(E)^{n} + o(k^{n}) $$
(avec \'egalit\'e si $q=n$), o\`u $X(\leq q,E)$ d\'esigne l'ouvert
de $X$ des points lisses de la m\'etrique d'indice inf\'erieur \`a $q$.
}

\bigskip

Ce r\'esultat est \`a mettre en parall\`ele avec la g\'en\'eralisation
du th\'eor\`eme de Ka\-wa\-ma\-ta-Viehweg donn\'ee par A.\ Nadel.

Nous montrons ensuite, g\'en\'eralisant un r\'esultat de
S.\ Ji et B.\ Shiffman obtenu ind\'ependamment
et simultan\'ement au n\^otre, que les crit\`eres de J.-P.\ Demailly et
Y.-T.\ Siu deviennent,
dans ce cadre plus souple, n\'ecessaires et
suffisants. Donnons par exemple le~:

\bigskip

\noindent {\bf Th\'eor\`eme } {\em Une vari\'et\'e compacte $X$ de
dimension $n$
est de Moishezon si et
seulement s'il existe sur $X$ un courant ferm\'e $T$ de bi-degr\'e $(1,1)$
tel que~:

 (i) $ \displaystyle{\{ T \} \in H^{2}(X,\ZZ) }$,

 (ii) $\displaystyle{ T= \frac{i}{\pi} \partial \overline{\partial} \varphi
+ \alpha }$,
o\`u $\varphi$ est une fonction r\'eelle \`a singularit\'es analytiques
et o\`u
$\alpha$ est un repr\'esentant
$ {\cal C}^{\infty}$ de $\{ T \}$,

(iii) $\displaystyle{\int_{X(\leq 1,T)} T^{n} > 0}$ o\`u
l'int\'egrale est prise
sur les points lisses du courant $T$.
}

\bigskip

Comme nous l'avons d\'ej\`a mentionn\'e, ce type de crit\`ere
a la propri\'et\'e d'\^etre invariant par morphisme bim\'eromorphe.

\bigskip

\bigskip

\noindent{\bf \'Etude alg\'ebrique}

\medskip

La deuxi\`eme partie de cette th\`ese consiste \`a \'etudier
en d\'etail la classe des vari\'et\'es de Moishezon dont le groupe
de Picard est $\ZZ$, et dont le fibr\'e canonique $K_X$
est gros (``big").

Une de nos motivations provient d'un r\'esultat de J.\ Koll\'ar affirmant
qu'en dimension $3$, et sous les hypoth\`eses pr\'ec\'edentes, le
fibr\'e canonique est alors num\'eriquement effectif (nef).
Il n'est donc pas possible de trouver dans cette classe des
vari\'et\'es de Moishezon de dimension $3$ ne v\'erifiant pas
les crit\`eres de J.-P.\ Demailly et Y.-T.\ Siu.
Nous montrons que ceci n'est plus vrai en dimension sup\'erieure
ou \'egale \`a $4$ en construisant explicitement des exemples~:

\bigskip

\noindent {\bf Th\'eor\`eme }
{\em Pour tout entier $n$ sup\'erieur
ou \'egal \`a $4$, il existe des vari\'et\'es de Moishezon $X$,
non projectives,
de dimension $n$ v\'erifiant~:

(i) $\Pic (X) = \ZZ$, (ii) $K_X$ est gros, (iii) $K_X$ n'est pas nef.}

\bigskip

La construction donnant ce r\'esultat montre que les vari\'et\'es
$X$ obtenues deviennent projectives apr\`es un \'eclatement
le long d'une sous-vari\'et\'e isomorphe \`a $\displaystyle{\PP ^{n-2}}$.

Plus g\'en\'eralement, un r\'esultat fondamental de B.\ Moishezon
affirme qu'une vari\'et\'e
de Moishezon peut \^etre rendue projective apr\`es une succession
finie d'\'eclatements le long de sous-vari\'et\'es lisses. Ce r\'esultat
difficile ne donne malheureusement aucune indication quant au choix explicite
des sous-vari\'et\'es en question.
Gr\^ace \`a l'utilisation de la c\'el\`ebre th\'eorie de Mori
sur un mod\`ele projectif, nous avons \'etudi\'e le centre
de l'\'eclatement en toutes dimensions~:

\bigskip

\noindent {\bf Th\'eor\`eme }
{\em Soit $X$ une vari\'et\'e
de Moishezon (non projective) de dimension $n$, avec $\Pic (X) = \ZZ$
et $K_X$ gros.
Supposons de plus que $X$ est rendue projective apr\`es \'eclatement
le long d'une sous-vari\'et\'e $Y$ lisse.

\noindent Alors, si $K_X$ n'est pas nef, on a
$\displaystyle{ \dim Y > \frac{n-1}{2}.}$
}

\bigskip

En dimension
$4$, ce r\'esultat peut \^etre pr\'ecis\'e, y compris dans
le cas o\`u le fibr\'e canonique est nef~:

\bigskip

\noindent {\bf Th\'eor\`eme }
{\em Soit $X$ une vari\'et\'e
de Moishezon (non projective) de dimension $4$, avec $\Pic (X) = \ZZ$
et $K_X$ gros.
Supposons de plus que $X$ est rendue projective apr\`es \'eclatement
le long d'une sous-vari\'et\'e $Y$ lisse.

\noindent Alors $Y$ est n\'ecessairement
une surface.
Autrement dit, et dans cette situation particuli\`ere,
il ne suffit pas d'\'eclater une courbe pour ``rentrer dans
le monde projectif".}

\bigskip

Nous avons vu pr\'ec\'edemment que $K_X$ n'est pas n\'ecessairement
nef \`a partir de la dimension $4$. Le r\'esultat suivant montre
que l'exemple que nous avons construit
est, en un sens, le seul possible en dimension $4$ dans le cas
o\`u $K_X$ n'est pas nef.

\bigskip

\noindent {\bf Th\'eor\`eme }
{\em Sous les hypoth\`eses pr\'ec\'edentes et si $K_{X}$ n'est
pas nef, alors le couple $\displaystyle{(Y,N_{Y/X})}$ est isomorphe \`a
$\displaystyle{(\PP ^2, {\cal O}_{\PP ^{2}}(-1)^{\oplus 2})}$.}

\bigskip

Ces r\'esultats pr\'ecis sont accessibles en dimension $4$ car les
contractions de Mori ont \'et\'e \'etudi\'ees par T.\ Ando, Y.\ Kawamata
et M.\ Beltrametti.
Nos r\'esultats peuvent \^etre vus comme un premier pas vers
une analyse du caract\`ere non projectif des vari\'et\'es de Moishezon
de dimension sup\'erieure ou \'egale \`a $4$~; la situation
en dimension $3$ \'etant maintenant assez bien comprise suite aux
travaux de J.\ Koll\'ar et T.\ Peternell.

\bigskip

\bigskip

\noindent{\bf Plan du texte}

\medskip

- le chapitre 1 est un chapitre de pr\'eliminaires~; il contient
une description pr\'ecise des motivations et des objets utilis\'es
dans le reste de la th\`ese. Nous y d\'etaillons en particulier
une d\'emonstration du th\'eor\`eme de Siegel et donnons
les constructions de I.\ Nakamura et K.\ Oguiso montrant que les
crit\`eres analytiques de J.-P.\ Demailly et Y.-T.\ Siu
ne sont pas n\'ecessaires en g\'en\'eral.

\medskip

- le chapitre 2 est consacr\'e \`a l'\'etude analytique. Nous rappelons
les premiers r\'esultats li\'es aux m\'etriques singuli\`eres,
\'enon\c cons et d\'emontrons notre version des in\'egalit\'es
de Morse holomorphes. On en d\'eduit les caract\'erisations
analytiques
des vari\'et\'es de Moishezon.
Enfin, nous donnons une version alg\'ebrique singuli\`ere
des in\'egalit\'es de Morse.

\medskip

- le chapitre 3 est consacr\'e \`a l'\'etude alg\'ebrique.
Nous commen\c cons par rappeler le r\'esultat de J.\ Koll\'ar
en dimension $3$, puis nous faisons une \'etude des vari\'et\'es
de Moishezon \`a groupe de Picard $\ZZ$, \`a fibr\'e canonique
gros et devenant projectives apr\`es un seul \'eclatement de
centre lisse et projectif. Nous obtenons
une restriction sur la dimension du centre de l'\'eclatement.
En dimension $4$, cette restriction implique que ce dernier
est n\'ecessairement une surface.
Nous d\'ecrivons alors notre exemple et montrons qu'en dimension
$4$, cette construction est essentiellement la seule dans le cas
o\`u le fibr\'e canonique n'est pas nef.

\medskip

Mentionnons que les chapitres 2 et 3 sont dans une large
mesure ind\'ependants et peuvent \^etre lus dans un ordre
quelconque.


\chapter{Pr\'eliminaires}

Ce chapitre a pour but d'introduire les principales
notions utilis\'ees par la suite, de pr\'esenter
les premi\`eres motivations en d\'etail et de rappeler
un certain nombre de r\'esultats auxquels nous nous r\'ef\'erons
dans les chapitres suivants.

\section{\! Quelques rappels de g\'eom\'etrie
analytique complexe}

\subsection{Vari\'et\'es, fibr\'es vectoriels}
Pr\'ecisons tout d'abord que dans toute cette th\`ese,
et sauf mention explicite
du contraire, le mot {\bf vari\'et\'e} sera utilis\'e pour d\'esigner
une {\bf vari\'et\'e analytique
complexe {\em non singuli\`ere}} suppos\'ee de plus connexe.
Pour toutes les notions introduites ici, nous renvoyons
de fa\c con g\'en\'erale \`a \cite{G-H78}.

\medskip

Un {\bf fibr\'e vectoriel} complexe $F$ au dessus d'une vari\'et\'e
$X$ est la donn\'ee d'une vari\'et\'e $F$
et d'une application $\displaystyle{\pi : F \to X}$ de sorte
qu'il existe un recouvrement de $X$ par des
ouverts trivialisants $U_{\alpha}$ et des isomorphismes
(appel\'es trivialisations)
$$\theta _{\alpha} : \pi ^{-1}(U_{\alpha}) \to U_{\alpha} \times \CC ^r$$
respectant la structure d'espace vectoriel des fibres, i.e
$$\theta _{\alpha \beta}(x,\xi)
:= \theta _{\alpha} \circ \theta _{\beta} ^{-1}
(x,\xi) = (x, g_{\alpha \beta}(x)\xi )$$ o\`u
$g_{\alpha \beta}$ est une application holomorphe
sur $U_{\alpha} \cap U_{\beta}$ \`a valeurs dans le
groupe des matrices complexes inversibles de taille $r$.
Nous notons $\xi _x$ un point de $F$ au dessus
du point $x$ de $X$ (i.e tel que $\pi (\xi _x) =x$).
L'entier $r$ est le {\bf rang} du fibr\'e $F$.
Si $r = 1$, on parle de {\bf fibr\'e en droites}.

Un exemple important de fibr\'e en droites
est le fibr\'e ${\cal O}(D)$
associ\'e \`a un diviseur $D$ sur $X$~: si $D$ est
un diviseur irr\'eductible donn\'e sur $U_{\alpha}$
par l'\'equation $f_{\alpha}=0$, le fibr\'e
${\cal O}(D)$ est le fibr\'e associ\'e au cocycle
$$g_{\alpha \beta} = \frac{f_{\alpha}}{f_{\beta}} \ .$$

\medskip

Toutes les constructions d'alg\`ebre lin\'eaire
s'\'etendent aux fibr\'es~: dual,
produit tensoriel, produit ext\'erieur.
Ainsi, un exemple fondamental de fibr\'e
en droites
sur une vari\'et\'e $X$ de dimension $n$ est le
{\bf fibr\'e canonique} d\'efini par
$$K_X := \det (T^{\ast}X) = \bigwedge ^n T^{\ast}X,$$ o\`u $T^{\ast}X$
d\'esigne le fibr\'e cotangent, dual du fibr\'e
tangent holomorphe $TX$ de $X$.

Un autre exemple important de fibr\'e en droites
est le fibr\'e ${\cal O}_{\PP ^n}(1)$ sur l'espace
projectif $\PP ^n$~: en associant \`a un point $[x]$
de $\PP ^n$ la droite $\CC x$, on construit
un sous-fibr\'e en droites du fibr\'e trivial $\PP ^n \times \CC ^{n+1}$~;
le dual de ce fibr\'e en droites est par d\'efinition
le fibr\'e ${\cal O}_{\PP ^n}(1)$.

L'ensemble des
fibr\'es en droites, modulo isomorphisme, sur
une vari\'et\'e $X$ est naturellement muni
d'une structure de groupe pour le produit
tensoriel~: on l'appelle {\bf groupe de Picard}
de $X$ et on le note $\Pic (X)$. Mentionnons ici
que nous identifions suivant l'usage un fibr\'e
en droites $E$ au faisceau inversible
${\cal O}(E)$ des germes de sections holomorphes de $E$.
La $k$-i\`eme puissance tensorielle d'un fibr\'e
en droites $E$ sera not\'ee indiff\'eremment
$E^{\otimes k}$, $E^k$,${\cal O}(kE)$ ou m\^eme $kE$.

\medskip

Un fibr\'e vectoriel $E$ peut \^etre muni d'une m\'etrique
hermitienne ${\cal C}^{\infty}$, on parle
alors de fibr\'e vectoriel {\bf hermitien} et on note
g\'en\'eralement $h$ une telle m\'etrique~: elle correspond
\`a la donn\'ee d'une forme hermitienne sur chaque
fibre $E_x$ de $E$, d\'ependant de fa\c con ${\cal C}^{\infty}$
de $x$.

Dans le cas particulier d'un fibr\'e en droites,
une m\'etrique hermitienne est
donn\'ee localement
sur un ouvert trivialisant $U_{\alpha}$
par
$$h(\xi _{x}) = ||\xi _x||_h := |\xi|\exp(-\varphi _{\alpha} (x))$$
(la fonction $\exp(-\varphi _{\alpha})$
est appel\'ee {\bf poids} de la m\'etrique $h$ dans la trivialisation
$\theta _{\alpha}$)
o\`u la fonction r\'eelle $\varphi _{\alpha}$ est
de classe ${\cal C}^{\infty}$ sur $U_{\alpha}$.

Lorsque le fibr\'e tangent $TX$ est muni d'une m\'etrique hermitienne,
on dit que la vari\'et\'e $X$ est hermitienne.
Comme il est d'usage, nous identifions toujours
la donn\'ee d'une m\'etrique hermitienne sur une vari\'et\'e
$X$ \`a celle de la $(1,1)$ forme r\'eelle, g\'en\'eralement
not\'ee $\omega$, qui lui est naturellement associ\'ee
($\omega$ est \`a un facteur $-2$ pr\`es la partie imaginaire de la
m\'etrique). Ainsi, une vari\'et\'e est {\bf k\"ahl\'erienne}
si elle poss\`ede une m\'etrique hermitienne pour laquelle $\omega$
est une forme ferm\'ee.

Une vari\'et\'e {\bf projective} est une vari\'et\'e isomorphe
\`a une sous-vari\'et\'e ferm\'ee d'un espace projectif $\PP ^N$.

\medskip

Pour un fibr\'e en droites hermitien $(E,h)$, on note $\Theta (E)$
la $(1,1)$ {\bf forme de courbure} de $(E,h)$~: c'est la forme
r\'eelle d\'efinie globalement sur $X$ et
donn\'ee localement par
$$\displaystyle{ \Theta (E)=
\frac{i}{\pi} \partial \overline{\partial} \varphi _{\alpha}};$$
c'est aussi la courbure de
la connexion de Chern du fibr\'e hermitien $E$.
La classe de cohomologie de $\Theta (E)$
appartient \`a $H^2(X,\ZZ)$ et ne d\'epend
pas de la m\'etrique $h$~; c'est la premi\`ere classe
de Chern de $E$ et elle est not\'ee $c_1(E)$. Remarquons
qu'il n'y a pas de sens \`a parler des valeurs propres
de la forme de courbure, mais que la {\bf signature} de la courbure
(i.e le nombre de ``valeurs propres" nulles, strictement positives
et strictement n\'egatives) est une notion bien d\'efinie sans donn\'ee
suppl\'ementaire.
Par exemple, le fibr\'e ${\cal O}_{\PP ^n}(1)$
muni de la m\'etrique induite de celle de $\CC ^{n+1}$
est \`a courbure strictement positive~: la forme de courbure
est la m\'etrique de Fubini-Study de $\PP ^n$.

\subsection{Th\'eor\`eme de Kodaira}
Nous sommes en mesure d'\'enoncer maintenant le c\'el\`ebre th\'eor\`eme
de plongement de Kodaira \cite{Kod54}~:

\bigskip

\noindent{\bf Th\'eor\`eme (K.\ Kodaira, 1954)}
{\em Une vari\'et\'e compacte $X$ est projective si et seulement si
elle poss\`ede un fibr\'e en droites hermitien $E$ \`a courbure
strictement positive.}

\bigskip

Signalons \'evidemment qu'un sens est ais\'e~: si une
vari\'et\'e est projective, la restriction \`a $X$ du fibr\'e
$\displaystyle{ {\cal O}_{\PP ^N} (1) }$ muni
de sa m\'etrique naturelle ayant pour courbure la forme
de Fubini-Study de $\PP ^N$ donne
le fibr\'e souhait\'e. L'autre sens consiste  \`a montrer
que pour $k$ entier assez grand, il est possible de plonger
$X$ dans l'espace projectif des hyperplans de $H^0(X,E^{\otimes k})$,
o\`u $H^0(X,E^{\otimes k})$
d\'esigne l'espace vectoriel des sections holomorphes globales
de $E^{\otimes k}$.

Rappelons ici qu'un fibr\'e en droites pouvant \^etre
muni d'une m\'etrique \`a courbure strictement positive
est dit {\bf ample}.

C'est le th\'eor\`eme de Kodaira que H.\ Grauert et
O.\ Riemenschneider \cite{GrR70} se
proposaient de g\'en\'eraliser aux vari\'et\'es de Moishezon,
vari\'et\'es que nous introduisons dans le paragraphe suivant.

\section{Quelques rappels sur les vari\'et\'es de Moishezon}

Les vari\'et\'es de Moishezon sont, parmi les vari\'et\'es
compactes,
celles qui poss\`edent le
``plus" de fonctions m\'eromorphes alg\'ebriquement
ind\'ependantes. Cette d\'efinition heuristique est justifi\'ee
par le th\'eor\`eme de Siegel que nous rappelons maintenant.

\subsection{Th\'eor\`eme de Siegel}
En 1955, C.L.\ Siegel d\'emontre le r\'esultat suivant \cite{Sie55}~:

\bigskip

\noindent{\bf Th\'eor\`eme (C.L.\ Siegel, 1955) }
{\em Si $X$ est une vari\'et\'e
compacte
de dimension $n$, alors $X$ poss\`ede au plus $n$ fonctions
m\'eromorphes alg\'ebriquement ind\'ependantes.}

\bigskip

La d\'emonstration originale de C.L.\ Siegel repose sur une
application \'el\'ementaire du
lemme de Schwarz. Il existe maintenant plusieurs d\'emonstrations
diff\'erentes, certaines g\'en\'eralisant cet \'enonc\'e
aux espaces complexes compacts.
Nous en donnons ici une preuve ``moderne" dans le cas
non singulier. Pour cela, nous utilisons un
r\'esultat de P.\ Gauduchon \cite{Gau77}~: {\em toute vari\'et\'e
analytique
complexe compacte
de dimension $n$ poss\`ede une m\'etrique hermitienne
$\omega$ de classe
${\cal C}^{\infty}$ et d'excentricit\'e nulle, i.e
telle que $\displaystyle{ \partial \overline{\partial} (\omega ^{n-1}) =0}$.}

Commen\c cons par montrer le lemme suivant~:

\medskip

\noindent{\bf Lemme} {\em Soit $X$ une vari\'et\'e compacte
que l'on munit d'une m\'etrique de Gauduchon $\omega$
et soit $x_0$ un point de $X$.
Alors, il
existe une constante $C := C(X,x_0,\omega )$ telle que
pour tout fibr\'e en droites hermitien $(E,h)$ au dessus de $X$
et pour toute section holomorphe $s$ de $E$
non identiquement nulle, on ait~:
$$ \mult \, (s,x_0) \leq C \int _X \omega ^{n-1} \wedge \Theta (E).$$}

\medskip

\noindent {\bf D\'emonstration}

Soit $r$ un r\'eel strictement positif fix\'e ``petit"
(de sorte qu'il existe une carte centr\'ee en $x_0$
et contenant la boule $B(x_o,r)$).
Alors, la multiplicit\'e de $s$ en $x_0$ v\'erifie~:
$$  \mult \, (s,x_0) \leq
\frac{ \vol _{n-1}(Z_s \cap B(x_o,r))}{\vol _{n-1}(B_{n-1}(x_o,r))} +o(r) $$
o\`u $\vol _{n-1}$ est le volume $(n-1)$-dimensionnel
mesur\'e avec la m\'etrique $\omega$, et o\`u
$Z_s$ d\'esigne le lieu des z\'eros de $s$. En effet,
la multiplicit\'e de $s$ en $x_0$ est en fait \'egale
\`a la limite d\'ecroissante lorsque $r$ tend vers $0$
de la quantit\'e du membre de droite de l'in\'egalit\'e.
Comme $X$ est compacte, on a {\em a fortiori}~:
$$ \mult \, (s,x_0) \leq C \vol _{n-1}(Z_s)
= C \int_{Z_s} \omega ^{n-1}.$$
Mais l'\'equation de Lelong-Poincar\'e affirme que~:
$$ \frac{i}{\pi}\partial \overline{\partial}\log ||s||
= [ Z_s ] - \Theta (E),$$
o\`u $[ Z_s ]$ d\'esigne le courant d'int\'egration sur
l'ensemble $Z_s$.
Comme $\omega ^{n-1}$ est $\partial \overline{\partial}$-ferm\'ee,
la formule de Stokes donne de suite~:
$$ \int_{Z_s} \omega ^{n-1} = \int _X \omega ^{n-1} \wedge \Theta (E),$$
ce qui prouve le lemme.\finpreuve

\medskip

Ce lemme implique le r\'esultat suivant~:

\medskip

\noindent{\bf Proposition} {\em Soient $X$ une vari\'et\'e compacte
de dimension $n$, $\omega$ une m\'etrique de Gauduchon sur $X$
et $E$ un fibr\'e en droites hermitien au dessus de $X$.
Alors~:

(i) $\dim H^0(X,E) \leq {n+p \choose n}$ o\`u $p$
est la partie enti\`ere
de $\displaystyle{C \int _X \omega ^{n-1} \wedge \Theta (E)}$
et o\`u $C$ est la constante du lemme pr\'ec\'edent,

(ii) $\displaystyle{\dim H^0(X,E^{\otimes k}) \leq
A \left(\int _X \omega ^{n-1} \wedge \Theta (E)\right)^n k^n + o(k^n)}$,
o\`u $A$ est une constante ind\'ependante de $E$ et $k$.
}

\medskip

\noindent {\bf D\'emonstration}

Il suffit de remarquer que
${n+p \choose n}$ est la dimension de l'espace vectoriel
des polyn\^omes
de $n$ variables et
de degr\'e inf\'erieur ou \'egal \`a $p$~: le lemme implique
en effet que l'application lin\'eaire qui \`a une section
holomorphe de $E$ associe son $p$-i\`eme jet en $x_0$
est injective, d'o\`u le point (i). Le point (ii) est cons\'equence du
fait que la forme de courbure
du fibr\'e $(E^{\otimes k},h^k)$ est donn\'ee par
$\displaystyle{ \Theta (E^{\otimes k}) = k \Theta (E)}$.
On applique alors (i) au fibr\'e $(E^{\otimes k},h^k)$.\finpreuve

\medskip

Remarquons que la proposition
pr\'ec\'edente affirme que la dimension de
l'espace vectoriel des sections holomorphes des
puissances $E^{\otimes k}$ d'un fibr\'e en droites sur une vari\'et\'e
compacte de dimension $n$ cro\^\i t au plus comme $k^n$.
Ce fait est bien classique et nous avons estim\'e
la dimension en fonction d'int\'egrales de courbure. Des estimations
bien plus pr\'ecises, valables pour la dimension de tous les groupes de
cohomologie seront donn\'ees par les in\'egalit\'es de Morse holomorphes
de J.-P.\ Demailly dans le paragraphe suivant.

\medskip

\noindent {\bf D\'emonstration du th\'eor\`eme de Siegel}

L'argument est standard~: soient
$\displaystyle{(f_{i})_{1\leq i\leq N}}$ $N$ fonctions m\'eromorphes
alg\'ebriquement
ind\'ependantes sur $X$.
Notons $D$ la somme des diviseurs des p\^oles des $f_{i}$ et
${\cal O}(D)$ le fibr\'e en droites associ\'e.
Rappelons que
$H^{0}(X,{\cal O}(D))$ est isomorphe \`a l'espace vectoriel
des fonctions m\'eromorphes sur $X$ v\'erifiant
$\divi (f) + D \geq 0$. Alors,
si $P$ est un polyn\^ome \`a coefficients complexes
en $N$ variables, de degr\'e total inf\'erieur ou \'egal
\`a $k$, la fonction m\'eromorphe $P(f_{1},\ldots,f_{N})$ est une
section holomorphe
de ${\cal O}(kD)$, et comme les $f_{i}$ sont
alg\'ebriquement ind\'ependantes, on a ${N+k \choose N}$
telles sections
lin\'eairement ind\'ependantes. De l\`a~:
$$  \dim H^{0}(X,{\cal O}(D)^{\otimes k})
\geq {N+k \choose N} \sim_{k \to +\infty} \frac{k^N}{N!}.$$
Avec la proposition, il vient $N \leq n$.\finpreuve

\subsection{\'Eclatements et vari\'et\'es de Moishezon}

Nous commen\c cons ce paragraphe par quelques rappels
sur les \'eclatements. Ces derniers jouent un r\^ole important
dans la th\'eorie des vari\'et\'es de Moishezon et la
construction d'exemples explicites.

\subsub{\'Eclatements}
Une r\'ef\'erence standard est \`a nouveau \cite{G-H78}.

\noindent Si $X$ est une vari\'et\'e, et $Y$ une sous-vari\'et\'e
de $X$ de codimension sup\'erieure ou
\'egale \`a $2$, on construit une vari\'et\'e $\tilde{X}$
appel\'ee {\bf \'eclatement de $X$ le long de $Y$}
en rempla\c cant les points de $Y$
par l'espace des directions normales \`a $Y$ dans $X$.

On note g\'en\'eralement $\pi : \tilde{X} \to X$
l'\'eclatement et $E := \pi ^{-1}(Y)$ le
{\bf diviseur exceptionnel}. La sous-vari\'et\'e $Y$
est appel\'ee {\bf centre de l'\'eclatement}. Par construction,
$\pi$ induit un isomorphisme
$$\pi _{| \tilde{X} \backslash E } : \tilde{X} \backslash E
\to X \backslash Y.$$
\noindent Signalons que le centre de l'\'eclatement
$Y$ peut \^etre
r\'eduit \`a un point.

La restriction de $\pi$ au diviseur exceptionnel $E$
munit $E$ d'une structure
de fibr\'e en espaces projectifs au dessus de $Y$~:
plus pr\'ecis\'ement, $E$ est isomorphe \`a
$\displaystyle{\PP (N^{\ast}_{Y/X})}$ (projectivis\'e
en droites du fibr\'e normal $N_{Y/X}$
suivant la convention de Grothendieck). De plus, le
fibr\'e normal
$\displaystyle{N_{E/\tilde{X}} = {\cal O}(E)_{| E}}$ est
isomorphe au fibr\'e
${\cal O}_{\PP (N_{Y/X}^{\ast})}(-1)$.

Mentionnons aussi que le groupe de Picard
de $\tilde{X}$ est \'egal \`a
$\pi ^*\Pic (X) \oplus \ZZ \cdot {\cal O}(E)$. Par exemple,
le fibr\'e canonique est donn\'e par
$$K_{\tilde{X}} = \pi ^{\ast}K_X + (r-1){\cal O}(E)$$ o\`u
$r$ est la codimension du centre de l'\'eclatement.
Enfin, si $Z$ est une sous-vari\'et\'e de $X$,
non incluse dans le centre de l'\'eclatement $Y$,
alors l'adh\'erence de $\pi ^{-1}(Z \cap (X \backslash Y))$
dans $\tilde{X}$ est appel\'ee {\bf transform\'ee stricte}
de $Z$. Si $Z'$ est la transform\'ee stricte de $Z$, alors
$$\pi _{|Z'}~: Z' \to Z$$
\noindent est l'\'eclatement de $Z \cap Y$ dans $Y$.

\medskip

Le r\'esultat suivant, d\^u \`a A.\ Fujiki et S.\ Nakano \cite{FuN72}
donne un crit\`ere pour qu'un diviseur soit
le diviseur exceptionnel d'un \'eclatement~; ce
crit\`ere est une extension du crit\`ere
de Castelnuovo sur les surfaces. Nous l'utiliserons
souvent lors de la construction d'exemples explicites de
vari\'et\'es de Moishezon non projectives.

\bigskip

{\bf Th\'eor\`eme (A.\ Fujiki, S.\ Nakano, 1972) }
{\em Soit $Z$ une vari\'et\'e et $D$ une sous-vari\'et\'e
de $Z$ de codimension $1$. On suppose que
$D$ est isomorphe \`a $\PP (G)$
o\`u $G$ est un fibr\'e vectoriel sur une vari\'et\'e
$Y$ ; on note $p : \PP (G) \to Y$ la projection.
On suppose enfin que $N_{D/Z}$ est
isomorphe \`a ${\cal O}_{\PP (G)}(-1)$.

Alors, il existe une vari\'et\'e $Z'$ contenant
$Y$ comme sous-vari\'et\'e et une application $\pi : Z \to Z'$
de sorte
que $\pi$ soit l'\'eclatement de $Z'$ le long de $Y$ et
que la restriction de $\pi$ \`a $D$ soit \'egale
\`a $p$.}

\bigskip

Notons qu'au vu de ce qui pr\'ec\`ede, les hypoth\`eses
faites sur le diviseur $D$ dans cet \'enonc\'e sont
\'evidemment n\'ecessaires~: ce crit\`ere remarquable
montre qu'elles sont suffisantes.

\subsub{Vari\'et\'es de Moishezon}

Le th\'eor\`eme de Siegel
motive les d\'efinitions suivantes~:

\medskip

\noindent {\bf D\'efinition } {\em
Un fibr\'e en droites sur une vari\'et\'e
compacte de dimension $n$ est dit {\bf gros} ({\bf big}
en anglais) si la dimension de
l'espace vectoriel des sections holomorphes globales de ses
puissances $E^{\otimes k}$ cro\^\i t exactement comme $k^n$.}

\medskip

\noindent {\bf Remarque }
De fa\c con g\'en\'erale, pour un fibr\'e en droites
$E$ sur une vari\'et\'e compacte $X$, il existe
un entier $\kappa (E)$ tel que la dimension de
$H^0(X,E^k)$ cro\^\i t comme $\displaystyle{k^{ \kappa (E)}}$.
Cet entier (\'egal \`a $-\infty$ si tous les $H^0(X,E^k)$ sont nuls)
est appel\'e {\bf dimension de Kodaira-Iitaka}. Cet entier
est inf\'erieur ou \'egal \`a la dimension de $X$, et selon
ce qui pr\'ec\`ede, le fibr\'e $E$ est gros si et seulement si
$\kappa (E) = \dim X$.

\medskip

\noindent {\bf D\'efinition } {\em Une vari\'et\'e
compacte de dimension $n$ est {\bf de Moishezon} si elle
poss\`ede exactement $n$ fonctions m\'eromorphes alg\'ebriquement
ind\'ependantes ou, de fa\c con \'equivalente, si elle
poss\`ede un fibr\'e en droites gros.}

\medskip

Les premiers exemples de vari\'et\'es de Moishezon sont les
vari\'et\'es projectives. En particulier, toutes les courbes
sont de Moishezon. En fait, le th\'eor\`eme suivant,
difficile et fondamental, montre
qu'une vari\'et\'e de Moishezon n'est pas tr\`es loin d'\^etre
projective. Ce r\'esultat est d\^u \`a B.\ Moishezon \cite{Moi67}~:

\bigskip

\noindent{\bf Th\'eor\`eme (B.\ Moishezon, 1967)}
{\em Une vari\'et\'e compacte est de Moishezon si
et seulement si elle peut
\^etre rendue projective apr\`es un nombre
fini d'\'eclatements de centres
lisses. On peut m\^eme choisir les centres projectifs.}

\bigskip

\noindent {\bf Remarque }
Mentionnons d\`es \`a pr\'esent que ce th\'eor\`eme ne
donne aucune m\'ethode pour choisir les sous-vari\'et\'es
le long desquelles il faut \'eclater. Ce probl\`eme figure
dans une liste de 100 probl\`emes ouverts en g\'eom\'etrie
\'etablie par S.\ T.\ Yau \cite{Yau93}. Nous donnerons
(modestement) quelques r\'eponses dans cette direction
dans la deuxi\`eme partie de cette th\`ese.

\medskip

Une cons\'equence du th\'eor\`eme pr\'ec\'edent est le
th\'eor\`eme de Chow-Kodaira~: {\em une surface
complexe compacte lisse est de Moishezon si et seulement
si elle est projective}. En effet, on ne peut qu'\'eclater
des points en dimension $2$. Or, de fa\c con g\'en\'erale,
si $X$ est une vari\'et\'e compacte et si $\tilde{X}$
est la vari\'et\'e $X$ \'eclat\'ee au point $x$, alors
$X$ est projective si et seulement si $\tilde{X}$ l'est
(voir par exemple \cite{Kle66}).

Ceci explique le fait que tous les exemples
de vari\'et\'es de Moishezon qui figurent dans notre travail
sont de dimension sup\'erieure ou
\'egale \`a $3$. Le premier exemple de vari\'et\'e de
Moishezon non projective a \'et\'e construit par
H.\ Hironaka (voir par exemple \cite{Har77}). Nous donnons deux constructions
dues \`a I.\ Nakamura et K.\ Oguiso \`a la fin de ces
pr\'eliminaires.

\section{Les in\'egalit\'es de Morse de J.-P.\ Demailly}

Nous rappelons ici les in\'egalit\'es de Morse
holomorphes de J.-P.\ Demailly~: on renvoie \`a \cite{Dem85}
pour la d\'emonstration de ce r\'esultat.

Dans ce qui suit, $X$ d\'esigne une vari\'et\'e
compacte
de dimension $n$ et $(E,h)$ un fibr\'e en droites hermitien sur $X$.
Au couple $(E,h)$,
on associe pour tout entier $q$ compris entre $0$ et $n$
l'ouvert $X(q,E)$
form\'e des points
$x$ de $X$ pour lesquels $\Theta (E)(x)$ poss\`ede
exactement $q$ valeurs propres strictement n\'egatives et $n-q$ valeurs
propres strictement positives ;
finalement on pose $X(\leq q,E) = X(0,E)\cup \cdots \cup X(q,E)$.

Les in\'egalit\'es
de Morse holomorphes donnent une estimation des groupes de cohomologie
\`a valeurs dans les puissances tensorielles $E^{\otimes k}$
en fonction d'int\'egrales de courbure sur $X$.

L'\'enonc\'e pr\'ecis est le suivant~:

\bigskip

\noindent{\bf Th\'eor\`eme (J.-P.\ Demailly, 1985)}
{\em Pour tout $q$ compris entre
$0$ et $n$, et si $F$ est un fibr\'e vectoriel holomorphe
de rang $r$ sur $X$, on a~:

(i) $\displaystyle{ \sum _{j=0}^{q}(-1)^{q-j} \dim H^{j}(X,E^{k}\otimes F)
 \leq r\frac{k^{n}}{n!} \int _{X(\leq q,E)} (-1)^{q} \Theta
(E)^{n} + o(k^{n})}$
(avec \'e\-ga\-li\-t\'e si $q=n$),

(ii) $ \displaystyle{ \dim H^{q}(X,E^{k}\otimes F)
\leq r\frac{k^{n}}{n!} \int _{X(q,E)} (-1)^{q} \Theta (E)^{n}
+ o(k^{n}).} $}

\bigskip

Le point (i) est d\'esign\'e sous le nom plus pr\'ecis
d'in\'egalit\'es de Morse fortes, alors que le point (ii),
qui est une cons\'equence imm\'ediate de (i), est d\'esign\'e
sous le nom
d'in\'egalit\'es de Morse faibles.

Ces in\'egalit\'es ont de nombreuses applications dont certaines
tr\`es r\'ecentes~: elles sont utilis\'ees dans les travaux
de J.-P.\ Demailly et Y.-T.\ Siu en direction de la conjecture de Fujita
\cite{Dem94}, \cite{Siu94}.

En g\'en\'eral, ces in\'egalit\'es peuvent se substituer
aux th\'eor\`emes
d'annulation lorsque la signature de la forme de courbure
n'est pas constante. En particulier, J.-P.\ Demailly
les a utilis\'ees initialement pour renforcer
un r\'esultat de Y.T.\ Siu donnant un crit\`ere analytique
suffisant
pour qu'une vari\'et\'e soit de Moishezon. Ces r\'esultats
fournissent la solution \`a la conjecture de H.\ Grauert
et O.\ Riemenschneider \cite{GrR70}~:

\bigskip

\noindent{\bf Th\'eor\`eme (J.-P.\ Demailly, Y.T.\ Siu, 1985)}
{\em Une vari\'et\'e compacte $X$ est de Moishezon
d\`es que $X$ poss\`ede
un fibr\'e $E$ en droites
muni d'une m\'etrique hermitienne lisse dont la forme de courbure
$\Theta(E)$ v\'erifie l'une des conditions suivantes :

(i) $\Theta(E)$ est partout semi-positive et d\'efinie positive en au moins
un point (``crit\`ere de Siu"),

(ii) $\displaystyle{ \int_{X(\leq 1,E)} \Theta(E)^{n} > 0}$
(``crit\`ere de Demailly").

}

\bigskip

Les deux \'enonc\'es d\'ecoulent des in\'egalit\'es de Morse~:
elles impliquent dans les deux cas que le fibr\'e $E$ est gros.
Avant de commenter ce r\'esultat,
rappelons qu'un fibr\'e en droites
$E$
sur une vari\'et\'e compacte $X$ est dit {\bf num\'eriquement effectif}
(en abr\'eg\'e {\bf nef}) si pour toute m\'etrique hermitienne
$\omega$ sur
$X$ et pour tout $\varepsilon > 0$, le fibr\'e
$E$ poss\`ede une m\'etrique lisse $h_{\varepsilon}$ telle que
$\Theta_{h_{\varepsilon}}(E) \geq -\varepsilon\omega$. Cette notion
a \'et\'e introduite par J.-P.\ Demailly, T.\ Peternell et
M.\ Schneider \cite{DPS94} et admet une formulation
\'equivalente sur les vari\'et\'es projectives~: sur une vari\'et\'e
projective $X$,
un fibr\'e en droites est nef si et seulement si
son intersection avec toute courbe de $X$ est semi-positive.
Sur une vari\'et\'e quelconque, il peut ne pas y avoir
de courbes, cependant il y en a ``suffisamment" sur une vari\'et\'e
de Moishezon et
Mihai Paun a \'etendu le r\'esultat pr\'ec\'edent \`a
ces derni\`eres~: {\em sur une vari\'et\'e de Moishezon,
un fibr\'e en droites est nef si et seulement si
son intersection avec toute courbe est semi-positive} \cite{Pau95}.

\medskip

\noindent {\bf Remarque-exemple }
Un fibr\'e en droites satisfaisant
au crit\`ere de Siu est simultan\'ement gros et nef.
Par ailleurs, des exemples de vari\'et\'es de Moishezon
satisfaisant les crit\`eres (i) et (ii) sont donn\'es
par les vari\'et\'es de Moishezon qui admettent
un morphisme g\'en\'eriquement fini vers une vari\'et\'e
projective.

\section{Exemples explicites}

Nous donnons dans ce paragraphe deux constructions,
l'une utilis\'ee par I.\ Nakamura \cite{Nak87}
et J.\ Koll\'ar \cite{Kol91} et l'autre due \`a K.\ Oguiso.
Nous \'etudions en d\'etail la premi\`ere
et expliquons plus bri\`evement celle de K.\ Oguiso.
Ces constructions nous permettent de montrer
qu'aucun des deux crit\`eres analytiques
pr\'ec\'edents pour qu'une vari\'et\'e soit
de Moishezon n'est n\'ecessaire~: ces \'enonc\'es n'admettent
donc
pas de r\'eciproque dans le cadre des fibr\'es hermitiens
\`a m\'etrique lisse.

\subsection{La premi\`ere construction}

La construction qui suit exhibe une famille de
vari\'et\'es de dimension
$3$ complexe d\'epen\-dant d'un param\`etre entier $m$.
L'origine de cette construction n'est pas tr\`es claire~;
elle est utilis\'ee dans \cite{Nak87} et
mentionn\'ee
dans \cite{Kol91} \S5.
Une des motivations de J.\ Koll\'ar, lorsqu'il mentionne
cet exemple, est de construire une vari\'et\'e
de Moishezon, dont le groupe de Picard est $\ZZ$
et dont le g\'en\'erateur gros du groupe de Picard
est d'auto-intersection n\'egative. En particulier,
ce g\'en\'erateur n'est pas nef.

\subsub{Construction explicite}
La construction est tr\`es simple~: elle
consiste \`a \'eclater $\PP ^3$ le long d'une
courbe contenue dans une quadrique, et \`a
contracter lorsque ceci est possible la transform\'ee
stricte de la quadrique sur une courbe rationnelle lisse.

Soit donc ${\cal Q} \subset \PP ^{3}$ une quadrique lisse,
donn\'ee par exemple par
l'\'equation homog\`ene $xy=zt$, o\`u $\lbrack x:y:z:t \rbrack$ sont les
coordonn\'ees homog\`enes sur $\PP ^{3}$.
La quadrique ${\cal Q}$ est isomorphe \`a $\PP^{1} \times \PP^{1}$
et nous
notons $$L_{1} = \{ \ast \} \times \PP ^{1} \ \mbox{et} \
L_{2} = \PP ^{1} \times \{ \ast \}$$ les g\'en\'erateurs
de $H_2(\cal Q,\ZZ) \simeq \ZZ ^2$.
On a \'evidemment
$$L_{i} \cdot L_{i}=0 \ \mbox{et} \  L_{1} \cdot L_{2}=1.$$
De plus, tout diviseur $D$ de ${\cal Q}$ est num\'eriquement
caract\'eris\'e par un couple d'entiers $(a,b)$ donn\'e par
l'intersection de $D$ avec $L_{1}$ et $L_{2}$~: $$(a,b)
=(D \cdot L_{1},D \cdot L_{2})
\in \ZZ^{2}.$$ Ce couple est appel\'e le type de $D$.
Par exemple, le diviseur canonique $K_{{\cal Q}}$ est de
type $(-2,-2)$.

\medskip

\noindent {\bf Affirmation } {\em  Pour tous $n$ et $m$ entiers positifs,
il existe
une courbe lisse $C_{n,m}$ incluse dans ${\cal Q}$ et
de type $(n,m)$.
Une telle courbe est de genre $g_{n,m}=(n-1)(m-1)$
et de degr\'e $n+m$.}

\medskip

Soit $C_{n,m}$ une telle courbe.
Nous \'eclatons alors $\PP ^{3}$ le long
de $C_{n,m}$~: on obtient une vari\'et\'e
projective $\tilde{X}$, et un morphisme
$$\pi_{1} : \tilde{X} \to \PP ^{3}.$$
Notons $E_{n,m}$ le diviseur exceptionnel de l'\'eclatement~;
il est isomorphe \`a $\PP (N^{\ast}_{C_{n,m}/\PP ^{3}})$.
Comme le groupe de Picard de $\PP ^3$ est $\ZZ$,
celui de $\tilde{X}$ est $\ZZ ^{2}$.
Si $\tilde{{\cal Q}}$ d\'esigne la transform\'ee stricte
de ${\cal Q}$ et
$\tilde{L_{i}}$
celle de $L_{i}$, alors $\tilde{{\cal Q}}$ et $\tilde{L_{i}}$
sont respectivement isomorphes \`a ${\cal Q}$
et $L_{i}$ car $C_{n,m}$ est incluse dans ${\cal Q}$.
De plus, le type du fibr\'e normal
$N_{\tilde{{\cal Q}}/\tilde{X}}$ de $\tilde{{\cal Q}}$
dans $\tilde{X}$ est donn\'e par l'affirmation suivante que
nous d\'emontrons plus loin~:

\medskip

\noindent {\bf Affirmation } {\em On a
$N_{\tilde{{\cal Q}}/\tilde{X}} \cdot
\tilde{L_{1}}
= 2-n$ et $N_{\tilde{{\cal Q}}/\tilde{X}} \cdot \tilde{L_{2}} = 2-m$.
}

\medskip

Comme cas particulier de l'affirmation pr\'ec\'edente,
consid\'erons le
cas o\`u $n=3$. La restriction \`a $\tilde{L_{1}}$ du fibr\'e
$N_{\tilde{{\cal Q}}/\tilde{X}}$ est alors isomorphe au fibr\'e
${\cal O}_{\PP^{1}}(-1)$.
Par le crit\`ere de contraction de Fujiki-Nakano,
il existe donc une vari\'et\'e $X_{m}$ et une application
$$\pi_{2} : \tilde{X} \to X_{m}$$ de sorte que $\pi_{2}$ soit
l'\'eclatement d'une courbe lisse rationnelle $C_{m}$, de fibr\'e normal
projectivement trivial (\'egal \`a ${\cal O}_{\PP ^1}(-m)^{\oplus 2}$)
tel que le diviseur exceptionnel de $\pi_{2}$ est
exactement
$\tilde{{\cal Q}}$.
Evidemment, $X_{m}$ est bim\'eromorphiquement \'equivalente \`a $\PP^{3}$
donc est de Moishezon. De plus, le groupe de Picard
de $X_{m}$ est $\ZZ$.

\bigskip

\vspace{+3mm}
\setlength{\unitlength}{0.0100in}
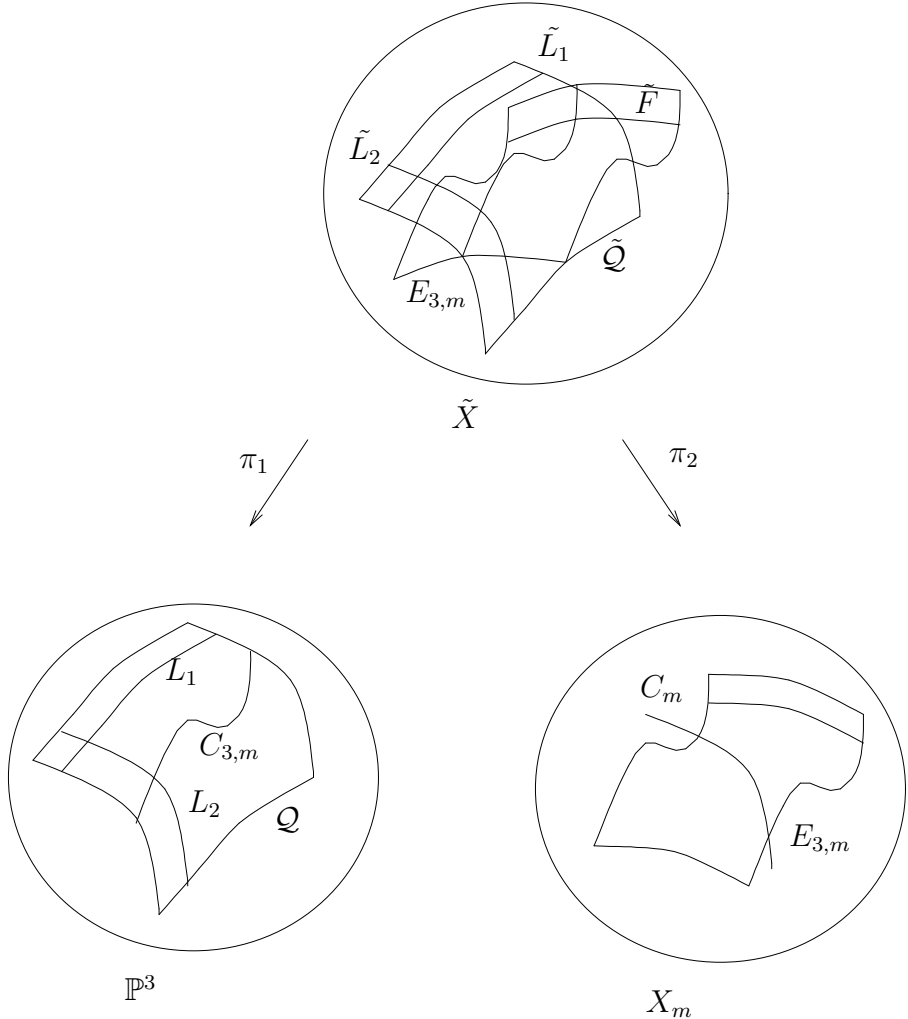
\begin{picture}(470,544)(0,-10)
\put(373,117){\ellipse{194}{182}}
\put(271,429){\ellipse{212}{200}}
\path(157,300)(127,255)
\path(129.774,262.766)(127.000,255.000)(133.102,260.547)
\path(322,300)(352,255)
\path(345.898,260.547)(352.000,255.000)(349.226,262.766)
\path(307,87)	(308.728,91.454)
	(310.345,95.593)
	(311.857,99.431)
	(313.272,102.984)
	(314.595,106.268)
	(315.833,109.298)
	(318.081,114.656)
	(320.069,119.183)
	(321.850,123.000)
	(325.000,129.000)

\path(325,129)	(326.547,131.899)
	(328.570,135.473)
	(334.000,141.000)

\path(334,141)	(338.893,141.000)
	(344.322,138.895)
	(349.840,137.093)
	(355.000,138.000)

\path(355,138)	(360.838,143.296)
	(362.959,146.772)
	(364.596,150.952)
	(365.784,155.951)
	(366.557,161.882)
	(366.799,165.233)
	(366.950,168.861)
	(367.016,172.778)
	(367.000,177.000)

\path(367,177)	(371.390,176.923)
	(375.473,176.836)
	(379.264,176.739)
	(382.779,176.630)
	(386.032,176.508)
	(389.041,176.372)
	(394.382,176.053)
	(398.925,175.666)
	(402.794,175.200)
	(409.000,174.000)

\path(409,174)	(415.135,172.017)
	(418.821,170.528)
	(423.071,168.637)
	(427.999,166.291)
	(433.721,163.438)
	(436.915,161.805)
	(440.350,160.026)
	(444.040,158.093)
	(448.000,156.000)

\path(367,162)	(371.390,161.923)
	(375.473,161.836)
	(379.264,161.739)
	(382.779,161.630)
	(386.032,161.508)
	(389.041,161.372)
	(394.382,161.054)
	(398.925,160.666)
	(402.794,160.200)
	(409.000,159.000)

\path(409,159)	(415.135,157.017)
	(418.821,155.528)
	(423.071,153.637)
	(427.999,151.291)
	(433.721,148.439)
	(436.915,146.805)
	(440.350,145.026)
	(444.040,143.093)
	(448.000,141.000)

\path(307,87)	(311.390,86.923)
	(315.473,86.836)
	(319.264,86.739)
	(322.779,86.630)
	(326.032,86.508)
	(329.041,86.372)
	(334.382,86.054)
	(338.925,85.666)
	(342.794,85.200)
	(349.000,84.000)

\path(349,84)	(355.135,82.017)
	(358.821,80.528)
	(363.071,78.637)
	(367.999,76.291)
	(373.721,73.439)
	(376.915,71.805)
	(380.350,70.026)
	(384.040,68.093)
	(388.000,66.000)

\path(388,66)	(389.728,70.454)
	(391.345,74.593)
	(392.857,78.431)
	(394.272,81.984)
	(395.595,85.268)
	(396.833,88.298)
	(399.081,93.656)
	(401.069,98.183)
	(402.850,102.000)
	(406.000,108.000)

\path(406,108)	(407.547,110.899)
	(409.570,114.473)
	(415.000,120.000)

\path(415,120)	(419.893,120.000)
	(425.322,117.895)
	(430.840,116.093)
	(436.000,117.000)

\path(436,117)	(441.838,122.296)
	(443.959,125.772)
	(445.596,129.952)
	(446.784,134.951)
	(447.557,140.882)
	(447.799,144.233)
	(447.950,147.861)
	(448.016,151.778)
	(448.000,156.000)

\path(184,426)	(189.926,423.721)
	(195.420,421.547)
	(200.501,419.467)
	(205.189,417.471)
	(209.504,415.546)
	(213.466,413.683)
	(217.094,411.870)
	(220.409,410.096)
	(226.175,406.622)
	(230.924,403.174)
	(234.813,399.662)
	(238.000,396.000)

\path(238,396)	(240.821,391.325)
	(243.149,385.576)
	(244.147,382.408)
	(245.043,379.101)
	(245.845,375.700)
	(246.560,372.248)
	(247.195,368.788)
	(247.757,365.364)
	(248.253,362.019)
	(248.691,358.796)
	(249.419,352.892)
	(250.000,348.000)

\path(250,348)	(250.000,345.000)

\path(250,345)	(254.247,349.854)
	(257.908,354.037)
	(261.043,357.620)
	(263.715,360.674)
	(267.915,365.474)
	(271.000,369.000)

\path(271,369)	(275.498,374.379)
	(278.244,377.726)
	(281.157,381.243)
	(284.112,384.730)
	(286.987,387.986)
	(292.000,393.000)

\path(292,393)	(297.433,397.137)
	(300.943,399.495)
	(305.134,402.141)
	(310.121,405.144)
	(316.018,408.574)
	(319.344,410.472)
	(322.940,412.503)
	(326.821,414.676)
	(331.000,417.000)

\path(199,420)	(203.247,424.854)
	(206.908,429.037)
	(210.043,432.620)
	(212.715,435.674)
	(216.915,440.474)
	(220.000,444.000)

\path(220,444)	(224.498,449.379)
	(227.244,452.726)
	(230.157,456.243)
	(233.112,459.730)
	(235.987,462.986)
	(241.000,468.000)

\path(241,468)	(246.433,472.137)
	(249.943,474.495)
	(254.134,477.141)
	(259.121,480.144)
	(265.018,483.574)
	(268.344,485.472)
	(271.940,487.503)
	(275.821,489.676)
	(280.000,492.000)

\path(184,426)	(188.247,430.854)
	(191.908,435.037)
	(195.043,438.620)
	(197.715,441.674)
	(201.915,446.474)
	(205.000,450.000)

\path(205,450)	(209.498,455.379)
	(212.244,458.726)
	(215.157,462.243)
	(218.112,465.730)
	(220.987,468.986)
	(226.000,474.000)

\path(226,474)	(231.433,478.137)
	(234.943,480.495)
	(239.134,483.141)
	(244.121,486.144)
	(250.018,489.574)
	(253.344,491.472)
	(256.940,493.503)
	(260.821,495.676)
	(265.000,498.000)

\path(292,393)	(293.728,397.454)
	(295.345,401.593)
	(296.857,405.431)
	(298.272,408.984)
	(299.595,412.268)
	(300.833,415.298)
	(303.081,420.656)
	(305.069,425.183)
	(306.850,429.000)
	(310.000,435.000)

\path(310,435)	(311.547,437.899)
	(313.570,441.473)
	(319.000,447.000)

\path(319,447)	(323.893,447.000)
	(329.322,444.895)
	(334.840,443.093)
	(340.000,444.000)

\path(340,444)	(345.838,449.296)
	(347.959,452.772)
	(349.596,456.952)
	(350.784,461.951)
	(351.557,467.882)
	(351.799,471.233)
	(351.950,474.861)
	(352.016,478.778)
	(352.000,483.000)

\path(202,384)	(205.676,385.460)
	(209.100,386.805)
	(215.246,389.167)
	(220.544,391.121)
	(225.099,392.702)
	(229.017,393.945)
	(232.403,394.885)
	(238.000,396.000)

\path(238,396)	(241.793,396.383)
	(246.115,396.567)
	(251.123,396.543)
	(256.976,396.303)
	(260.268,396.099)
	(263.832,395.837)
	(267.685,395.517)
	(271.849,395.137)
	(276.343,394.697)
	(281.186,394.195)
	(286.398,393.629)
	(292.000,393.000)

\path(13,132)	(18.926,129.721)
	(24.420,127.547)
	(29.501,125.467)
	(34.189,123.471)
	(38.504,121.546)
	(42.466,119.683)
	(46.094,117.870)
	(49.409,116.096)
	(55.175,112.622)
	(59.924,109.174)
	(63.813,105.662)
	(67.000,102.000)

\path(67,102)	(69.821,97.325)
	(72.149,91.576)
	(73.147,88.408)
	(74.043,85.101)
	(74.845,81.700)
	(75.560,78.248)
	(76.195,74.788)
	(76.757,71.364)
	(77.253,68.019)
	(77.691,64.796)
	(78.419,58.892)
	(79.000,54.000)

\path(79,54)	(79.000,51.000)

\path(79,51)	(83.247,55.854)
	(86.908,60.037)
	(90.043,63.620)
	(92.715,66.674)
	(96.915,71.474)
	(100.000,75.000)

\path(100,75)	(104.498,80.379)
	(107.244,83.726)
	(110.157,87.243)
	(113.112,90.730)
	(115.987,93.986)
	(121.000,99.000)

\path(121,99)	(126.433,103.137)
	(129.943,105.495)
	(134.134,108.141)
	(139.121,111.144)
	(145.018,114.574)
	(148.344,116.472)
	(151.940,118.503)
	(155.821,120.676)
	(160.000,123.000)

\path(94,204)	(99.926,201.721)
	(105.420,199.547)
	(110.501,197.467)
	(115.189,195.471)
	(119.504,193.546)
	(123.466,191.683)
	(127.094,189.870)
	(130.409,188.096)
	(136.175,184.622)
	(140.924,181.174)
	(144.813,177.662)
	(148.000,174.000)

\path(148,174)	(150.821,169.325)
	(153.149,163.576)
	(154.147,160.408)
	(155.043,157.101)
	(155.845,153.700)
	(156.560,150.248)
	(157.195,146.788)
	(157.757,143.364)
	(158.253,140.019)
	(158.691,136.796)
	(159.419,130.892)
	(160.000,126.000)

\path(160,126)	(160.000,123.000)

\path(28,126)	(32.247,130.854)
	(35.908,135.037)
	(39.043,138.620)
	(41.715,141.674)
	(45.915,146.474)
	(49.000,150.000)

\path(49,150)	(53.498,155.379)
	(56.244,158.726)
	(59.157,162.243)
	(62.112,165.730)
	(64.987,168.986)
	(70.000,174.000)

\path(70,174)	(75.433,178.137)
	(78.943,180.495)
	(83.134,183.141)
	(88.121,186.144)
	(94.018,189.575)
	(97.344,191.472)
	(100.940,193.503)
	(104.821,195.676)
	(109.000,198.000)

\path(13,132)	(17.247,136.854)
	(20.908,141.037)
	(24.043,144.620)
	(26.715,147.674)
	(30.915,152.474)
	(34.000,156.000)

\path(34,156)	(38.498,161.379)
	(41.244,164.726)
	(44.157,168.243)
	(47.112,171.730)
	(49.987,174.986)
	(55.000,180.000)

\path(55,180)	(60.433,184.137)
	(63.943,186.495)
	(68.134,189.141)
	(73.121,192.144)
	(79.018,195.575)
	(82.344,197.472)
	(85.940,199.503)
	(89.821,201.676)
	(94.000,204.000)

\path(67,99)	(68.728,103.454)
	(70.345,107.593)
	(71.857,111.431)
	(73.272,114.984)
	(74.595,118.268)
	(75.833,121.298)
	(78.081,126.656)
	(80.069,131.183)
	(81.850,135.000)
	(85.000,141.000)

\path(85,141)	(86.547,143.899)
	(88.570,147.473)
	(94.000,153.000)

\path(94,153)	(98.893,153.000)
	(104.322,150.895)
	(109.840,149.093)
	(115.000,150.000)

\path(115,150)	(120.838,155.296)
	(122.959,158.772)
	(124.596,162.952)
	(125.784,167.951)
	(126.557,173.882)
	(126.799,177.233)
	(126.950,180.861)
	(127.016,184.778)
	(127.000,189.000)

\path(28,147)	(33.926,144.721)
	(39.420,142.547)
	(44.501,140.467)
	(49.189,138.471)
	(53.504,136.546)
	(57.466,134.683)
	(61.094,132.870)
	(64.409,131.096)
	(70.175,127.622)
	(74.924,124.174)
	(78.813,120.662)
	(82.000,117.000)

\path(82,117)	(84.821,112.325)
	(87.149,106.576)
	(88.147,103.408)
	(89.043,100.101)
	(89.845,96.700)
	(90.560,93.248)
	(91.195,89.788)
	(91.757,86.364)
	(92.253,83.019)
	(92.691,79.796)
	(93.419,73.892)
	(94.000,69.000)

\path(94,69)	(94.000,66.000)

\put(97,123){\ellipse{194}{182}}
\path(334,156)	(339.926,153.721)
	(345.420,151.547)
	(350.501,149.467)
	(355.189,147.471)
	(359.504,145.546)
	(363.466,143.683)
	(367.094,141.870)
	(370.409,140.096)
	(376.175,136.622)
	(380.924,133.174)
	(384.813,129.662)
	(388.000,126.000)

\path(388,126)	(390.821,121.325)
	(393.149,115.576)
	(394.147,112.408)
	(395.043,109.101)
	(395.845,105.700)
	(396.560,102.248)
	(397.195,98.788)
	(397.757,95.364)
	(398.253,92.019)
	(398.691,88.796)
	(399.419,82.892)
	(400.000,78.000)

\path(400,78)	(400.000,75.000)

\put(61,6){\makebox(0,0)[lb]{\raisebox{0pt}[0pt][0pt]{\shortstack[l]{{\egtrm
$\PP^{3}$}}}}}
\path(199,444)	(204.926,441.721)
	(210.420,439.547)
	(215.501,437.467)
	(220.189,435.471)
	(224.504,433.546)
	(228.466,431.683)
	(232.094,429.870)
	(235.409,428.096)
	(241.175,424.622)
	(245.924,421.174)
	(249.813,417.662)
	(253.000,414.000)

\path(253,414)	(255.821,409.325)
	(258.149,403.576)
	(259.147,400.408)
	(260.043,397.101)
	(260.845,393.700)
	(261.560,390.248)
	(262.195,386.788)
	(262.757,383.364)
	(263.253,380.019)
	(263.691,376.796)
	(264.419,370.892)
	(265.000,366.000)

\path(265,366)	(265.000,363.000)

\path(265,498)	(270.926,495.721)
	(276.420,493.547)
	(281.501,491.467)
	(286.189,489.471)
	(290.504,487.546)
	(294.466,485.683)
	(298.094,483.870)
	(301.409,482.096)
	(307.175,478.622)
	(311.924,475.174)
	(315.813,471.662)
	(319.000,468.000)

\path(319,468)	(321.821,463.325)
	(324.149,457.576)
	(325.147,454.408)
	(326.043,451.101)
	(326.845,447.700)
	(327.560,444.248)
	(328.195,440.788)
	(328.757,437.364)
	(329.253,434.019)
	(329.691,430.796)
	(330.419,424.892)
	(331.000,420.000)

\path(331,420)	(331.000,417.000)

\path(202,384)	(203.728,388.454)
	(205.345,392.593)
	(206.857,396.431)
	(208.272,399.984)
	(209.595,403.268)
	(210.833,406.298)
	(213.081,411.656)
	(215.069,416.183)
	(216.850,420.000)
	(220.000,426.000)

\path(220,426)	(221.547,428.899)
	(223.570,432.473)
	(229.000,438.000)

\path(229,438)	(233.893,438.000)
	(239.322,435.895)
	(244.840,434.093)
	(250.000,435.000)

\path(250,435)	(255.838,440.296)
	(257.959,443.772)
	(259.596,447.952)
	(260.784,452.951)
	(261.557,458.882)
	(261.799,462.233)
	(261.950,465.861)
	(262.016,469.778)
	(262.000,474.000)

\path(262,474)	(265.676,475.460)
	(269.100,476.805)
	(275.246,479.167)
	(280.544,481.121)
	(285.099,482.702)
	(289.017,483.945)
	(292.403,484.885)
	(298.000,486.000)

\path(298,486)	(301.793,486.383)
	(306.115,486.567)
	(311.123,486.543)
	(316.976,486.303)
	(320.268,486.099)
	(323.832,485.837)
	(327.685,485.517)
	(331.849,485.137)
	(336.343,484.697)
	(341.186,484.195)
	(346.398,483.629)
	(352.000,483.000)

\path(262,456)	(265.676,457.460)
	(269.100,458.805)
	(275.246,461.167)
	(280.544,463.121)
	(285.099,464.702)
	(289.017,465.945)
	(292.403,466.885)
	(298.000,468.000)

\path(298,468)	(301.793,468.383)
	(306.115,468.567)
	(311.123,468.543)
	(316.976,468.303)
	(320.268,468.099)
	(323.832,467.837)
	(327.685,467.517)
	(331.849,467.137)
	(336.343,466.697)
	(341.186,466.195)
	(346.398,465.629)
	(352.000,465.000)

\path(238,396)	(239.728,400.454)
	(241.345,404.593)
	(242.857,408.431)
	(244.272,411.984)
	(245.595,415.268)
	(246.833,418.298)
	(249.081,423.656)
	(251.069,428.183)
	(252.850,432.000)
	(256.000,438.000)

\path(256,438)	(257.547,440.899)
	(259.570,444.473)
	(265.000,450.000)

\path(265,450)	(269.893,450.000)
	(275.322,447.895)
	(280.840,446.093)
	(286.000,447.000)

\path(286,447)	(291.838,452.296)
	(293.959,455.772)
	(295.596,459.952)
	(296.784,464.951)
	(297.557,470.882)
	(297.799,474.233)
	(297.950,477.861)
	(298.016,481.778)
	(298.000,486.000)

\put(331,165){\makebox(0,0)[lb]{\raisebox{0pt}[0pt][0pt]{\shortstack[l]{{\tenrm
$C_{m}$}}}}}
\put(334,0){\makebox(0,0)[lb]{\raisebox{0pt}[0pt][0pt]{\shortstack[l]{{\tenrm
$X_{m}$}}}}}
\put(121,285){\makebox(0,0)[lb]{\raisebox{0pt}[0pt][0pt]{\shortstack[l]{{\tenrm
$\pi_{1}$}}}}}
\put(346,288){\makebox(0,0)[lb]{\raisebox{0pt}[0pt][0pt]{\shortstack[l]{{\tenrm
$\pi_{2}$}}}}}
\put(82,174){\makebox(0,0)[lb]{\raisebox{0pt}[0pt][0pt]{\shortstack[l]{{\tenrm
$L_{1}$}}}}}
\put(100,135){\makebox(0,0)[lb]{\raisebox{0pt}[0pt][0pt]{\shortstack[l]{{\tenrm
$C_{3,m}$}}}}}
\put(139,96){\makebox(0,0)[lb]{\raisebox{0pt}[0pt][0pt]{\shortstack[l]{{\tenrm
${\cal Q}$}}}}}
\put(94,105){\makebox(0,0)[lb]{\raisebox{0pt}[0pt][0pt]{\shortstack[l]{{\tenrm
$L_{2}$}}}}}
\put(232,306){\makebox(0,0)[lb]{\raisebox{0pt}[0pt][0pt]{\shortstack[l]{{\tenrm
$\tilde{X}$}}}}}
\put(310,390){\makebox(0,0)[lb]{\raisebox{0pt}[0pt][0pt]{\shortstack[l]{{\tenrm
$\tilde{{\cal  Q}}$}}}}}
\put(328,471){\makebox(0,0)[lb]{\raisebox{0pt}[0pt][0pt]{\shortstack[l]{{\tenrm
$\tilde{F}$}}}}}
\put(277,501){\makebox(0,0)[lb]{\raisebox{0pt}[0pt][0pt]{\shortstack[l]{{\tenrm
$\tilde{L_{1}}$}}}}}
\put(178,447){\makebox(0,0)[lb]{\raisebox{0pt}[0pt][0pt]{\shortstack[l]{{\tenrm
$\tilde{L_{2}}$}}}}}
\put(208,372){\makebox(0,0)[lb]{\raisebox{0pt}[0pt][0pt]{\shortstack[l]{{\tenrm
$E_{3,m}$}}}}}
\put(409,87){\makebox(0,0)[lb]{\raisebox{0pt}[0pt][0pt]{\shortstack[l]{{\tenrm
$E_{3,m}$}}}}}
\end{picture}


\medskip

Avant d'\'etudier plus en d\'etail la vari\'et\'e
$X_m$, d\'emontrons les deux affirmations n\'ecessaires
\`a sa construction.

\medskip

\noindent {\bf D\'emonstration des affirmations}

L'existence de $C_{n,m}$ r\'esulte du fait que
$\displaystyle{{\cal O}(n,m)=
\pr_{1}^{\ast}{\cal O}(n)\otimes \pr_{2}^{\ast}{\cal O}(m)}$
est tr\`es ample sur
$\PP ^{1} \times \PP ^{1}$. Enfin, le calcul du genre
est donn\'e par la formule classique
$$\displaystyle{2g_{n,m}-2 = C_{n,m} \cdot (C_{n,m}+K_{{\cal Q}})}.$$
Ici, $$C_{n,m} \cdot C_{n,m}=2nm
\ \mbox{et} \  C_{n,m} \cdot K_{{\cal Q}}= -2(n+m).$$ Ceci d\'emontre
la premi\`ere affirmation.\finpreuve

\medskip

Pour la deuxi\`eme, la suite exacte
$$ 0 \rightarrow T \tilde{{\cal Q}} \to T\tilde{X}_{| \tilde{{\cal Q}}}
\to
N_{\tilde{{\cal Q}}/\tilde{X}} \to 0$$
donne $N_{\tilde{{\cal Q}}/\tilde{X}}=K_{\tilde{{\cal Q}}}-
K_{\tilde{X}|\tilde{{\cal Q}}}$ o\`u $K_{\tilde{X}}=\pi_{1}^{\ast}K_{\PP^{3}}+
{\cal O}(E_{n,m})$.
De l\`a~:

\vspace{-3mm}
\begin{eqnarray*}
N_{\tilde{{\cal Q}}/\tilde{X}} \cdot \tilde{L_{i}} & = &
K_{\tilde{{\cal Q}}} \cdot \tilde{L_{i}}-
\pi_{1}^{\ast}K_{\PP^{3}} \cdot \tilde{L_{i}}-
{\cal O}(E_{n,m}) \cdot \tilde{L_{i}}\\
&  = & K_{{\cal Q}} \cdot L_{i}-
K_{\PP^{3}} \cdot L_{i}-C_{n,m} \cdot L_{i}.
\end{eqnarray*}

Or, $K_{{\cal Q}}={\cal O}_{\PP^{3}}(-2)_{|{\cal Q}}$ et $K_{\PP^{3}}=
{\cal O}_{\PP^{3}}(-4)$. Ceci conclut le calcul.\finpreuve
\subsub{Deux propri\'et\'es de $X_m$}

Nous montrons ici que les crit\`eres de Demailly et Siu
ne sont pas satisfaits pour la vari\'et\'e de Moishezon $X_m$.

\bigskip

\noindent {\bf Th\'eor\`eme A}
{\em La vari\'et\'e $X_m$ v\'erifie les deux propri\'et\'es
suivantes~:

(i) si $m$ est strictement plus grand que $3$,
$X_{m}$ ne poss\`ede pas de fibr\'e en droites \`a la fois gros et nef,
et donc ne satisfait pas au crit\`ere de Siu,

(ii) si $m$ est strictement plus grand que $5$, $X_{m}$
ne poss\`ede pas de fibr\'e en droites $E$
muni d'une m\'etrique
hermitienne lisse $h$ telle
que la forme de courbure $\Theta (E)$ v\'erifie :
$$ \int_{X(\leq 1,E)}\Theta (E)^{3} > 0.$$
}

\medskip

L'affirmation (i) est due \`a J.\ Koll\'ar, nous en donnons
une preuve \'el\'ementaire, tandis que (ii) est nouveau \`a notre
connaissance.

\medskip

\noindent {\bf D\'emonstration de (i)}

Soit $E$ un fibr\'e holomorphe de rang $1$
sur $X_{m}$, que l'on suppose non trivial (le fibr\'e trivial, bien que nef,
n'est pas gros ! ).
Il existe alors des entiers $k$ et $l$ tels que~:
$$\pi_{2}^{\ast}E = \pi_{1}^{\ast}{\cal O}_{\PP^{3}}(l)-{\cal O}(kE_{3,m}).$$
\noindent Comme $\tilde{L_{1}}$ est une fibre de $\pi_{2}$, on a~:
$\pi_{2}^{*}E \cdot \tilde{L_{1}} = 0$. On en d\'eduit la relation $l=3k$
et donc
$$\pi_{2}^{*}E = k(3\pi_{1}^{\ast}{\cal O}_{\PP^{3}}(1)-{\cal O}(E_{3,m})),$$
o\`u $k$ est un entier non nul.
En particulier, si $\tilde{F}$ est une fibre non triviale de $\pi_{1}$
dans $\tilde{X}$,
on a les
nombres d'intersection suivants~:

$$\left\{
\begin{array}{l}
\pi_{2}^{\ast}E \cdot \tilde{L_{2}}=k(3-m) \\
\pi_{2}^{\ast}E \cdot \tilde{F}=k
\end{array}
\right.
$$

\noindent On en d\'eduit que pour $m>3$, le fibr\'e
$\pi_{2}^{\ast}E$ n'est pas nef (sinon son
intersection avec toute courbe serait positive ou nulle), et donc
$E$ n'est pas nef.\finpreuve

\noindent {\bf D\'emonstration de (ii)}

Notons dans la suite ${\cal O}_{X_m}(1)$ le g\'en\'erateur du
groupe de Picard de
$X_{m}$
tel que $\pi_{2}^{\ast}{\cal O}_{X_m}(1) =
\pi_{1}^{\ast}{\cal O}_{\PP^{3}}(3)-{\cal O}(E_{3,m})$.

\bigskip

\noindent {\bf Affirmation } {\em Le fibr\'e canonique
$K_{X_{m}}$ est \'egal \`a ${\cal O}_{X_m}(-2)$.}

\medskip

En effet, on a
$$K_{\tilde{X}}=
\pi_{2}^{\ast}K_{X_{m}}+{\cal O}(\tilde{{\cal Q}}) =
\pi_1^*{\cal O}_{\PP^{3}}(-4)+ {\cal O}(E_{3,m})$$
\noindent par
construction. Or,
$$ {\cal O}(\tilde{{\cal Q}}) = \pi_1^*{\cal O}_{\PP^{3}}(2) - {\cal
O}(E_{3,m}),$$
d'o\`u l'affirmation.

\bigskip

\noindent {\bf Affirmation } {\em Les espaces de sections
holomorphes
$H^{0}(X_m,{\cal O}_{X_m}(k))$ sont nuls pour tout entier
$k < 0$ et le fibr\'e ${\cal O}_{X_m}(1)$ est gros.}

\medskip

Par invariance bim\'eromorphe des plurigenres, les groupes
$H^{0}(X_m,{\cal O}_{X_m}(k))$ sont nuls pour tout entier
$k$ strictement n\'egatif et pair. Comme
$X_{m}$ est de Moishezon de groupe de Picard $\ZZ$, $X_{m}$ poss\`ede
un fibr\'e gros qui n'est donc pas ${\cal O}_{X_m}(-1)$, c'est donc que
${\cal O}_{X_m}(1)$
est gros et que les $H^{0}(X_m,{\cal O}_{X_m}(k))$ sont nuls pour tout entier
$k$ strictement n\'egatif. Ceci d\'emontre l'affirmation.

\bigskip

Par dualit\'e de Serre, on d\'eduit des deux affirmations
pr\'ec\'edentes que {\em les groupes de cohomologie
$H^{3}(X_{m},{\cal O}_{X_m}(k))$ sont nuls pour tout
entier
$k > -2$.}

\medskip

Nous sommes maintenant en mesure de passer \`a la preuve
de (ii) proprement dit~:
raisonnons par l'absurde et supposons que ${\cal O}_{X_m}(1)$ poss\`ede
une telle
m\'etrique.
D'apr\`es les in\'egalit\'es de Morse holomorphes de J.-P.\ Demailly
appliqu\'ees pour $q=1$ au fibr\'e hermitien $({\cal O}_{X_m}(1),h)$, on a~:
$$\dim H^{0}(X_{m},{\cal O}_{X_m}(k))-
\dim H^{1}(X_{m},{\cal O}_{X_m}(k)) \geq
\frac{1}{6} \left(\int_{X(\leq 1,E)}\Theta (E)^{3}\right)k^{3}+
o(k^{3}).$$
\noindent Comme les groupes de cohomologie
$H^{3}(X_{m},{\cal O}_{X_m}(k))$ sont nuls pour tout
entier
$k > -2$, on a successivement~:

\vspace{-3mm}
\begin{eqnarray*}
c_{1}({\cal O}_{X_m}(1))^{3}\frac{k^{3}}{6} + o(k^{3})
& = &
\sum _{i=0}^3 (-1)^i \dim H^{i}(X_{m},{\cal O}_{X_m}(k)) \\
& = & \sum _{i=0}^2 (-1)^i \dim H^{i}(X_{m},{\cal O}_{X_m}(k)) \\
& \geq & \dim H^{0}(X_{m},{\cal O}_{X_m}(k))-
\dim H^{1}(X_{m},{\cal O}_{X_m}(k)) \\
& \geq &  \left(\int_{X(\leq 1,E)}\Theta (E)^{3}\right)\frac{k^{3}}{6}+
o(k^{3}).
\end{eqnarray*}

\noindent On en d\'eduit que
$\displaystyle{c_{1}({\cal O}_{X_m}(1))^{3} \geq
\int_{X(\leq 1,E)}\Theta (E)^{3} >0}$.
Il suffit donc de montrer pour obtenir la contradiction cherch\'ee
que pour $m > 5$, on a $c_{1}({\cal O}_{X_m}(1))^{3} \leq 0$.
Or, cette derni\`ere quantit\'e est ais\'ement calculable~:

\medskip

\noindent {\bf Affirmation } {\em
La quantit\'e $c_{1}({\cal O}_{X_m}(1))^{3}$ est \'egale \`a $6-m$.}

\medskip

En effet
\vspace{-3mm}
\begin{eqnarray*}
c_{1}({\cal O}_{X_m}(1))^{3} & = &
c_{1}(\pi_{1}^{*}{\cal O}_{\PP^{3}}(3)-{\cal O}(E_{3,m}))^{3} \\
& = &
c_{1}({\cal O}_{\PP^{3}}(3))^{3}
-3\int_{E_{3,m}} c_{1}(\pi_{1}^{*}{\cal O}_{\PP^{3}}(3))^{2} \\
& & +3\int_{\tilde{X}}c_{1}(\pi_{1}^{*}{\cal O}_{\PP^{3}}(3))
\wedge c_{1}({\cal O}(E_{3,m}))^{2}
-E_{3,m}^{3}.
\end{eqnarray*}

\noindent On a \'evidemment $c_{1}({\cal O}_{\PP^{3}}(3))^{3}=27$ et
$\displaystyle{\int_{E_{3,m}} c_{1}(\pi_{1}^{*}{\cal O}_{\PP^{3}}(3))^{2}= 0.}$

\noindent Pour les deux termes restants, on commence par remarquer que
$\displaystyle{c_{1}({\cal O}(E_{3,m}))_{|E_{3,m}}}$ est \'egale \`a $-\xi$
o\`u
$\displaystyle{\xi=c_{1}({\cal O}_{\PP (N_{C_{3,m}/\PP ^{3}}^{*})}(1))}$
d\'esigne la
classe fondamentale de l'\'eclatement.
On en d\'eduit que~:
\vspace{-3mm}
\begin{eqnarray*}
\int_{\tilde{X}}c_{1}(\pi_{1}^{*}{\cal O}_{\PP^{3}}(3))
\wedge c_{1}({\cal O}(E_{3,m}))^{2} & = &
-\int_{E_{3,m}} \pi_{1}^{*}c_{1}({\cal O}_{\PP^{3}}(3))
\wedge \xi \\
& = & -\int_{C_{3,m}}c_{1}({\cal O}_{\PP^{3}}(3)) \\
& = & -3(3+m),
\end{eqnarray*}

\noindent la derni\`ere \'egalit\'e venant du fait
que $C_{3,m}$ est de degr\'e
$3+m$ dans $\PP ^{3}$ (rappelons que
$\displaystyle{{\cal O}_{\PP ^{3}}(1)_{|Q}=
{\cal O}(1,1)}$).
Finalement, il reste \`a calculer $E_{3,m}^{3}$. Pour cela, rappelons
que $\xi$ v\'erifie la formule fondamentale suivante~:
$$ \xi^{2} - \pi_{1}^{*}c_{1}(N_{C_{3,m}/\PP^{3}}^{*})\xi +
\pi_{1}^{*}c_{2}(N_{C_{3,m}/\PP^{3}}^{*})=0 $$ qui se r\'eduit ici \`a
$\displaystyle{\xi^{2} - \pi_{1}^{*}c_{1}(N_{C_{3,m}/\PP^{3}}^{*})\xi =0}$.

\noindent Il vient alors
$\displaystyle{E_{3,m}^{3}= \int_{E_{3,m}}\xi^{2}=\int_{E_{3,m}}
\pi_{1}^{*}c_{1}(N_{C_{3,m}/\PP^{3}}^{*})\xi =
\int_{C_{3,m}}c_{1}(N_{C_{3,m}/\PP^{3}}^{*}).}$

\noindent Or la suite exacte~:
$$ 0 \rightarrow TC_{3,m} \rightarrow T\PP^{3}_{| C_{3,m}} \rightarrow
N_{C_{3,m}/\PP^{3}} \rightarrow 0$$
donne de suite :
$$ \int_{C_{3,m}}c_{1}(N_{C_{3,m}/\PP^{3}}^{*})=
\int_{C_{3,m}}c_{1}({\cal O}_{\PP^{3}}(-4))-2g_{3,m}+2
=-6-8m .$$

\noindent Il reste finalement~:
$$c_{1}({\cal O}_{X_m}(1))^{3} = 27 - 27 - 9m + 6 + 8m = 6-m.$$
\noindent Ceci ach\`eve la preuve de l'affirmation et
par suite celle du th\'eor\`eme.\finpreuve

\subsection{La construction de K.\ Oguiso}

Pour ce paragraphe, la r\'ef\'erence est \cite{Ogu94}.
Dans cet article, K.\ Oguiso construit une vari\'et\'e
de Moishezon non projective, de dimension $3$ et qui est de plus
de Calabi-Yau~: ceci signifie que cette vari\'et\'e
est simplement connexe et \`a fibr\'e canonique trivial.
Mettons en garde le lecteur sur l'usage fait ici de l'expression
Calabi-Yau. En effet, la vari\'et\'e en question
n'est pas k\"ahl\'erienne. Rappelons plus g\'en\'eralement
qu'un th\'eor\`eme de B.\ Moishezon  affirme que {\em toute vari\'et\'e
de Moishezon k\"ahl\'erienne est projective} (voir
le ``survey" de T.\ Peternell \cite{Pet95} pour une d\'emonstration
rapide de ce r\'esultat).

Le r\'esultat de K.\ Oguiso est le suivant~:

\bigskip

\noindent {\bf Th\'eor\`eme (K.\ Oguiso, 1994)}
{\em Il existe une vari\'et\'e de Moishezon $Y$, de dimension $3$
et de Calabi-Yau telle que
$\displaystyle{H^2(Y,\ZZ) = \Pic (Y) = \ZZ \cdot L}$
o\`u $L$ satisfait $L^3 := c_1(L)^3 =0$. Autrement dit,
la forme cubique d'intersection sur $H^2(Y,\ZZ)$
est identiquement nulle.}

\bigskip

Nous donnons dans la derni\`ere partie de cette th\`ese
une construction, diff\'erente de celle d'Oguiso, permettant
de retrouver ce r\'esultat. Dans \cite{Ogu94},
K.\ Oguiso obtient ce th\'eor\`eme
comme cons\'equence du r\'esultat suivant~:

\bigskip

\noindent {\bf Th\'eor\`eme (K.\ Oguiso, 1994)}
{\em Soit $d$ un entier positif. Alors il existe
une vari\'et\'e projective $X_d$ de dimension $3$,
intersection compl\`ete d'une quadrique et d'une
quartique dans $\PP ^5$ et contenant une courbe rationnelle
lisse $C_d$
de degr\'e $d$
dont le fibr\'e normal $N_{C_d/X_d}$ est ${\cal O}_{C_d}(-1)^{\oplus 2}$.}

\bigskip

Montrons que ce dernier r\'esultat implique le
premier.

\noindent Soit donc $X_d$ comme ci dessus et soit
$\tilde{X}_d$ la vari\'et\'e projective obtenue
en \'eclatant $X_d$ le long de la courbe $C_d$. Un argument
identique \`a celui d\'evelopp\'e dans la construction
pr\'ec\'edente assure que l'on peut contracter
le diviseur exceptionnel $E \simeq C_d \times \PP ^1
= \PP ^1 \times \PP ^1$ dans l'autre direction~:
on note $Y_d$ la vari\'et\'e obtenue. Par construction,
$Y_d$ est de Moishezon, de Calabi-Yau et son groupe de
Picard est $\ZZ$.
Enfin, si $\displaystyle{{\cal O}_{Y_d}(1)}$ est le g\'en\'erateur
gros de $\displaystyle{\Pic (Y_d)}$, on montre comme pr\'ec\'edemment
que $\displaystyle{c_1({\cal O}_{Y_d}(1))^3 = 8 - d^3}$.
Le premier r\'esultat est d\'emontr\'e
en choisissant $d=2$ i.e $Y := Y_2$.

\noindent Comme pour les vari\'et\'es $X_m$
pr\'ec\'edemment construites, {\em les vari\'et\'es $Y_d$ pour
$d$ sup\'erieur ou \'egal \`a $2$ ne satisfont pas
les crit\`eres de Siu et Demailly.}

\medskip

\noindent {\bf Remarque } Les vari\'et\'es $Y_d$ ont \'et\'e
obtenues apr\`es
une transformation birationnelle classique en th\'eorie
de Mori appel\'ee ``flop". Nous revenons sur ce type de
construction dans la deuxi\`eme partie de cette th\`ese.

\subsection{Quelques commentaires}

Ces pr\'eliminaires illustrent une partie
des motivations de cette th\`ese. En effet, alors que
les deux constructions pr\'ec\'edentes montrent
que les crit\`eres de J.-P.\ Demailly et Y.-T.\ Siu ne sont
pas n\'ecessaires dans le cadre des fibr\'es hermitiens
munis d'une m\'etrique ${\cal C}^{\infty}$, nous montrons
dans le deuxi\`eme chapitre de cette th\`ese que ces crit\`eres
deviennent n\'ecessaires et suffisants dans le cadre
plus souple des m\'etriques singuli\`eres. Pour cela, nous
\'etendons les in\'egalit\'es de Morse en autorisant
un certain type de singularit\'es aux m\'etriques des fibr\'es
consid\'er\'es.

Remarquons ensuite que les constructions pr\'ec\'edentes
donnent	des exem\-ples de vari\'et\'es de Moishezon de dimension $3$,
de groupe de Picard $\ZZ$ avec respectivement $-K_X$ gros et
$K_X$ trivial. Une des motivations du dernier chapitre de cette
th\`ese est de r\'epondre \`a la question suivante~: peut-on construire
des exemples analogues avec $K_X$ gros ?
Un r\'esultat de J.\ Koll\'ar montre
que ce n'est pas possible en dimension $3$~:
nous montrons en revanche que de tels exemples
existent en dimension sup\'erieure.
Ceci nous conduira naturellement \`a \'etudier
la structure du centre de l'\'eclatement projectif donn\'e
abstraitement par le th\'eor\`eme de Moishezon lorsque
la vari\'et\'e devient projective apr\`es un \'eclatement
seulement.

\chapter{In\'egalit\'es de Morse singuli\`eres}

Le but central de ce chapitre est d'\'etendre les in\'egalit\'es
de Morse
holomorphes de J.-P.\ Demailly au cas d'un
fibr\'e en droites $E$ muni d'une m\'etrique singuli\`ere au dessus
d'une vari\'et\'e complexe compacte $X$.
Ces in\'egalit\'es nous permettent ensuite de caract\'eriser
analytiquement les vari\'et\'es de Moishezon.
Enfin, nous donnons une version alg\'ebrique
singuli\`ere de nos in\'egalit\'es de Morse.

\section{M\'etriques singuli\`eres}

La notion de m\'etrique singuli\`ere pour des fibr\'es
en droites
a \'et\'e introduite par
J.-P.\ Demailly, A.\ Nadel et H.\ Tsuji. Nous commen\c cons
ce paragraphe en rappelant les premi\`eres
d\'efinitions et les exemples classiques.
Une r\'ef\'erence est \cite{Dem90}.

\subsection{Premi\`eres d\'efinitions}

Soit $X$ une vari\'et\'e et $E$ un fibr\'e en droites
sur $X$. Une m\'etrique hermitienne {\bf singuli\`ere}
sur $E$ est donn\'ee localement
sur un ouvert trivialisant $U_{\alpha}$ par
$$h(\xi _{x}) = ||\xi _x||_h := |\xi|\exp(-\varphi _{\alpha} (x))$$ o\`u la
fonction r\'eelle $\varphi _{\alpha}$ est
seulement suppos\'ee {\bf localement int\'egrable}.

Cette derni\`ere hypoth\`ese suffit \`a donner
encore un sens \`a la notion de courbure
d'un tel fibr\'e~: en effet, on pose
toujours $\displaystyle{\Theta(E) :=
\frac{i}{\pi} \partial \overline{\partial} \varphi _{\alpha}}$
o\`u le $\partial \overline{\partial}$ est pris au sens
des distributions. L'objet ainsi d\'efini
n'est plus une forme ${\cal C}^{\infty}$ mais
un {\bf courant} (appel\'e courant de courbure)
de bi-degr\'e $(1,1)$. Le lemme de Dolbeault-Grothendieck
\'etant vrai pour les courants, la cohomologie de
De Rham est calculable aussi bien avec les courants
qu'avec les formes et la classe de cohomologie
du courant de courbure $\Theta(E)$ est \'egale
comme dans le cas lisse \`a la premi\`ere classe de Chern
du fibr\'e $E$.

\medskip

\noindent {\bf Exemples } Les deux exemples suivants
jouent un r\^ole essentiel.

(i) Si $\displaystyle{D = \sum \alpha _j D_j}$ est un diviseur
de $X$ et si $g_j$ est l'\'equation locale de $D_j$
sur un ouvert $U_{\alpha}$, alors la fonction
$\displaystyle{\varphi _{\alpha} = \sum \alpha _j \log |g_j|}$
d\'efinit une m\'etrique singuli\`ere naturelle sur
le fibr\'e en droites ${\cal O}(D)$ associ\'e au
diviseur $D$.

\noindent Pour cette m\'etrique, l'\'equation de
Lelong-Poincar\'e est $\displaystyle{\Theta({\cal O}(D)) = [D]}$
o\`u $[D]$ d\'esigne le courant d'int\'egration
sur le diviseur $D$.

\noindent Remarquons que dans le cas o\`u $X$ est une courbe complexe
et $E$ le fibr\'e tangent de $X$, la construction pr\'ec\'edente
donne une m\'etrique plate avec des masses de Dirac de courbure~:
ce type de m\'etrique est bien connu en g\'eom\'etrie
riemannienne sous le nom de {\bf m\'etrique plate
\`a singularit\'es coniques}.

(ii) Soit $E$ un fibr\'e en droites sur une
vari\'et\'e $X$ et soient $s_1,\ldots \!,s_p$
des sections de $E^{\otimes k}$.
Alors,
$$||\xi_x||^2 :=
\left( \frac{|\theta _{\alpha}(\xi)|^2}{|\theta _{\alpha}(s_1(x))|^2+
\cdots+|\theta _{\alpha}(s_p(x))|^2} \right) ^{1/k}$$ o\`u
$\theta _{\alpha}$ est une trivialisation locale de $E$ et $E^{\otimes k}$,
d\'efinit une m\'etrique singuli\`ere sur le fibr\'e $E$.
La fonction $\varphi _{\alpha}$ est ici $\displaystyle{\varphi _{\alpha}(x) =
\frac{1}{2k} \log(|\theta _{\alpha}(s_1(x))|^2+
\cdots+|\theta _{\alpha}(s_p(x))|^2)}$.

\medskip

Comme dans le cas ${\cal C}^{\infty}$, nous d\'esirons
pouvoir donner un sens \`a l'expression ``\^etre
\`a courbure positive ou strictement positive"
pour un fibr\'e en droites muni d'une m\'etrique
singuli\`ere.

Pour cela, rappelons qu'un courant
$T$ de bi-degr\'e $(1,1)$ est {\bf positif}
si pour toutes formes $\alpha_1,\ldots \!,\alpha_{n-1}$ de type
$(1,0)$ et
${\cal C}^{\infty}$ \`a support compact, on a~:
$$ < T , i \alpha_1 \wedge \overline{\alpha}_1\wedge \ldots \wedge
i \alpha_{n-1} \wedge \overline{\alpha}_{n-1} > \ \  \geq \ \ 0,$$
o\`u $<\cdot, \cdot >$ est la dualit\'e entre les courants et les formes.
On note alors $T \geq 0$.

De m\^eme, on dira qu'un courant
$T$ de bi-degr\'e $(1,1)$ est {\bf strictement positif}
s'il existe une m\'etrique $\omega$ hermitienne ${\cal C}^{\infty}$
sur $X$ telle que $T - \omega$ est un courant positif.
On note dans ce cas $T >0$.

Nous introduisons alors la d\'efinition suivante~:

\medskip

\noindent {\bf D\'efinition } {\em Un fibr\'e en droites muni
d'une m\'etrique singuli\`ere est positif (respectivement
strictement positif) si le courant de courbure associ\'e
est positif (respectivement strictement positif).}

\medskip

En reprenant les notations de l'exemple ci-dessus,
le fibr\'e en droites ${\cal O}(D)$ muni de sa
m\'etrique singuli\`ere est positif si et seulement
si le diviseur $D$ est effectif (i.e tous les coefficients
$\alpha _j$ sont positifs ou nuls). L'exemple (ii)
d\'efinit toujours un courant de courbure
positif car une fonction de la forme
$\varphi = \log \sum |f_j|^2$ o\`u les
$f_j$ sont des fonctions holomorphes est plurisousharmonique
(en abr\'eg\'e {\bf psh}).

\subsection{Faisceaux d'id\'eaux multiplicateurs de Nadel}

L'\'etude d'un fibr\'e en droites muni d'une
m\'etrique singuli\`ere est facilit\'ee par l'utilisation
d'un outil pertinent associ\'e aux singularit\'es de
la m\'etrique~: il s'agit d'un faisceau d'id\'eaux
introduit par A.\ Nadel.

\medskip

\noindent {\bf D\'efinition } {\em
Soit $(E,h)$ un fibr\'e en droites sur une vari\'et\'e $X$.
On appelle ``faisceau multiplicateur de Nadel" le faisceau
d'id\'eaux ${\cal I}(h)$ des germes de fonctions holomorphes
$L^2$ par rapport au poids de la m\'etrique singuli\`ere,
i.e l'ensemble des germes $f \in {\cal O}_{X,x}$ tels
que $|f|^2 \exp (-2\varphi)$ est int\'egrable par rapport
\`a la mesure de Lebesgue dans des coordonn\'ees
locales au voisinage de $x$.

Plus g\'en\'eralement, si $\varphi$ est une fonction
r\'eelle
sur un ouvert $\Omega$, nous notons aussi
${\cal I}(\varphi)$ le faisceau des germes
$f \in {\cal O}_{\Omega}$ tels
que $|f|^2 \exp (-2\varphi)$ est localement int\'egrable.}

\medskip

La propri\'et\'e essentielle satisfaite par ce faisceau
d'id\'eaux, due \`a Nadel, est que
{\em si $\varphi$ est une fonction plurisousharmonique,
alors ${\cal I}(\varphi)$ est un faisceau coh\'erent.}

\medskip

Avant de donner des exemples, mentionnons
le r\'esultat suivant qui g\'en\'eralise
le th\'eor\`eme de Kawamata-Viehweg~: il s'agit du
th\'eor\`eme d'annulation de Nadel \cite{Nad89} \cite{Dem89}.

\bigskip

\noindent {\bf Th\'eor\`eme d'annulation de Nadel}
{\em Soit $X$ une vari\'et\'e k\"ahl\'erienne compacte,
et $E$ un fibr\'e en droites muni d'une m\'etrique
singuli\`ere \`a courbure
strictement positive. Alors, pour tout $q \geq 1$, on a~:
$$ H^q (X,{\cal O}(E+K_X) \otimes {\cal I}(h)) = 0.$$
}

\medskip

\vspace{-2mm}
Ce r\'esultat montre que pour g\'en\'eraliser un th\'eor\`eme
d'annulation
dans le contexte des fibr\'es en droites munis
de m\'etrique singuli\`ere, une bonne m\'ethode
consiste \`a ``tordre" la cohomologie
du fibr\'e par le faisceau de Nadel.

\medskip

Donnons quelques exemples de faisceaux multiplicateurs,
qui bien que tr\`es simples jouent un r\^ole important.

\medskip

\noindent {\bf Exemples}

(i) Soit $\varphi$ une fonction r\'eelle sur un ouvert $\Omega$
de $\CC ^n$
contenant l'origine. Si $\varphi$ est minor\'ee au voisinage
de l'origine, alors pour tout $x$ proche de $0$, on
a clairement ${\cal I}(\varphi)_x = {\cal O}_{\Omega,x}$.
En particulier, si $\varphi$ est continue sur $\Omega$,
alors ${\cal I}(\varphi) = {\cal O}_{\Omega}$.

(ii) Pla\c cons nous dans $\displaystyle{\CC ^n}$ au voisinage de
l'origine et, pour $\alpha_{1},\ldots \!,\alpha_{p}$ des
r\'eels positifs et $k$ un entier naturel, posons~:
$$ \varphi _k (z)=
\frac{k}{2}\log(|z_{1}|^{2\alpha_{1}}+ \cdots +|z_{p}|^{2\alpha_{p}}).$$
Alors, ${\cal I}(\varphi _k)_{\CC ^n,0}$ est
${\cal O}_{\CC ^n,0}$-engendr\'e
par les $\displaystyle{\prod _{j=1}^{p}z_{j}^{\beta_{j}}}$ tels que~:
$$ \sum_{j=1}^{p}\frac{\beta_{j}+1}{\alpha_{j}} > k .$$
\indent En particulier, si tous les $\alpha_{i}$ sont \'egaux
\`a $\alpha$, alors ${\cal I}(\varphi _k)_{\CC ^n,0}$ est
${\cal O}_{\CC ^n,0}$-engendr\'e
par les $\displaystyle{\prod _{j=1}^{p}z_{j}^{\beta_{j}}}$ o\`u
$$\displaystyle{\sum_{j=1}^{p}\beta_{j} \geq \lfloor k\alpha - p \rfloor +1}$$
($\lfloor x \rfloor$ d\'esigne la partie enti\`ere du r\'eel $x$).
Autrement dit,
$\displaystyle{{\cal I}(\varphi _k) = {\cal I}_{Y}^{
\lfloor k\alpha \rfloor -p+1}}$
o\`u $Y$ est la sous-vari\'et\'e $\{ z_{1}= \cdots =z_{p}=0 \}$
et ${\cal I}_{Y}$ son id\'eal annulateur.

(iii) Soit $\varphi$ une fonction de la forme
$\sum \alpha _j \log |g_j|$ o\`u les $\alpha _j$
sont des r\'eels positifs et o\`u les fonctions
holomorphes $g_j$ sont telles que les $D_j := g_j^{-1}(0)$
soient des diviseurs irr\'eductibles lisses
se coupant transversalement (on dit alors
que $\displaystyle{D = \sum \alpha _j D_j}$ est un diviseur
{\bf lisse \`a croisements normaux}).
Dans ce cas, le faisceau multiplicateur
${\cal I}(\varphi)$
s'identifie au faisceau inversible de rang un
${\cal O}(-\sum_{j} \lfloor \alpha _j \rfloor D_{j})$.

\medskip

\noindent {\bf D\'emonstration}

Le point (i) est trivial.

Pour (ii), il s'agit d'estimer l'int\'egrale suivante
sur un voisinage de $0$~:
$$ I := \int_{D(0,\varepsilon) ^n}
\frac{|\sum a_{\beta}z_1 ^{\beta _1} \ldots z_n ^{\beta _n}|^2}{
(|z_{1}|^{2\alpha_{1}}+ \cdots +|z_{p}|^{2\alpha_{p}})^k} \ d\lambda(z).$$
Le passage en coordonn\'ees polaires $z_j = \rho _j e^{i\theta _j}$
donne
$$ I = (2\pi)^n \sum |a_{\beta}|^2
\int_{[0,\varepsilon]^n}
\frac{\rho_1^{2\beta_{1}+1}\cdots\rho_n^{2\beta_{n}+1}}
{(|\rho_{1}|^{2\alpha_{1}}+
\cdots +|\rho_{p}|^{2\alpha_{p}})^k}\ d\rho_1 \ldots d\rho_n,$$
puis
$$ I = (2\pi)^n \sum |a_{\beta}|^2
\frac{\varepsilon^{2\beta_{p+1}+2}\ldots \varepsilon^{2\beta_{n}+2}}
{(2\beta_{p+1}+2)\ldots (2\beta_{n}+2)}
\int_{[0,\varepsilon]^p}
\frac{\rho_1^{2\beta_{1}+1}\cdots\rho_p^{2\beta_{p}+1}}
{(\rho_{1}^{2\alpha_{1}}+
\cdots +\rho_{p}^{2\alpha_{p}})^k}\ d\rho_1 \ldots d\rho_p.$$
Par le changement de variables $u_j = \rho_j^{\alpha_j}$, l'int\'egrale
ci-dessus est \'egale \`a
$$ \int_{[0,\varepsilon]^p}
\frac{u_1^{((2\beta _1 +2)/\alpha_1) -1}\ldots u_p^{((2\beta _p
+2)/\alpha_p)-1}}
{(u_1^2+ \cdots + u_p^2)^k} \ du_1 \ldots du_p,$$
Par homog\'en\'eit\'e, cette
derni\`ere int\'egrale converge si et seulement si
l'int\'egrale
$$\displaystyle{ \int_0^{\varepsilon}
\frac{ t^{ 2\sum_{j=1}^p ((\beta _j +1)/\alpha_j) - p}}{t^{2k}}t^{p-1} \ dt}$$
converge, soit
$\displaystyle{2\sum_{j=1}^p \frac{\beta _1 +1}{\alpha_1} - p -2k + p - 1 >
-1}$.
Ceci donne bien la condition annonc\'ee.

Pour le point (iii), il s'agit de d\'eterminer
le crit\`ere pour qu'une fonction de la forme
$$\displaystyle{
\frac{|f|^{2}}{|g_1|^{2\alpha_1}\ldots |g_n|^{2\alpha_n}}}$$
\noindent soit dans $L_{\mbox{\scriptsize  loc}}^{1}$.
Comme les $g_j$ fournissent
des coordonn\'ees locales, et si $p_{j}$ d\'esigne l'ordre d'annulation de
$f$ le long
de $D_j = \{ g_{j}=0 \}$, la condition n\'ecessaire
et suffisante est que $2p_{j}-2\alpha _j > -2$, soit
$p_{j} \geq \lfloor \alpha _j \rfloor$. Comme les sections du fibr\'e
${\cal O}(-D_j)$ s'identifient aux fonctions
holomorphes s'annulant \`a un ordre sup\'erieur ou
\'egal \`a un le long de $D_j$, le r\'esultat en d\'ecoule.\finpreuve

\section{In\'egalit\'es de Morse singuli\`eres}

Dans cette partie, nous \'enon\c cons et d\'emontrons notre
version singuli\`ere des in\'e\-galit\'es de Morse
holomorphes.

\subsection{\'Enonc\'e du r\'esultat principal}

Dans tout ce qui suit, $X$ d\'esigne une vari\'et\'e
compacte et $(E,h)$ un fibr\'e en droites sur $X$
muni d'une m\'etrique singuli\`ere.

A $(E,h)$, on associe pour tout entier $k$ positif le
faisceau d'id\'eaux de Nadel donn\'e par la
m\'etrique singuli\`ere $h^k$ du fibr\'e $E^k$~;
il s'agit ici du faisceau
des germes de fonctions holomorphes $f$ telles que
$|f|^{2}\exp(-2k\varphi)$
est $L_{\mbox{\scriptsize  loc}}^{1}$ (o\`u $\exp (-\varphi)$ est
le poids local de $h$). Nous notons ce faisceau ${\cal I}_{k} (h)$.

Pour des raisons qui apparaitront plus loin, nous sommes
contraints de n'accepter qu'un certain type
de singularit\'es, que nous appelons
{\bf ``singularit\'es analytiques"}. Plus pr\'ecis\'ement, nous
faisons l'hypoth\`ese suivante sur $h$ (c'est-\`a-dire localement
sur $\varphi$ )~:

\bigskip

\noindent {\bf Hypoth\`ese ({\cal S}) :}
{\em la fonction $\varphi$ s'\'ecrit localement~:
$$ \varphi=\frac{c}{2}\log(\sum \lambda_{j}|f_{j}|^{2})+\psi,$$
o\`u les $f_{j}$ sont holomorphes,
les $\lambda_{j}$ sont des fonctions r\'eelles positives ${\cal C}^{\infty}$
sans z\'eros communs,
$\psi$ est ${\cal C}^{\infty}$
et $c$ est un rationnel positif ou nul. Mentionnons que
la somme $\sum \lambda_{j}|f_{j}|^{2}$ peut poss\'eder une
infinit\'e de termes.}

\medskip

Cette hypoth\`ese implique que la fonction $\varphi$ est {\bf quasi
plurisousharmonique}
(ce qui signifie que son hessien complexe est minor\'e par une
$(1,1)$-forme \`a coefficients continus) et donc en particulier que le
faisceau $ {\cal I}_{k} (h)$
est un faisceau coh\'erent d'apr\`es le r\'esultat de A.\ Nadel
pr\'ec\'edemment rencontr\'e.

En terme de courbure, l'hypoth\`ese ({\cal S}) implique en particulier
que le
courant de courbure
$\displaystyle{ \Theta (E)=
\frac{i}{\pi} \partial \overline{\partial} \varphi}$
est quasi positif (c'est-\`a-dire minor\'e
par une $(1,1)$-forme \`a coefficients continus).

\bigskip

\noindent {\bf Remarque }
Signalons d\`es maintenant que cette hypoth\`ese est
``raisonnable" car c'est pr\'ecis\'ement par des fonctions
ayant ce type de
singularit\'es que J.-P.\ Demailly approche une fonction
quasi plurisousharmonique
quelconque \cite{Dem92}.
Ceci aura une importance capitale dans les applications.
Remarquons aussi que tous les exemples rencontr\'es
jusqu'\`a pr\'esent satisfont l'hypoth\`ese ({\cal S}).
Nous renvoyons ici au paragraphe 2.3.1 pour une discussion
plus compl\`ete de l'origine de ce type de
singularit\'es autoris\'ees.

\medskip

Dans notre contexte, nous introduisons les notations
suivantes~:

\medskip

\noindent {\bf Notations }
Le symbole $X(q,E)$ d\'esigne l'ouvert de $X$ form\'e des points
$x$ au voisinage desquels $\varphi$ est born\'ee et $\Theta (E) (x)$
poss\`ede
exactement $q$ valeurs propres strictement n\'egatives et $n-q$ valeurs
propres strictement positives~;
finalement et comme dans le cas lisse $X(\leq q,E)
= X(0,E)\cup \cdots \cup X(q,E)$.

\medskip

Notre r\'esultat est le suivant~:

\bigskip

\noindent {\bf Th\'eor\`eme B } {\em
Soit $(E,h)$ \`a singularit\'es analytiques sur
une vari\'et\'e compacte $X$ de dimension $n$ et soit $F$
un fibr\'e
de rang $r$ sur $X$. Alors pour tout $q$ compris entre
$0$ et $n$~:

\noindent \ (i) $\displaystyle{\
\sum _{j=0}^{q}(-1)^{q-j} \dim H^{j}(X,{\cal O}(E^{k}\otimes F) \otimes
{\cal I}_{k}(h)) \leq r\frac{k^{n}}{n!} \int _{X(\leq q,E)} (-1)^{q} \Theta
(E)^{n} + o(k^{n}) }$,\\
\noindent (a\-vec \'egalit\'e si $q=n$),

\noindent \ (ii) $\displaystyle{ \  \dim H^q(X,{\cal O}(E^{k}\otimes F) \otimes
{\cal I}_{k}(h)) \leq r\frac{k^{n}}{n!} \int _{X(q,E)} (-1)^{q} \Theta
(E)^{n} + o(k^{n}).}$}

\bigskip

Comme dans le cas o\`u la m\'etrique est lisse, le point
(ii) est cons\'equence imm\'ediate du point (i).

\medskip

Remarquons que comme dans le th\'eor\`eme d'annulation de Nadel,
le ph\'enom\`ene nouveau par rapport au cas ${\cal C}^{\infty}$
est la pr\'esence des faisceaux d'id\'eaux~; les groupes de
cohomologie que nous estimons sont les groupes de cohomologie
\`a valeurs dans les grandes puissances de $E$, tordues
par la suite des faisceaux multiplicateurs naturellement
associ\'ee aux singularit\'es de la m\'etrique.

Evidemment, si la m\'etrique est lisse, nous retrouvons les
in\'egalit\'es de J.-P.\ Demailly car tous
les ${\cal I}_{k}(h)$ sont \'egaux au faisceau structural ${\cal O}_X$.
Cependant, mentionnons d\`es maintenant que notre d\'emonstration
repose sur les in\'egalit\'es dans le cas ${\cal C}^{\infty}$ !

\subsection{D\'emonstration du th\'eor\`eme B}

\subsub{Plan de la preuve}

La d\'emarche suivie pour d\'emontrer nos
in\'egalit\'es est la suivante~:

a) apr\`es \'eclatement de $X$ le long de sous-vari\'et\'es d\'efinies par les
singularit\'es de $h$, on se ram\`ene \`a un fibr\'e muni d'une m\'etrique
lisse~; on peut alors appliquer les in\'egalit\'es de Morse holomorphes
dans le cas $h$ lisse \`a la vari\'et\'e obtenue $\tilde{X}$,

b) on relie les groupes de cohomologie sur $\tilde{X}$
\`a ceux de $X$. Pour cela, nous \'etudions
le comportement des faisceaux multiplicateurs
par rapport aux \'eclatements en montrant que
les dimensions des groupes de cohomologie
associ\'es sont asymptotiquement de m\^eme
dimension.

\subsub{R\'eduction au cas lisse}

Commen\c cons par expliquer la premi\`ere partie
de la preuve.

Pour cela, il est bon de remarquer que l'hypoth\`ese ({\cal S})
implique
que les singularit\'es de la m\'etrique
$h$ sont localis\'ees le long d'un ensemble analytique, d\'efini
localement par $ \{ x|\ \forall j,\ f_{j}(x)=0 \}$.
Comme nous l'avons vu dans les exemples, cet ensemble analytique
n'est pas n\'ecessairement irr\'eductible et ses composantes irr\'eductibles
sont en g\'en\'eral
de dimension quelconque.

La r\'eduction au cas lisse consiste dans un premier temps \`a se ramener
\`a une vari\'et\'e $\tilde{X}$, obtenue en \'eclatant $X$ le long de centres
lisses
$\pi : \tilde{X} \to X$ de telle sorte que la m\'etrique
$\tilde{h}=\pi^{*}h$ sur le fibr\'e $\tilde{E}=\pi^{*}E$ n'ait ses
singularit\'es qu'en codimension $1$ (ou, de fa\c con
\'equivalente, de telle
sorte que le faisceau $\tilde{{\cal I}}_{k} (\tilde{h})$ soit inversible).

Dans un deuxi\`eme temps, nous appliquerons les in\'egalit\'es de Morse
dans le cas lisse.

\bigskip

\noindent {\bf a) D\'esingularisation de ${\cal I}_{k} (h)$}

\medskip

Pour rendre les faisceaux ${\cal I}_{k} (h)$ localement libres,
l'id\'ee (classique) est d'\'eclater l'id\'eal
``engendr\'e par les $f_j$". Cependant, une difficult\'e
appara\^{\i}t ici car la donn\'ee des $f_j$ est locale sur $X$
et la notion d'id\'eal engendr\'e par ces fonctions
n'a donc pas de sens {\em a priori}. Nous avons cependant la proposition
suivante~:

\medskip

\noindent {\bf Proposition } {\em
Il existe un faisceau d'id\'eaux global
${\cal J}$, qui co\"{\i}ncide avec la cl\^oture int\'egrale
du faisceau
d'id\'eaux engendr\'e par les $f_{j}$ sur chaque ouvert
o\`u $\varphi$ s'\'ecrit
comme dans l'hypoth\`ese ({\cal S}).
}

\medskip

\noindent {\bf D\'emonstration}

En effet, c'est un fait bien connu (r\'esultant par
exemple du th\'eor\`eme de
Brian\c con-Skoda \cite{BSk74}, voir aussi \cite{L-T74}) que la cl\^oture
int\'egrale du faisceau
d'id\'eaux engendr\'e par les $f_{j}$ est donn\'ee par~:
$${\cal J}_{x} = \{ f \in {\cal O}_{X,x} \ ; \ \exists C >0 \ ,
\ |f(z)| \leq C\exp(\frac{1}{c}\varphi(z)) \ \mbox{au voisinage de}\ x \}.$$
Si maintenant $\varphi _{\alpha}$ et $\varphi _{\beta}$
d\'esignent les poids de la m\'etrique sur des ouverts
trivialisants $U_{\alpha}$ et $U_{\beta}$, alors
sur l'intersection $U_{\alpha} \cap U_{\beta}$,
on a $\varphi _{\alpha} = \varphi _{\beta} + O(1)$.
La caract\'erisation ci-dessus implique donc bien que
${\cal J}$ est d\'efini globalement sur $X$.\finpreuve

\medskip

Nous sommes en mesure de d\'emontrer la proposition
suivante~:

\medskip

\noindent {\bf Proposition } {\em
Sous les hypoth\`eses pr\'ec\'edentes, il existe une vari\'et\'e
$\tilde{X}$
et $\pi : \tilde{X} \rightarrow X$ une compos\'ee d'un nombre fini
d'\'eclate\-ments de centres lisses tels que le
fibr\'e $\tilde{E} :=\pi^{*}E$ muni
de la m\'etrique singuli\`ere $\tilde{h}=\pi^{*}h$
de poids local $\exp(-\tilde{\varphi})$
v\'erifie la propri\'et\'e suivante :
pour tout $x_{0} \in \tilde{X}$, il existe des coordonn\'ees holomorphes
$w_{1},\ldots \!,w_{n}$ centr\'ees en
$x_{0}$ et une fonction $\tilde{\psi}$ de classe ${\cal C}^{\infty}$
telles que~:
$$ \tilde{\varphi} (w) = c \sum_{j}a_{j}\log|g_{j}(w)|+\tilde{\psi}(w)$$
o\`u les $a_{j}$ sont des entiers positifs ou nuls et o\`u
les $g_{j}$ sont
irr\'eductibles dans ${\cal O}_{\tilde{X},x_{0}}$
et d\'efinissent un diviseur lisse \`a croisements normaux.}

\medskip

\noindent {\bf D\'emonstration}

\'Eclatons l'id\'eal ${\cal J}$ de sorte que l'image inverse
$\pi^{-1}{\cal J}.{\cal O}_{X'}$ soit un faisceau inversible.
Par le th\'eor\`eme d'aplatissement d'Hironaka \cite{Hir75},
on peut dominer
cet \'eclatement par
une vari\'et\'e $\tilde{X}$ obtenue par une suite d'\'eclatements de
centres lisses dans $X$, et l'image inverse de ${\cal J}$ est
toujours inversible !

\noindent Mais alors, la m\'etrique image r\'eciproque sur $\pi^{*}E$ est
donn\'ee
par

\begin{eqnarray*}
\ \tilde{\varphi} & = &
\frac{c}{2}\log(\sum (\lambda_{j}\circ\pi)|f_{j}\circ\pi|^{2})+
\psi\circ\pi \\
 & = & \frac{c}{2}\log(|g|^{2}) +
\frac{c}{2}\log(\sum (\lambda_{j}\circ\pi)|h_{j}|^{2})+
\psi\circ\pi \\
& = & \frac{c}{2}\log(|g|^{2}) + \tilde{\psi},
\end{eqnarray*}

\noindent o\`u on a not\'e $g$ le g\'en\'erateur local du faisceau d'id\'eaux
engendr\'e
par les $f_{j}\circ\pi$. Si la d\'ecomposition de $g$ en facteurs
irr\'eductibles
dans $\displaystyle{{\cal O}_{\tilde{X},x_{0}}}$ s'\'ecrit~:
$\displaystyle{g = \prod_{j}g_{j}^{a_{j}}}$,
on a le r\'esultat apr\`es
application du th\'eor\`eme de d\'esingularisation
d'Hironaka \cite{Hir64}
pour rendre le diviseur d\'efini par les $g_{j}$ \`a croisements
normaux.\finpreuve

\medskip

Dans tout ce qui suit, les objets utilis\'es sont ceux obtenus
apr\`es application de la proposition pr\'ec\'edente.

Notons $\displaystyle{\tilde{{\cal I}}_{k} (\tilde{h})}$ le faisceau
d'id\'eaux des germes
de fonctions holomorphes sur $\tilde{X}$
telles que $\displaystyle{|f|^{2}e^{-2k\tilde{\varphi}}}$
est $L_{\mbox{\scriptsize loc}}^{1}$ (c'est le faisceau multiplicateur
de Nadel associ\'e au fibr\'e $\pi^*E$).
Si nous notons $b_{j,k} = \lfloor ca_{j}k \rfloor$ et $\tilde{D}_{j}$
le diviseur
d\'efini localement par $g_{j}=0$, alors la d\'etermination
du faisceau $\displaystyle{\tilde{{\cal I}}_{k} (\tilde{h})}$
rel\`eve des exemples
pr\'ec\'edents si bien que le lemme suivant en d\'ecoule~:

\medskip

\noindent {\bf Lemme }{\em Sous les conditions pr\'ec\'edentes,
le faisceau d'id\'eaux
$\tilde{{\cal I}}_{k} (\tilde{h})$
s'identifie au faisceau inversible de rang un
${\cal O}(-\sum_{j}b_{j,k}\tilde{D}_{j})$.}

\medskip

\noindent {\bf b) Exemples }

\medskip

Illustrons ce qui pr\'ec\`ede en reprenant les notations
de l'exemple (ii) du 2.1.2.
Dans ces cas \'evidemment simples, il n'est nul besoin d'appliquer les
th\'eor\`emes d'Hironaka~: on explicite directement le choix des
\'eclatements !

\medskip

Si on suppose que tous les $\alpha_{i}$ sont \'egaux
\`a $\alpha$, nous avons vu alors que
${\cal I}(\varphi_{k} )$ est \'egal \`a
${\cal I}_{Y}^{ \lfloor k\alpha \rfloor -p+1}$
o\`u $Y$ est la sous-vari\'et\'e de codimension
$p$ donn\'ee par $\{ z_{1}= \cdots =z_{p}=0 \}$.

\noindent Si $p=1$, le faisceau d'id\'eaux est d\'ej\`a inversible, sinon
\'eclatons $\CC ^n$ le long de $Y$. L'expression
de la nouvelle m\'etrique est donn\'ee dans la premi\`ere carte
par $$\displaystyle{\tilde{\varphi}(w)= \alpha\log(|w_{1}|) +
\frac{1}{2}\log(1+|w_{2}|^{2\alpha}+\cdots +|w_{p}|^{2\alpha})}$$ si bien que
$\tilde{{\cal I}}(\varphi_{k}) = {\cal O}(- \lfloor k\alpha \rfloor D)$ o\`u
$D$
est le diviseur exceptionnel de l'\'eclatement.
Il suffit donc dans ce cas d'un \'eclatement en codimension
$p$ pour obtenir le r\'esultat
souhait\'e.

\medskip

De m\^eme, si $\alpha$ est un entier positif et si
$\displaystyle{\varphi(z) =
\frac{1}{2}\log(|z_{1}|^{2}+|z_{2}|^{2\alpha})}$
dans
$\CC ^n$, il faut cette fois $\alpha$ \'eclatements en codimension~2.

\noindent En effet, \'eclatons $\CC ^n$ le long de $\{ z_{1}=z_{2}=0 \}$.
L'expression
de la nouvelle m\'etrique est donn\'ee dans la premi\`ere carte
par $$\displaystyle{\tilde{\varphi}(w)= \log(|w_{1}|) +
\frac{1}{2}\log(1+|w_{1}w_{2}|^{2\alpha})}$$ qui est de la forme voulue alors
qu'on obtient dans la deuxi\`eme carte
$$\displaystyle{\tilde{\varphi}(w)= \log(|w_{2}|) +
\frac{1}{2}\log(|w_{1}|^{2}+|w_{2}|^{2(\alpha-1)})}.$$
On \'eclate alors dans la deuxi\`eme carte le long de
$\{ w_{1}=w_{2}=0 \}$. En r\'ep\'etant ce proc\'ed\'e $\alpha$ fois,
on obtient une m\'etrique de la forme souhait\'ee en tout point.

\noindent D\'ecrivons le faisceau d'id\'eaux obtenu~: pour
tout $j$ compris entre $1$ et $\alpha$, notons $D_{j}$ la transform\'ee
stricte dans $\tilde{X}$ du diviseur exceptionnel du $j$-i\`eme \'eclatement.
Alors $D_{j}$ et $D_{j+1}$ se coupent transversalement et on a
$\tilde{{\cal I}}_{k} (\tilde{\varphi}) = {\cal O}(-kD)$ o\`u $D$
d\'esigne le diviseur \`a croisements normaux
$\displaystyle{D=\sum_{j=1}^{\alpha}jD_{j}}$.

\medskip

\noindent {\bf c) In\'egalit\'es de Morse sur $\tilde{X}$ }

\medskip

On montre maintenant comment appliquer
les in\'egalit\'es de Morse holomorphes
dans le cas ${\cal C}^{\infty}$ sur $\tilde{X}$ au
fibr\'e en droites $\tilde{E} = \pi ^{\ast}E$ muni de
la m\'etrique image r\'eciproque. Pour cela, nous devons
montrer que $(\tilde{E})^k \otimes \tilde{{\cal I}}_{k} (\tilde{h})$
peut essentiellement s'\'ecrire comme la $k$-i\`eme
puissance tensorielle d'un fibr\'e en droites hermitien fixe.

Notons $\displaystyle{c=\frac{u}{v}}$ et supposons que $k = vk'$ est
un multiple du
d\'enomina\-teur de $c$. Comme $b_{j,k}=ca_{j}k$,
le faisceau multiplicateur $\displaystyle{\tilde{{\cal I}}_k (\tilde{h})}$
est \'egal au faisceau inversible
$\displaystyle{{\cal O}(-k'\tilde{D})}$
o\`u $\displaystyle{\tilde{D}=u\sum_{j}a_{j}\tilde{D}_{j}}$.

Avec ces notations, on a la~:

\medskip

\noindent {\bf Proposition } {\em Notons $\hat{E}$
le fibr\'e $\tilde{E}^{v} \otimes {\cal O}(-\tilde{D})$. Alors~:

(i)  pour tout $k = vk'$, on a
$(\hat{E})^{k'} = \tilde{E}^{k} \otimes \tilde{{\cal I}}_{k} (\tilde{h}),$

(ii) la m\'etrique hermitienne sur $\hat{E}$, produit
de la m\'etrique $\tilde{h} ^v$ et de la m\'etrique singuli\`ere
naturelle sur ${\cal O}(-\tilde{D})$, est une m\'etrique
hermitienne ${\cal C}^{\infty}$. De plus, et en
dehors des singularit\'es de la m\'etrique $\tilde{h}$,
on a l'\'egalit\'e $\Theta(\hat{E}) =
v\Theta(\tilde{E})=\pi^{\ast}\Theta(E)$.
}

\medskip

\noindent {\bf D\'emonstration }

Le point (i) est \'evident. Pour le point (ii), la m\'etrique
produit naturelle est donn\'ee localement
par le poids
$\displaystyle{\tilde{\chi}(z)=
v \tilde{\phi} (z) -\sum_{j=1}^{n}ua_{j}\log|g_{j}|
= v \tilde{\psi}(z).}$ Ainsi, cette m\'etrique est lisse
et l'\'egalit\'e de courbure d\'ecoule de suite du fait
que $\displaystyle{\partial \overline{\partial} \log|g_j|
= 0}$ l\`a o\`u $g_j$ ne s'annule pas. \finpreuve

\medskip

Cette proposition nous permet d'estimer les groupes de cohomologie qui nous
int\'eressent~:

\medskip

\noindent {\bf Proposition }{\em
On a pour tout $k$ :
$$\sum_{j=0}^{q}(-1)^{q-j}\dim
H^{j}(\tilde{X},{\cal O}(\tilde{E}^{k} \otimes \tilde{F}) \otimes
\tilde{{\cal I}}_{k} (\tilde{h})) \leq
r\frac{k^{n}}{n!} \int _{X(\leq q,E)} (-1)^{q} \Theta
(E)^{n} + o(k^{n}) $$
o\`u l'int\'egrale est prise sur les points lisses de
la m\'etrique de $E$.
}

\medskip

\newpage

\noindent {\bf D\'emonstration}

D'apr\`es la proposition pr\'ec\'edente, on peut
appliquer les in\'egalit\'es de Morse de Demailly au
fibr\'e $\hat{E}$,
si bien que pour $k=k'v$, on a~:

$$\sum_{j=0}^{q}(-1)^{q-j}\dim
H^{j}(\tilde{X},{\cal O}(\tilde{E}^{k} \otimes \tilde{F}) \otimes
\tilde{{\cal I}}_{k} (\tilde{h})) \leq
r\frac{k'{}^{n}}{n!} \int _{\tilde{X}(\leq q,\hat{E})} (-1)^{q} \Theta
(\hat{E})^{n} + o(k'{}^{n}).$$

\noindent Relions alors l'int\'egrale de courbure sur $\tilde{X}$ \`a une
int\'egrale
de courbure sur $X$.

\noindent Comme $\Theta(\hat{E})=v\Theta(\tilde{E})$ sur
les points lisses de la m\'etrique de $\tilde{E}$ si $k=k'v$,
on a
$$k'{}^{n} \int _{\tilde{X}(\leq q,\hat{E})} \Theta (\hat{E})^{n}
+ o(k'{}^{n})= k^{n} \int _{\tilde{X}(\leq q,\tilde{E})}
 \Theta (\tilde{E})^{n} + o(k^{n}).$$

\noindent Notons alors $S$ la r\'eunion des diviseurs exceptionnels de
$\pi : \tilde{X} \rightarrow X $~; $S$ est n\'egligeable pour la
mesure de Lebesgue et on a

\begin{eqnarray*}
\int _{\tilde{X}(\leq q,\tilde{E})} \Theta
(\tilde{E})^{n} & = & \int _{\tilde{X}(\leq q,\tilde{E}) \setminus S} \Theta
(\tilde{E})^{n}  = \int _{\tilde{X}(\leq q,\tilde{E}) \setminus S} \Theta
(\pi ^{*}E)^{n} \\
 & = & \int _{\tilde{X}(\leq q,\tilde{E}) \setminus S} \pi^{*}(\Theta
(E)^{n}) =  \int _{X(\leq q,E) \setminus \pi(S)} \Theta
(E)^{n} \\
& = & \int _{X(\leq q,E)} \Theta
(E)^{n}.
\end{eqnarray*}

\noindent On a donc pour les entiers $k$ multiples d'un entier fixe~:
$$\sum_{j=0}^{q}(-1)^{q-j}\dim
H^{j}(\tilde{X},{\cal O}(\tilde{E}^{k} \otimes \tilde{F}) \otimes
\tilde{{\cal I}}_{k} (\tilde{h})) \leq
r\frac{k^{n}}{n!} \int _{X(\leq q,E)} (-1)^{q} \Theta
(E)^{n} + o(k^{n}) $$
o\`u l'int\'egrale est prise sur les points lisses de
la m\'etrique de $E$.

\medskip

Pour terminer cette preuve, il suffit de montrer que l'estimation
pr\'ec\'edente est valable sans restriction sur $k$.

\noindent En reprenant les notations du lemme et en
posant $\displaystyle{c=\frac{u}{v}}$,
on a, pour $k=k'v+r$
$$ b_{j,k} = ua_{j}k'+r' $$
o\`u $r'$ ne prend qu'un nombre fini de valeurs enti\`eres. On en d\'eduit
alors

\begin{eqnarray*}
\! \tilde{E} ^{k} \otimes \tilde{F} \otimes
{\cal O}(-\sum_{j}b_{j,k}D_{j}) & = & (\tilde{E}^{v})^{k'}\otimes
\tilde{E}^{r'} \otimes \tilde{F} \otimes {\cal O}(-k'u\sum_{j}a_{j}D_{j})
\otimes {\cal O}(-r'\sum_{j}D_{j})\\
& = &
(\tilde{E}^{v})^{k'}\otimes
{\cal O}(-k'u\sum_{j}a_{j}D_{j}) \otimes \hat{F}_{r'}.
\end{eqnarray*}
On raisonne alors comme pr\'ec\'edemment~: on munit
$\hat{F}_{r'}$ d'une m\'etrique lisse quelconque tandis que
$(\tilde{E}^{v})^{k'}\otimes {\cal O}(-k'u\sum_{j}a_{j}D_{j})$
est muni de la m\'etrique lisse naturelle donn\'ee localement
par $v\tilde{\psi} (z)$.
Ceci d\'emontre la proposition.\finpreuve

\subsub{Lien entre cohomologie sur $\tilde{X}$ et cohomologie sur $X$}
Pour achever la preuve du th\'eor\`eme B, il reste \`a relier les
groupes $$H^{q}(\tilde{X},{\cal O}(\tilde{E}^{k} \otimes \tilde{F}) \otimes
\tilde{{\cal I}}_{k} (\tilde{h}))$$ de la proposition pr\'ec\'edente et les
$$H^{q}(X,{\cal O}(E^{k}\otimes F) \otimes
{\cal I}_{k}(h))$$ qui nous int\'eressent directement.
Nous adoptons ici une d\'emarche un peu diff\'erente de celle
de notre travail \cite{Bo93b}.

Soient $X$ une vari\'et\'e compacte et $\mu : X' \to X$
l'\'eclatement de $X$ le long d'une sous-vari\'et\'e lisse $Y$.
Soit $E$ un fibr\'e en droites sur $X$, muni d'une
m\'etrique hermitienne singuli\`ere $h$ de poids local
$\exp (-\varphi)$. On note ${\cal I}_k(\varphi)$ le
faisceau multiplicateur associ\'e \`a la m\'etrique
sur $E^k$. Soit enfin $F$ un fibr\'e
vectoriel sur $X$. Nous montrons dans ce paragraphe la proposition
suivante~:

\bigskip

\noindent {\bf Proposition } {\em Supposons qu'au voisinage
de tout point du diviseur exceptionnel de $\mu$, la fonction
$\varphi \circ \mu$ s'\'ecrive~:
$$ \varphi \circ \mu = \alpha \log |f| + \psi,$$
o\`u $\alpha $ est strictement positif, $f$ d\'esigne
une \'equation locale du diviseur exceptionnel de $\mu$ et
$\psi$ est une fonction psh.

\noindent Alors, si $k \alpha >1$, on a pour tout entier
$q \geq 0$~:
$$ H^q(X',K_{X'} \otimes \mu^*(E^k \otimes F) \otimes {\cal I}_k(\varphi \circ
\mu))
\simeq H^q(X,K_X \otimes E^k \otimes F \otimes {\cal I}_k(\varphi)) .$$}

\vspace*{-3mm}

\noindent {\bf Commentaires } Nous ne supposons plus dans ce dernier
r\'esultat que la m\'etrique est \`a singularit\'es analytiques,
ni que les faisceaux ${\cal I}_k(\varphi \circ \mu)$ sont inversibles.
Le point important est que l'hypoth\`ese faite sur
l'\'ecriture de $\varphi \circ \mu$ est automatiquement satisfaite
si $\varphi$ est \`a singularit\'es analytiques et le centre
de l'\'eclatement $Y$ est inclus dans le lieu singulier de la m\'etrique.

La fin de la d\'emonstration
de nos in\'egalit\'es de Morse se conclut par application
r\'ep\'et\'ee de cette proposition aux \'eclatements
successifs dont $\pi : \tilde{X} \to X$ est la compos\'ee.\finpreuve

\medskip

Avant de d\'emontrer la proposition, mentionnons le
probl\`eme suivant~:

\medskip

\noindent {\bf Question } {\em Si les singularit\'es de $\varphi$ sont
quelconques,
a-t-on toujours
$$ \dim H^q(X',K_{X'} \otimes \mu^*(E^k \otimes F)
\otimes {\cal I}_k(\varphi \circ \mu))
= \dim H^q(X,K_X \otimes E^k \otimes F \otimes {\cal I}_k(\varphi))
+ o(k^n) ?$$}

\vspace*{-3mm}

\noindent {\bf D\'emonstration de la proposition }

Elle repose tout d'abord sur le fait que les faisceaux
de Nadel se comportent bien par image directe. De fa\c con pr\'ecise,
si $(E,h)$ est un fibr\'e en droites muni d'une
m\'etrique singuli\`ere au dessus d'une vari\'et\'e
$X$, et si $\mu : \tilde{X} \to X$ est une modification,
alors \cite{Dem94}
$$\mu _*(K_{\tilde{X}}\otimes {\cal I}(\mu^*h))
= K_X \otimes {\cal I}(h).$$

\noindent \'Evidemment, ceci ne suffit pas pour d\'emontrer qu'il y a
isomorphisme en cohomologie~: l'obstruction est mesur\'ee
par les images directes sup\'erieures. Ainsi, il suffit
de montrer, gr\^ace au
th\'eor\`eme de Leray, que pour tout entier
$q \geq 1$ et pour tout $k$ assez grand,
$$R^{q}\mu_{*}(K_{X'} \otimes {\cal I}_k(\varphi \circ \mu))
= 0.$$

\noindent On note $r$ la codimension de $Y$, centre de l'\'eclatement.

\noindent Dans ce cas, le faisceau $q$-i\`eme image directe sup\'erieure
$R^{q}\mu_{*}(K_{X'} \otimes {\cal I}_k(\varphi \circ \mu))$
est un faisceau \`a support dans $Y$,
la fibre au dessus d'un point $y$ de $Y$ \'etant \'egale \`a~:
$$ F_{k,y} = \limind_{y \in U}
H^{q}(\mu^{-1}(U),(K_{X'}
\otimes {\cal I}_k(\varphi \circ \mu))_{|\mu^{-1}(U)}),$$
o\`u la limite porte sur les voisinages $U$ de $y$
dans $X$.

\noindent Soit donc $U$ un ouvert de Stein voisinage de $y$
et soit $\omega$ une m\'etrique hermitienne
sur $X'$. Il s'agit de montrer
que pour toute forme $u$ de type $(n,q)$ sur $\mu^{-1}(U)$,
\`a coefficients localement $L^2$ satisfaisant~:

(i) $\displaystyle{\overline{\partial}u = 0}$,

(ii) $\displaystyle{
I := \int _{\mu^{-1}(U)}|u|^2 \exp(-2k \varphi \circ \mu) dV_{\omega} < +\infty
,}$

\noindent il existe une forme $v$ de type $(n,q-1)$
\`a coefficients localement $L^2$ satisfaisant~:

(i)' $\displaystyle{\overline{\partial}v = u}$,

(ii)' $\displaystyle{
\int _{\mu^{-1}(U)}|v|^2 \exp(-2k \varphi \circ \mu) dV_{\omega} < +\infty .}$

\noindent R\'esoudre un tel probl\`eme est en g\'en\'eral possible
gr\^ace aux estimations $L^2$ de H\"ormander. La difficult\'e
ici est que la fonction
$k \varphi \circ \mu$
n'est pas strictement psh au voisinage du diviseur
exceptionnel $D$ de l'\'eclatement. C'est exactement ici que nous utilisons
l'hypoth\`ese faite sur $\varphi$.
\noindent En effet, l'\'egalit\'e
$$ \varphi \circ \mu = \alpha \log |f| + \psi$$
se traduit par~:
$$ {\cal I}_k(\varphi \circ \mu)
= {\cal O}(- \lfloor k\alpha \rfloor D) \otimes {\cal I}
\left( k\psi + (k\alpha - \lfloor k\alpha \rfloor)\log |f|) \right).$$

\noindent Il s'agit alors de montrer l'annulation des groupes
$$H^{n,q}\left (\mu^{-1}(U),{\cal O}(- \lfloor k\alpha \rfloor D)
\otimes {\cal I}(k\psi + (k\alpha - \lfloor k\alpha \rfloor)\log |f|\right ),$$
autrement dit de r\'esoudre le probl\`eme du
$\overline{\partial}$~:

(i)' $\displaystyle{\overline{\partial}v = u}$,

(ii)' $\displaystyle{
\int _{\mu^{-1}(U)}|v|^2 \exp(-2k \psi) dV_{\omega} < +\infty ,}$

\noindent pour des $(n,q)$ formes {\em \`a valeurs dans le fibr\'e
${\cal O}(- \lfloor k\alpha \rfloor D)$}. Or, le fibr\'e
${\cal O}(-D)_{|D} = {\cal O}_D (1)$ est strictement positif sur
les fibres de l'\'eclatement. On peut donc munir
${\cal O}(-\lfloor k\alpha \rfloor D)$
d'une m\'etrique \`a courbure strictement positive au voisinage du
diviseur exceptionnel, et comme $\psi$ peut \^etre suppos\'ee
strictement psh en dehors de $D$, nous sommes maintenant
dans les hypoth\`eses d'application des estimations $L^2$ de
H\"ormander.\finpreuve

\section{Caract\'erisations analytiques des vari\'et\'es \- de \-
Moishezon}
Dans la lign\'ee des conditions suffisantes donn\'ees par Y.-T.\ Siu
et J.-P.\ Demailly pour caract\'eriser les vari\'et\'es de Moishezon,
nous montrons la caract\'erisation suivante~:

\bigskip

\noindent {\bf Th\'eor\`eme C} {\em Une vari\'et\'e compacte $X$ de
dimension $n$
est de Moishezon si et
seulement si il existe sur $X$ un courant $T$ ferm\'e
de bidegr\'e $(1,1)$
tel que~:

(i) $ \displaystyle{\{ T \} \in H^{2}(X,\ZZ) }$,

(ii) $\displaystyle{ T= \frac{i}{\pi} \partial \overline{\partial} \varphi
+ \alpha }$,
o\`u $\varphi$ est une fonction r\'eelle \`a singularit\'es analytiques
et o\`u
$\alpha$ est un repr\'esentant
$ {\cal C}^{\infty}$ de $\{ T \}$,

(iii) $\displaystyle{\int_{X(\leq 1,T)} T^{n} > 0}$ o\`u
l'int\'egrale est prise
sur les points lisses du courant $T$.
}

\bigskip

L'id\'ee de donner une caract\'erisation analytique
des vari\'et\'es de Moishezon en terme de courant de
courbure est aussi pr\'esente dans
un travail de S.\ Ji et B.\ Shiffman \cite{JiS93} simultan\'e
au n\^otre. Nos in\'egalit\'es permettent ainsi de d\'emontrer
le r\'esultat suivant~:

\bigskip

\noindent {\bf Th\'eor\`eme D \cite{JiS93}, \cite{Bo93b}} {\em
Une vari\'et\'e compacte $X$ de
dimension $n$
est de Moishezon si et
seulement si il existe sur $X$ un courant $T$ ferm\'e de bidegr\'e $(1,1)$
tel que~:

(i) $ \displaystyle{\{ T \} \in H^{2}(X,\ZZ) }$,

(ii) le courant $T$ est strictement positif (i.e minor\'e
par une $(1,1)$-forme ${\cal C}^{\infty}$ hermitienne).
}

\bigskip

Ces deux \'enonc\'es fournissent l'extension naturelle aux vari\'et\'es
de Moishezon du th\'eor\`eme de plongement de Kodaira pour les vari\'et\'es
projectives.

Avant de d\'emontrer les r\'esultats ci-dessus,
nous rappelons deux r\'esultats sur les courants.

\subsection{Deux rappels}

Il est bien connu \cite{G-H78} que la seule obstruction
pour qu'une $(1,1)$-forme ferm\'ee r\'eelle ${\cal C}^{\infty}$
soit la forme de courbure d'un fibr\'e en droites hermitien
est que sa classe de cohomologie soit enti\`ere, i.e
appartienne \`a $H^2(X,\ZZ)$. Nous montrons dans la proposition suivante
que ce r\'esultat persiste pour les courants quasi positifs~:

\medskip

\noindent {\bf Proposition  }
{\em Soit $T$ un courant quasi-positif ferm\'e, de bi-degr\'e
$(1,1)$ dont la classe de cohomologie
$\{ T \}$ est dans $H^2(X,\ZZ)$.

Alors, il existe un fibr\'e en droites $E$
muni d'une m\'etrique singuli\`ere dont
le courant de courbure est \'egal \`a $T$.}

\medskip

\noindent {\bf D\'emonstration }

Elle suit la d\'emonstration du cas ${\cal C}^{\infty}$ \cite{S-S85}.
On recouvre la vari\'et\'e $X$ par des ouverts de Stein
contractiles dont les intersections mutuelles
sont elles aussi contractiles. Sur chaque
ouvert $U_{\alpha}$, on \'ecrit
$\displaystyle{T = \frac{i}{\pi}\partial \overline{\partial}
\varphi_{\alpha} }$. Comme $T$ est quasi positif, les fonctions
$\varphi_{\alpha}$ sont quasi plurisousharmoniques, donc
localement int\'egrables (c'est le seul point qui
diff\`ere du cas ${\cal C}^{\infty}$).

\noindent De l\`a, on \'ecrit successivement
$\varphi_{\alpha \beta} := \varphi_{\beta}-\varphi_{\alpha}$
sur l'intersection $U_{\alpha}\cap U_{\beta}$,
puis $$ i(\overline{\partial} - \partial)\varphi_{\alpha \beta}
= 2\pi d u_{\alpha \beta}.$$
Les fonctions $c_{\alpha \beta \gamma} := u_{\beta \gamma} - u_{\alpha \gamma}
+ u_{\alpha \beta}$ sont constantes. Comme $\{ T \}$ est
enti\`ere, la classe $\{ c_{\alpha \beta \gamma} \}$ l'est aussi.
Il existe donc une $1$-c\^ochaine \`a coefficients r\'eels
$\{ b_{\alpha \beta} \}$ telle que~:
$$ c_{\alpha \beta \gamma} + b_{\beta \gamma} - b_{\alpha \gamma} + b_{\alpha
\beta}
= m_{\alpha \beta \gamma} \in \ZZ.$$
On pose alors~:
$$ g_{\alpha \beta} :=
\exp (\varphi_{\alpha \beta} + 2i\pi(u_{\alpha \beta}
+ b_{\alpha \beta})).$$
Les $g_{\alpha \beta}$ sont holomorphes sans z\'ero et
forment un cocycle~: on note $E$ le fibr\'e en droites
associ\'e. Comme
$$ |g_{\alpha \beta}|\exp(-\varphi_{\alpha})=\exp(-\varphi_{\beta}),$$
les poids $\exp(-\varphi_{\alpha})$ d\'efinissent une
m\'etrique singuli\`ere sur $E$ dont le courant de courbure
est $T$.\finpreuve

\medskip

Comme nos in\'egalit\'es de Morse supposent que la
m\'etrique est \`a singularit\'es analytiques, nous avons besoin
d'un r\'esultat d'approximation permettant de
s'y ramener. De fa\c con g\'en\'erale, \'etant donn\'e
un courant ferm\'e $T$ sur une vari\'et\'e compacte,
c'est un probl\`eme classique que de vouloir le
r\'egulariser dans la m\^eme classe de cohomologie.
Malheureusement, il n'est pas possible en g\'en\'eral
de r\'egulariser ${\cal C}^{\infty}$ tout en perdant
aussi peu de positivit\'e que souhait\'e. L'obstruction
\`a le faire est mesur\'ee par les nombres de Lelong du
courant. Pour ne pas perdre de positivit\'e,
on r\'egularise en autorisant des singularit\'es
analytiques. Le r\'esultat que nous utilisons est le th\'eor\`eme
d'approximation des courants de J.-P.\ Demailly \cite{Dem92}~:

\bigskip

\noindent {\bf Th\'eor\`eme (J.-P.\ Demailly, 1992)}
{\em Soit $T$ un courant ferm\'e de bi-degr\'e $(1,1)$
de sorte que $T \geq \alpha$ o\`u $\alpha$ est une $(1,1)$
forme ${\cal C}^{\infty}$.
Alors pour toute m\'etrique hermitienne $\omega$
de classe ${\cal C}^{\infty}$ sur $X$, il existe
une suite de courants $T_{\varepsilon}$ telle que~:

(i) $\{ T_{\varepsilon} \} = \{ T \}$,

(ii) $T_{\varepsilon}$ tend (faiblement) vers $T$ lorsque
$\varepsilon$ tend vers $0$,

(iii) $T_{\varepsilon} \geq \alpha - \varepsilon \omega$,

(iv) $\displaystyle{ T_{\varepsilon}
= i \partial \overline{\partial} \varphi _{\varepsilon}
+ \beta }$,
o\`u $\varphi _{\varepsilon}$ est une fonction r\'eelle
\`a singularit\'es analytiques
et o\`u
$\beta$ est un repr\'esentant
$ {\cal C}^{\infty}$ de $\{ T \}$.
}

\bigskip

Il n'est pas question ici de donner la d\'emonstration
de ce r\'esultat, mais pour \'eclairer le lecteur,
mentionnons la premi\`ere \'etape de la d\'emonstration,
qui est une version locale du r\'esultat~:

\medskip

\noindent {\bf Proposition }
{\em Soit $\varphi$ une fonction psh dans
la boule unit\'e $B$ de $\CC ^n$. Pour $k$ entier
positif, notons ${\cal H}(k \varphi)$ l'espace
de Hilbert d\'efini de la fa\c con suivante~:
$$ {\cal H}(k \varphi) = \{ f \in {\cal O}(B) \ | \
\int _B |f|^2 \exp(-2k\varphi) d\lambda < +\infty \}.$$
Soit enfin $(\sigma_{j,k})_j$ une base orthonorm\'ee de
${\cal H}(k \varphi)$.
Alors la suite de fonctions
$$ \varphi _k := \frac{1}{2k} \log (\sum_j |\sigma_{j,k}|^2) $$
converge vers $\varphi$ simplement et dans $L_{\mbox{\scriptsize  loc}}^1$
lorsque $k$ tend vers $+\infty$.}

\medskip

La difficult\'e du th\'eor\`eme r\'eside dans le recollement
des diverses approximations locales donn\'ees par la proposition
pr\'ec\'edente. Les fonctions $\lambda _j$ figurant dans
la d\'efinition de l'hypoth\`ese {\cal S} proviennent
essentiellement de partition de l'unit\'e. En particulier,
elles ne s'annulent pas toutes simultan\'ement si bien que les
singularit\'es sont vraiment donn\'ees par les z\'eros communs
d'une famille de fonctions holomorphes.

\subsection{D\'emonstration des th\'eor\`emes C et D}
Nous d\'emontrons ici les deux caract\'erisations
analytiques. Pour cela, on commence par remarquer
que le sens faux dans le cadre des m\'etriques lisses
est v\'erifi\'e de fa\c con presque imm\'ediate dans le cadre plus souple
des vari\'et\'es de Moishezon.

Soit en effet $X$ une vari\'et\'e de Moishezon
et $\pi : \hat{X} \to X$ une modification projective.
Si $\hat{\omega}$ est une $(1,1)$ forme ${\cal C}^{\infty}$
d\'efinie positive sur $\hat{X}$ telle que $\{ \hat{\omega} \} \in
H^{2}(\hat{X},\ZZ)$, et si $\omega$ est une $(1,1)$ forme ${\cal C}^{\infty}$
d\'efinie positive sur $X$, alors la forme $\pi^{*}\omega$
est ${\cal C}^{\infty}$ et semi-positive. Il
existe donc une constante $A > 0$ telle que
$\hat{\omega} \geq A\pi^{*}\omega$.
Par cons\'equent,
le courant $T=\pi_{*} \hat{\omega}$ v\'erifie $T \geq A\omega$~:
en effet, si $\alpha$ est une $(n-1,n-1)$ forme positive
${\cal C}^{\infty}$ sur $X$, on a

\vspace{-3mm}
\begin{eqnarray*}
\ <T,\alpha> & = & \int_{\hat{X}} \hat{\omega} \wedge \pi^{*} \alpha \\
 & \geq & \int_{\hat{X}} A \pi^{*}\omega \wedge \pi^{*} \alpha \\
 & = &  A \int_{X} \omega \wedge \alpha \\
 & = & <A\omega,\alpha>.
\end{eqnarray*}

Le courant $T$ satisfait les points (i) et (ii) du th\'eor\`eme D,
et le th\'eor\`eme d'approximation des courants de Demailly
rappel\'e pr\'ec\'edemment implique qu'il existe
un courant $T' \in \{ T \}$
ayant localement les singularit\'es de l'hypoth\`ese {\cal S}
et
tel que $T' \geq  (A/2) \omega$. Ainsi, $T'$ v\'erifie la conclusion
du th\'eor\`eme C.

\medskip

Pour la r\'eciproque dans le th\'eor\`eme C, soit $T$
satisfaisant (i), (ii) et (iii). Alors, il existe un fibr\'e
en droites sur $X$ muni d'une m\'etrique
hermitienne singuli\`ere dont le courant de courbure est \'egal \`a
$T$. Les in\'egalit\'es de Morse singuli\`eres impliquent
que~:
$$ \dim H^0(X,E^k \otimes {\cal I}_k(h))
- \dim H^1(X,E^k \otimes {\cal I}_k(h)) \sim_{k \to +\infty} k^n.$$
{\em A fortiori},
$$ \dim H^0(X,E^k) \sim_{k \to +\infty} k^n ,$$
i.e le fibr\'e $E$ est gros et $X$ est de Moishezon.

Pour le th\'eor\`eme D, si $T$ est strictement positif,
on peut supposer par le th\'eor\`eme d'approximation des
courants que $T$ a localement les singularit\'es de l'hypoth\`ese {\cal S}.
On conclut alors comme ci-dessus.\finpreuve

\medskip

\noindent {\bf Remarque } La preuve du th\'eor\`eme D
donn\'ee par S.\ Ji et B.\ Shiffman suit la
d\'emarche suivante~: on approche, comme pr\'ec\'edemment,
le courant $T$ par un courant strictement positif
et singulier sur un ensemble analytique $S$. Puis
on montre directement gr\^ace aux estimations
$L^2$ de H\"ormander appliqu\'ees \`a la vari\'et\'e
k\"ahl\'erienne compl\`ete $X \backslash S$ que les grandes
puissances de $E$ engendrent les $1$-jets en dehors de $S$.

\subsection{Quelques commentaires}

Les th\'eor\`emes C et D sont satisfaisants car ils donnent
une caract\'erisation analytique des vari\'et\'es de Moishezon.
Cependant, on souhaiterait pouvoir se dispenser de
l'hypoth\`ese faite sur les singularit\'es
dans le th\'eor\`eme C ou affaiblir l'hypoth\`ese de stricte
positivit\'e du th\'eor\`eme D. On conjecture le
r\'esultat suivant~:

\medskip

\noindent {\bf Conjecture }
{\em Une vari\'et\'e compacte $X$
est de Moishezon si et
seulement si il existe sur $X$ un courant $T$ positif ferm\'e
de bi-degr\'e $(1,1)$
tel que~:

(i) $ \displaystyle{\{ T \} \in H^{2}(X,\ZZ) }$,

(ii) il existe un ouvert $U$ sur lequel $T$ est
strictement positif.}

\medskip

La condition (ii) est une fa\c con d'assurer que le support
de $T$ n'est pas n\'egligeable pour la mesure de Lebesgue.
Evidemment, l'hypoth\`ese {\cal S} nous a permis de
d\'efinir les int\'egrales de courbure en int\'egrant
simplement sur la partie lisse de la m\'etrique. Dans le cas
g\'en\'eral, d\'efinir un produit de courants $T^n$ est un probl\`eme
de type {\bf Monge-Amp\`ere}.

\medskip

Le r\'esultat suivant va dans le sens de la
conjecture~:

\medskip

\noindent {\bf Proposition }
{\em La conjecture est vraie dans le cas des
surfaces complexes (et dans ce cas, $X$
\'etant de Moishezon est donc projective).}

\medskip

\noindent {\bf D\'emonstration }

Soit $\omega$ une m\'etrique hermitienne sur $X$,
et pour $\varepsilon > 0$, soit $T_{\varepsilon}$
un courant donn\'e par le th\'eor\`eme d'approximation
des courants v\'erifiant $T_{\varepsilon} \geq - \varepsilon \omega$.
Par le th\'eor\`eme C, il suffit de montrer
que pour $\varepsilon$ assez petit, on a
$$ \displaystyle{\int_{X(\leq 1,T_{\varepsilon})} T_{\varepsilon}^{n} > 0}.$$
Comme il existe un ouvert $U$ sur lequel $T$ est
strictement positif (disons sup\'erieur \`a $C_U \omega _{| U}$
o\`u $C_U$ est une constante strictement positive),
il existe un ouvert plus petit $U'$ ind\'ependant de $\varepsilon$
sur lequel $T_{\varepsilon}$ est sup\'erieur \`a $(C_{U}/2) \omega _{| U'}$.
On en d\'eduit qu'il existe une constante $C>0$, ind\'ependante de
$\varepsilon$
telle que pour tout $\varepsilon$ petit, on ait
$$\displaystyle{\int_{X(0,T_{\varepsilon})} T_{\varepsilon}^{n}
\geq \int_{U'} T_{\varepsilon}^{n} > C}.$$
Il suffit donc de montrer que
$$ \lim _{\varepsilon \to 0}\int_{X(1,T_{\varepsilon})} T_{\varepsilon}^{n}
=0.$$
Or, sur l'ouvert $X(q,T_{\varepsilon})$, on a
$$ 0 \leq (-1)^q T_{\varepsilon}^{n} \leq \frac{n!}{q! (n-q)!}
\varepsilon ^q \omega ^q \wedge
(T_{\varepsilon} + \varepsilon \omega )^{n-q}.$$
Il suffit donc de contr\^oler les produits de
Monge-Amp\`ere et plus pr\'ecis\'ement de montrer que pour $q >0$
(et dans la perspective de la conjecture, $q=1$ suffit), on a~:

$$ \lim _{\varepsilon \to 0} \varepsilon ^q \int_X \omega ^q \wedge
(T_{\varepsilon} + \varepsilon \omega )^{n-q} =0.$$
C'est \'evidemment vrai pour $q=n$ et, si on suppose
de plus que la m\'etrique $\omega$ est une m\'etrique de Gauduchon,
alors la formule de Stokes donne
$$ \int_X \omega ^{n-1} \wedge
(T_{\varepsilon} + \varepsilon \omega ) = \mbox{Cste} + O(\varepsilon),$$
et donc c'est aussi vrai pour $q = n-1$.

Dans le cas des surfaces, on a l'estimation souhait\'ee
pour $q= n-1 =1$, donc $X$ est de Moishezon par le
th\'eor\`eme C.\finpreuve

\medskip

Comme le montre la preuve ci-dessus, le cas g\'en\'eral
pourrait d\'ecouler du contr\^ole des masses de Monge-Amp\`ere
des approximations $\varphi _{\varepsilon}$ utilis\'ees par
J.-P.\ Demailly dans la d\'emonstration de son th\'eor\`eme
d'approximation des courants.

\section{Une version alg\'ebrique singuli\`ere des in\'egalit\'es
de Morse}

Dans ce paragraphe, nous donnons quelques exemples ``alg\'e\-bri\-ques"
de faisceaux d'id\'eaux de Nadel. Leur origine en g\'eom\'etrie
alg\'ebrique se situe dans la version du th\'eor\`eme de Kawamata-Viehweg
pour les diviseurs \`a coefficients rationnels. Comme dans \cite{Dem90} ou
\cite{EsV92}, ce faisceau d'id\'eaux sert de terme correctif dans le cas o\`u
les diviseurs consid\'er\'es ne sont pas \`a croisements normaux.

Par ailleurs, il existe une version alg\'ebrique
des in\'egalit\'es
de Morse holomorphes de Demailly \cite{Dem94} dans le cas d'un
fibr\'e en droites
diff\'erence de deux fibr\'es amples. Il est alors naturel de donner
une version analogue dans le cadre singulier. \'Evidemment, il s'agit
essentiellement d'une reformulation dans un cas particulier de
notre theor\`eme~B.

\subsection{Th\'eor\`eme de Kawamata-Viehweg}
Soient $X$ une vari\'et\'e projective, et $M$ un diviseur rationnel
effectif de $X$. On note $M = \sum a_i D_i$ o\`u les
$a_i$ sont des rationnels positifs et les $D_i$ sont des diviseurs
irr\'eductibles. Nous notons ${\cal I}(M)$
le faisceau d'id\'eaux de Nadel associ\'e \`a la m\'etrique singuli\`ere
$\phi = \sum a_i \log|g_i|$
o\`u $g_i$ est un g\'en\'erateur local de $D_i$.
On rappelle que si $M$ est
\`a croisements normaux, alors ${\cal I}(M)$ n'est rien d'autre
que le faisceau
inversible ${\cal O}(- \lfloor M \rfloor)$ o\`u
$\lfloor M \rfloor := \sum  \lfloor a_i \rfloor D_i$.
Avec ces notations, rappelons le th\'eor\`eme
de Kawamata-Viehweg \cite{Kaw82}, \cite{Vie82}~:

\medskip

\noindent {\bf Th\'eor\`eme d'annulation de Kawamata-Viehweg  } {\em
Soient $X$ une vari\'et\'e projective et
$L$ un fibr\'e en droites sur $X$.
On suppose que $L = M + \sum_{j} \alpha _jE_j$
o\`u~:

(i) $M$ est un $\QQ$-diviseur effectif gros et nef,

(ii) les $\alpha _j$ sont des r\'eels v\'erifiant
$0 \leq \alpha _j <1$,

(iii) le diviseur $\sum_{j} \alpha _jE_j$ est
\`a croisements normaux.

\noindent Alors
$$ H^q(X,K_X + L)=0 \ \mbox{pour tout}\ q \geq 1 .$$}

Dans le cas o\`u $\sum_{j} \alpha _jE_j$
n'est pas \`a croisements
normaux dans l'\'enonc\'e du th\'eor\`eme de Kawamata-Viehweg,
le faisceau multiplicateur de Nadel sert de terme correctif.
Nous allons d\'etailler un exemple utilis\'e par
L.\ Ein et R.\ Lazarsfeld dans l'\'etude
des diviseurs \`a singularit\'e presque isol\'ee (voir
le papier de R.\ Lazarsfeld \cite{Laz93} pour plus de d\'etails).

Le contexte est le suivant : soient $X$ une vari\'et\'e projective,
$A$ un fibr\'e gros et nef sur $X$ et $D$ un diviseur dans le syst\`eme
lin\'eaire $|kA|$.

Choisissons $\pi : \tilde{X} \rightarrow X$ une compos\'ee d'un nombre fini
d'\'eclate\-ments de centres lisses de sorte que $\pi ^* D$ soit
un diviseur
\`a croisements normaux et soit $K_{\tilde{X}/X}$ la diff\'erence
des fibr\'es canoniques~:
$$K_{\tilde{X}/X} = K_{\tilde{X}} - \pi^*K_{X}.$$
Alors, pour $\lambda$ rationnel, $0\leq \lambda <1$, on pose~:
$${\cal I}_{\lambda}=
\pi_*\left(K_{\tilde{X}/X}-
\lfloor \pi^*(\lambda \frac{D}{k}) \rfloor \right).$$

La proposition suivante donne les propri\'et\'es de ${\cal I}_{\lambda}$~:

\medskip

\noindent {\bf Proposition }{\em On a

(i) le faisceau d'id\'eaux $\displaystyle{{\cal I}_{\lambda}}$
est \'egal au faisceau multiplicateur
$\displaystyle{{\cal I}(\lambda \frac{D}{k})}$, en particulier
${\cal I}_{\lambda}$ est ind\'ependant de la r\'esolution $\pi$ choisie,

(ii) les images directes sup\'erieures
$\displaystyle{R^q \pi_*\left(K_{\tilde{X}/X}-
 \lfloor \pi^*(\lambda \frac{D}{k}) \rfloor \right)}$ sont
nulles pour tout
$q \geq 1$,

(iii) les groupes de cohomologie
$\displaystyle{H^q(X,(K_X + A)\otimes{\cal I}_{\lambda})}$ sont nuls pour tout
$q \geq 1$.
}

\medskip

\noindent {\bf D\'emonstration}

Pour le point (i), nous utilisons \`a nouveau l'identit\'e
$$\mu _*(K_{\tilde{X}}\otimes {\cal I}(\mu^*h))
= K_X \otimes {\cal I}(h).$$
Dans notre situation, on a
\vspace{-2mm}
$$ \pi_*\left(K_{\tilde{X}}\otimes
{\cal I}(\pi^*(\lambda \frac{D}{k}))\right)
= K_X \otimes {\cal I}(\lambda \frac{D}{k}).$$
\vspace{-3mm}
Comme $\displaystyle{\pi^*(\lambda \frac{D}{k})}$ est \`a croisements normaux,
on a
\vspace{-2mm}
$$ {\cal I}(\pi^*(\lambda \frac{D}{k}))=
{\cal O}\left(- \lfloor \pi^*(\lambda \frac{D}{k}) \rfloor \right).$$
De l\`a, on d\'eduit successivement
\vspace{-3mm}
\begin{eqnarray*}
{\cal I}_{\lambda} & = &
\pi_*\left(K_{\tilde{X}/X}-
\lfloor \pi^*(\lambda \frac{D}{k}) \rfloor \right)\\
& = & \pi_*\left(K_{\tilde{X}}-\pi^*K_X -
\lfloor \pi^*(\lambda \frac{D}{k}) \rfloor \right)\\
& = & -K_X \otimes K_X \otimes {\cal I}(\lambda \frac{D}{k}) \\
& = & {\cal I}(\lambda \frac{D}{k}).
\end{eqnarray*}
\noindent Ceci d\'emontre le point (i).

Pour le point (ii), la d\'emonstration repose sur l'observation
classique suivante (d\'ej\`a observ\'ee par H.\ Grauert et
O.\ Riemenschneider \cite{GrR70})~: les faisceaux $R^q\pi_*{\cal F}$ sont nuls
si et
seulement si pour tout fibr\'e $L$ sur $X$ suffisamment ample
et tout $q \geq 1$,
on a $H^q(\tilde{X},{\cal F}\otimes \pi^*L) = 0$.
Comme
$$ R^q \pi_* \left(\pi^*(K_X + \pi^*A) +
K_{\tilde{X}/X}- \lfloor \pi^*(\lambda \frac{D}{k}) \rfloor \right) =$$
$$(K_X + \pi^*A)\otimes  R^q \pi_* \left(
 K_{\tilde{X}/X}- \lfloor \pi^*(\lambda \frac{D}{k}) \rfloor \right),$$
il suffit de montrer que pour tout fibr\'e $L$ sur $X$ suffisamment ample
et tout $q \geq 1$, on a
$$ H^q\left(\tilde{X},\pi^*L + \pi^*A +  K_{\tilde{X}}-
\lfloor \pi^*(\lambda \frac{D}{k})\rfloor \right) = 0. $$
Or ce dernier groupe est \'egal \`a
$$ H^q\left(\tilde{X}, K_{\tilde{X}} + \pi^*L + \left(1-\lambda\right)\pi^*A
+ \left(\lambda\pi^*A - \lfloor \lambda\pi^*A\right)\rfloor \right)$$
qui est nul gr\^ace au th\'eor\`eme de Kawamata-Viehweg
appliqu\'e au $\QQ$-diviseur effectif gros et nef
$$ M:= \pi^*L + (1-\lambda)\pi^*A .$$

Pour le point (iii), on a
$$ H^q\left(X,\left(K_X + A\right)\otimes{\cal I}_{\lambda}\right)=
H^q\left(\tilde{X}, K_{\tilde{X}} + \pi^*A-\lfloor \pi^*(\lambda
\frac{D}{k})\rfloor \right).$$
Ce dernier groupe s'\'ecrit encore
$\displaystyle{H^q\left(\tilde{X}, K_{\tilde{X}} +\left(1-\lambda\right)\pi^*A
+ \lambda\pi^*A - \lfloor \lambda\pi^*A \rfloor \right)}$
qui est nul \`a nouveau par le th\'eor\`eme de Kawamata-Viehweg.\finpreuve

\subsection{In\'egalit\'es de Morse alg\'ebriques singuli\`eres}

En nous inspirant de l'exemple pr\'ec\'edent, nous sommes en mesure
de donner une version alg\'ebrique des in\'egalit\'es de Morse
holomorphes singuli\`eres. Pour cela, rappelons au pr\'ea\-la\-ble la version
suivante des in\'egalit\'es de Morse
holomorphes de J.-P.\ Demailly \cite{Dem94}~:

\bigskip

\noindent {\bf Th\'eor\`eme }
{\em Soit $X$ une vari\'et\'e k\"ahl\'erienne de dimension $n$ et soient
$F$ et $G$ deux fibr\'es en droites nef sur $X$.
Alors, on a~:
$$ \sum _{j=0}^{q}(-1)^{q-j} \dim H^{j}(X,k(F-G))
\leq \frac{k^{n}}{n!}\sum _{j=0}^{q}(-1)^{q-j}{n \choose j} F^{n-j}\cdot G^j
+ o(k^{n}).$$
}

\bigskip

Ce r\'esultat a \'et\'e obtenu dans un premier temps par
J.-P.\ Demailly comme
cons\'equence des in\'egalit\'es de Morse, et plus r\'ecemment,
F.\ Angelini en a donn\'e une d\'emonstration purement alg\'ebrique
\cite{Ang95}. Auparavant, S.\ Trapani \cite{Tra91} et Y.-T.\ Siu
\cite{Siu93} avaient d\'emontr\'e le cas particulier
$q = 1$ en vue d'obtenir des crit\`eres num\'eriques pour
l'existence de sections. Le terme alg\'ebrique a ici une double origine~:
les estimations font intervenir des nombres d'intersection
\`a la place d'int\'egrales de courbure, et un cas particulier
du th\'eor\`eme est celui d'une vari\'et\'e projective et
d'un fibr\'e \'ecrit comme diff\'erence de deux fibr\'es amples.

Nous montrons le r\'esultat suivant~:

\bigskip

\noindent {\bf Th\'eor\`eme E } {\em
Soit $X$ une vari\'et\'e k\"ahl\'erienne de dimension $n$ et soient
$F$ et $G$ deux fibr\'es en droites sur $X$.
On suppose que $G$ est nef, et qu'il existe un entier
positif $m$, un fibr\'e en droites nef $A$ et un diviseur effectif
$D$ de sorte que~: $mF = A + D.$
Alors, on a~:
$$ \sum _{j=0}^{q}(-1)^{q-j}
\dim H^{j}((X,k(F-G)\otimes {\cal I}_{k}(m^{-1}D)) $$
$$ \leq
\frac{k^{n}}{n!}
\sum _{j=0}^{q}(-1)^{q-j}{n \choose j} m^{-n+j}A^{n-j}\cdot G^j
+ o(k^{n}).$$
}

\bigskip

\noindent {\bf D\'emonstration}

Soit $\pi : \tilde{X} \rightarrow X$ une
compos\'ee d'un nombre fini
d'\'eclate\-ments de centres lisses de sorte que $\pi ^* D$ soit un diviseur
\`a croisements normaux.

\noindent Comme dans le cadre purement analytique de nos in\'egalit\'es
de Morse singuli\`eres, on travaille sur
$\tilde{X}$ o\`u l'on applique simplement
l'\'enonc\'e pr\'ec\'edent.

\noindent D\'etaillons bri\`evement~: le faisceau d'id\'eaux
${\cal I}_{k}(m^{-1}D)$
est \'egal
\`a l'image directe $$\pi_*\left(K_{\tilde{X}/X}-\lfloor
\pi^*(km^{-1}D)\rfloor \right),$$
si bien qu'il s'agit d'estimer les dimensions
$$\dim H^{q}\left(\tilde{X},k\pi^*(F-G)-\lfloor
\pi^*(km^{-1}D)\rfloor\right).$$

\noindent Or pour $k$ multiple de $m$, on a
$$k\pi^*(F-G)-\lfloor \pi^*(km^{-1}D)\rfloor = k\pi^*(m^{-1}A-G),$$
\noindent et il suffit d'appliquer le th\'eor\`eme pr\'ec\'edent
aux fibr\'es nef $\pi^*A$ et $\pi^*G$.\finpreuve


\chapter{\'Etude de certaines vari\'et\'es
de Moishezon dont le groupe de Picard est
infini cyclique}

Le th\`eme central de ce chapitre
est l'\'etude d'une classe particuli\`ere
de vari\'et\'es de Moishezon~: celles dont
le groupe de Picard est $\ZZ$ et dont le fibr\'e
canonique est gros. En faisant l'hypoth\`ese suppl\'ementaire
que la vari\'et\'e $X$ devient projective apr\`es un seul
\'eclatement de centre lisse et projectif, nous \'etudions ce centre
{\em via} la th\'eorie de Mori
sur le mod\`ele projectif. Nous obtenons alors
une restriction sur la dimension du centre de l'\'eclatement
dans le cas o\`u le fibr\'e canonique n'est pas nef.
Apr\`es avoir donn\'e une nouvelle famille de vari\'et\'es
de Moishezon ne poss\'edant pas de fibr\'e en droites
gros et nef et s'inscrivant
dans ce cadre d'\'etude, nous nous restreignons \`a
la dimension $4$. Nous obtenons alors une description
pr\'ecise du centre de l'\'eclatement et montrons que
notre construction est essentiellement unique dans le cas
o\`u le fibr\'e canonique n'est pas nef. Enfin, nous
obtenons aussi des restrictions partielles en dimension $4$
dans le cas o\`u le fibr\'e canonique est nef.

\section{Un th\'eor\`eme de J.\ Koll\'ar}

\subsection{\'Enonc\'e du r\'esultat}

Nous avons vu dans les pr\'eliminaires de cette
th\`ese qu'il existe des vari\'et\'es de Moishezon
ne poss\'edant pas de fibr\'e en droites
simultan\'ement gros et num\'eriquement effectif.
Plus pr\'ecis\'ement, nous avons rencontr\'e le
r\'esultat suivant~:

\bigskip

\noindent {\bf Th\'eor\`eme (J.\ Koll\'ar, K.\ Oguiso) } {\em

(i) Il existe des vari\'et\'es de Moishezon $X$ de
dimension
$3$ dont le groupe de Picard est \'egal \`a $\ZZ$,
avec $-K_X$ gros et ne poss\'edant
pas de fibr\'e gros et nef,

(ii) il existe des vari\'et\'es de Moishezon $X$ de dimension
$3$ dont le groupe de Picard est \'egal \`a $\ZZ$,
dont le fibr\'e canonique est trivial et
ne poss\'edant pas de fibr\'e gros et nef.
}

\bigskip

Remarquons que pour une vari\'et\'e de Moishezon dont le
groupe de Picard est $\ZZ$,
un et un seul g\'en\'erateur de $\Pic(X)$ est gros. Suivant
J.\ Koll\'ar \cite{Kol91}, nous notons ce g\'en\'erateur
${\cal O}_X(1)$ et nous \'ecrivons $\Pic(X)
= \ZZ  \cdot {\cal O}_X(1)$.

Notons  aussi $m_X$ l'entier d\'efini par la relation $K_X = {\cal O}_X(m_X)$.
Remarquons ici que les trois cas $m_X < 0$
(respectivement $m_X=0$ et $m_X >0$)
correspondent aux trois possibilit\'es $\kappa(X) = -\infty$ (respectivement
$\kappa(X) = 0$ et $\kappa(X) = \dim X$), o\`u  $\kappa(X)$ d\'esigne la
dimension de Kodaira de $X$.

\medskip

Evidemment, il reste un cas non couvert par l'\'enonc\'e pr\'ec\'edent
et la question suivante est naturelle~:

\medskip

\noindent {\bf Question } {\em  Existe-t-il des vari\'et\'es de Moishezon $X$,
avec $\Pic(X) = \ZZ \cdot {\cal O}_X(1)$ et $m_X > 0$ ne poss\'edant pas
de fibr\'e gros et nef \  ?
}

\medskip

En dimension $3$, la r\'eponse \`a cette question est n\'egative
comme le montre le r\'esultat suivant~:

\bigskip

\noindent {\bf Th\'eor\`eme (J.\ Koll\'ar, 1991) } {\em
Soit $X$ une vari\'et\'e de Moishezon
de dimension $3$. On suppose que le groupe
de Picard $\Pic(X)$ est $\ZZ$ et que le
fibr\'e canonique $K_X$ est gros.
Alors $K_X$ est nef.}

\bigskip

\noindent {\bf Remarque  } Dans le cas o\`u le fibr\'e canonique
est gros et nef, l'un de ses multiples est globalement engendr\'e
par le ``base-point free theorem".
Au vu du r\'esultat pr\'ec\'edent, il n'existe
donc pas d'exemple
de vari\'et\'e de Moishezon de dimension $3$,
de groupe de Picard $\ZZ$ et \`a fibr\'e canonique gros
ne satisfaisant pas aux crit\`eres de J.-P.\ Demailly et Y.-T.\ Siu.

\medskip

La suite de ce paragraphe consiste \`a rappeler
la d\'emonstration de ce r\'esultat car les id\'ees
qu'elle contient seront pr\'esentes dans tout
le chapitre.

\subsection{D\'emonstration}

La r\'ef\'erence pr\'ecise est \cite{Kol91},
page 170 et suivantes.

\medskip

Le lemme suivant, bien qu'\'el\'ementaire est essentiel~:

\medskip

\noindent {\bf Lemme } {\em
Soit $X$ une vari\'et\'e de Moishezon, de dimension quelconque,
avec $\Pic(X) = \ZZ \cdot {\cal O}_X(1)$. Soit
$\pi : \tilde{X} \rightarrow X$ une modification projective,
de
lieu exceptionnel
$\tilde{S} \subset \tilde{X} \stackrel{\pi}{\rightarrow} S \subset X$.

\noindent Alors pour toute courbe $C$ de $X$, non incluse dans $S$, on a
${\cal O}_X(1) \cdot C > 0$.}

\bigskip

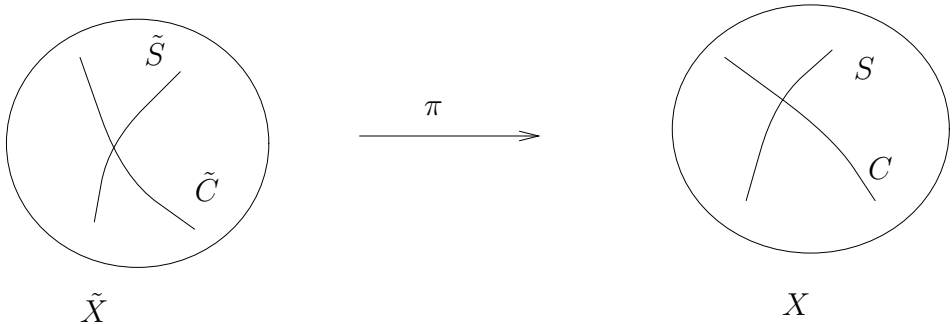
\setlength{\unitlength}{0.0125in}
\begin{picture}(395,148)(0,-10)
\put(337,81){\ellipse{116}{104}}
\path(215.000,76.000)(223.000,78.000)(215.000,80.000)
\path(223,78)(148,78)(148,78)
\spline(31,111)
(49,60)(79,39)
\spline(37,42)
(43,75)(73,105)
\spline(301,111)
(346,78)(364,51)
\spline(310,51)
(322,93)(346,114)
\put(55,75){\ellipse{110}{104}}
\put(79,51){\makebox(0,0)[lb]{\raisebox{0pt}[0pt][0pt]{\shortstack[l]{{\twlrm
$\tilde{C}$}}}}}
\put(325,3){\makebox(0,0)[lb]{\raisebox{0pt}[0pt][0pt]{\shortstack[l]{{\twlrm
$X$}}}}}
\put(58,108){\makebox(0,0)[lb]{\raisebox{0pt}[0pt][0pt]{\shortstack[l]{{\twlrm
$\tilde{S}$}}}}}
\put(31,0){\makebox(0,0)[lb]{\raisebox{0pt}[0pt][0pt]{\shortstack[l]{{\twlrm
$\tilde{X}$}}}}}
\put(361,60){\makebox(0,0)[lb]{\raisebox{0pt}[0pt][0pt]{\shortstack[l]{{\twlrm
$C$}}}}}
\put(355,102){\makebox(0,0)[lb]{\raisebox{0pt}[0pt][0pt]{\shortstack[l]{{\twlrm
$S$}}}}}
\put(175,87){\makebox(0,0)[lb]{\raisebox{0pt}[0pt][0pt]{\shortstack[l]{{\twlrm
$\pi$}}}}}
\end{picture}


\bigskip

\noindent {\bf Corollaire } {\em
Sous les hypoth\`eses du lemme, et si de plus $K_X$ est gros,
alors pour toute courbe $C$ (respectivement $\tilde{C}$)
de $X$ (respectivement de $\tilde{X}$) non incluse dans $S$
(respectivement $\tilde{S}$), on a $K_X\cdot C > 0$ (respectivement
$K_{\tilde{X}} \cdot \tilde{C} > 0$).}

\medskip

\noindent {\bf D\'emonstration du lemme }

Soit $\tilde{H}$ un diviseur ample dans $\tilde{X}$, et
$H := \pi _{*}(\tilde{H})$. Alors $H$ est gros, donc s'\'ecrit
${\cal O}_X(p)$, o\`u $p$ est un entier strictement positif.
Comme
$${\cal O}_X(1) \cdot C = \frac{1}{p}H \cdot C ,$$
il suffit de montrer que pour toute courbe $C$ non incluse
dans $S$, on a $H \cdot C > 0$.

\noindent Comme $C$ n'est pas incluse dans $S$, si $\tilde{C}$ d\'esigne
la transform\'ee stricte de $C$, l'\'egalit\'e suivante est v\'erifi\'ee~:
$$ H \cdot C = \pi^{*}(H) \cdot \tilde{C} =
(\tilde{H} + \sum a_iE_i) \cdot \tilde{C},$$
o\`u les $a_i$ sont des entiers positifs ou nuls, et les $E_i$ les composantes
irr\'eductibles de $\tilde{S}$. De l\`a, $H \cdot C > 0$
car $\tilde{H} \cdot \tilde{C} > 0$ et pour tout $i$, on a $E_i \cdot \tilde{C}
\geq 0$.
\finpreuve

\medskip

\noindent {\bf Remarque } Ce lemme affirme en particulier que les courbes
sur lesquelles ${\cal O}_X(1)$ est
n\'egatif ou nul sont incluses dans un ensemble analytique
de codimension sup\'erieure ou \'egale \`a $2$.

\medskip

\noindent {\bf D\'emonstration du th\'eor\`eme }

En dimension $3$, on d\'eduit du lemme pr\'ec\'edent qu'il n'y
a qu'un nombre fini de courbes sur lesquelles ${\cal O}_X(1)$ est
n\'egatif. Si $K_X$ est gros, et si $C$ est une courbe telle que
$K_X \cdot C < 0$, alors une telle courbe se d\'eforme dans $X$ car la formule
de Riemann-Roch donne
$$ \chi (N_{C/X}) = -K_X \cdot C + (n-3)(1-g) = -K_X \cdot C > 0,$$
o\`u $N_{C/X}$ est le fibr\'e normal de $C$ dans $X$. Mentionnons
que si la courbe $C$ est singuli\`ere, on d\'efinit
$N_{C/X}$ comme \'etant \'egal \`a $\displaystyle{ \nu ^{\ast}T X / T \tilde{C}
}$
o\`u $\nu : \tilde{C} \to C$ est la normalis\'ee de la courbe $C$.
Ceci donne bien la contradiction.\finpreuve

\subsection{Commentaires}

La d\'emonstration ci-dessus repose sur un argument
de d\'eformation. Comme nous utilisons dans la suite
ce type d'argument, il est sans doute bon de faire ici
un bref rappel.

\'Etant donn\'ees une vari\'et\'e $X$ et une sous-vari\'et\'e
$Y$ de $X$, c'est un probl\`eme classique et important de
d\'eterminer les d\'eformations de $Y$ dans $X$. Dans le
cadre analytique, ce probl\`eme a \'et\'e consid\'er\'e par
K.\ Kodaira \cite{Kod62} et par A.\ Grothendieck et
D.\ Mumford dans le cadre
alg\'ebrique o\`u la notion de {\bf sch\'ema de Hilbert}
joue un r\^ole essentiel~: une r\'ef\'erence importante est
le travail r\'ecent de J.\ Koll\'ar \cite{Kol94}.

La ``solution" au probl\`eme est donn\'ee par le~:

\bigskip

\noindent {\bf Th\'eor\`eme \cite{Gro62}, \cite{Kod62}  }{\em
Soit $Y$ une sous-vari\'et\'e d'une vari\'et\'e $X$.
Alors le sch\'ema de Hilbert $\Hilb (X)$ des sous-ensembles
analytiques de $X$ admet $H^0(Y,N_{Y/X})$ comme espace tangent
de Zariski en $[Y]$. La dimension de $\Hilb (X)$ en $[Y]$ satisfait~:
$$ \dim H^0(Y,N_{Y/X}) - \dim H^1(Y,N_{Y/X})
\leq \dim_{[Y]} \Hilb (X) \leq \dim H^0(Y,N_{Y/X}).$$
En particulier, si $H^1(Y,N_{Y/X}) =0$, alors $\Hilb (X)$ est lisse
au voisinage de $[Y]$.}

\bigskip

Dans le cadre analytique, la construction de K.\ Kodaira
consiste \`a trouver explicitement en coordonn\'ees
locales les s\'eries enti\`eres d\'efinissant les sous-ensembles
proches de $Y$, tandis que dans le cadre alg\'ebrique, l'id\'ee
est qu'une vari\'et\'e projective $Z$ est d\'etermin\'ee par
le sous-espace vectoriel des polyn\^omes de degr\'e suffisamment
grand qui s'annulent sur $Z$.

\section{Quelques rappels sur la th\'eorie de Mori}

La th\'eorie de Mori, n\'ee dans les ann\'ees 80,
consiste \`a \'etendre et approfondir en dimension sup\'erieure
ou \'egale \`a $3$ la classification des surfaces
complexes et l'\'etude des applications bim\'eromorphes
entre surfaces complexes.

Cependant, cette
th\'eorie n'est valable que sur les vari\'et\'es
projectives. Comme une vari\'et\'e de Moishezon
est domin\'ee par une vari\'et\'e projective, une
id\'ee naturelle est d'appliquer certains r\'esultats
de la th\'eorie de Mori pour obtenir des renseignements
concernant la structure des vari\'et\'es de Moishezon
ou de leur caract\`ere non projectif.

Nous rappelons dans ce paragraphe les r\'esultats
essentiels utilis\'es dans ce chapitre. Deux excellentes
r\'ef\'erences sont \cite{CKM88} et \cite{KMM87}.

Mentionnons aussi qu'une des grandes id\'ees de S.\ Mori
est, m\^eme pour l'\'etude des vari\'et\'es
non singuli\`eres, de quitter le monde lisse pour autoriser
certains types de singularit\'es~; cette analyse
fut en particulier mise en oeuvre par M.\ Reid. Il est cependant
bon de pr\'eciser ici que la plupart
des \'enonc\'es que nous utilisons sont d'une difficult\'e
moindre dans le cas lisse, cadre dans lequel nous les
appliquons.

\subsection{C\^one des courbes effectives}
Dans tout ce paragraphe, $X$ est une vari\'et\'e projective.

\subsub{Notations}

Rappelons tout d'abord que la notation
$N_1(X,\RR)$ d\'esigne l'espace vectoriel
des combinaisons lin\'eaires finies (\`a coefficients
r\'eels) de courbes (irr\'eductibles et
\'eventuellement singuli\`eres) de $X$, modulo l'\'equivalence
num\'erique~: deux courbes sont \'equivalentes si et
seulement si leurs intersections avec tout diviseur
sont \'egales.

\noindent Pour une vari\'et\'e projective (et m\^eme de
Moishezon), cet espace vectoriel est de dimension finie et est
en dualit\'e naturelle ({\em via} la forme d'intersection)
avec le groupe de N\'eron-Severi
$\displaystyle{(\Pic (X)/ \Pic ^0(X))\otimes _{\ZZ} \RR}$~;
la dimension de $N_1(X,\RR)$
est appel\'ee {\bf nombre de Picard de $X$}.
L'espace vectoriel
$N_1(X,\RR)$ est naturellement un sous-espace vectoriel
de $H_2(X,\RR)$.

Enfin, nous notons suivant l'usage $\NE (X)$ le
sous-c\^one convexe de $N_1(X,\RR)$ engendr\'e par les classes
d'homologie des courbes
effectives. L'adh\'erence de ce c\^one est not\'ee
$\overline{\NE }(X)$.

\subsub{Th\'eor\`eme du c\^one}

L'un des premiers succ\`es de la th\'eorie de Mori est de
montrer que {\em si le fibr\'e canonique $K_X$
n'est pas nef, alors il existe une courbe rationnelle
sur laquelle il est strictement n\'egatif}.

L'\'enonc\'e pr\'ecis est le suivant~:

\medskip

\noindent {\bf Th\'eor\`eme du c\^one } {\em
Soit $X$ une vari\'et\'e projective. Alors il existe
un ensemble minimal (fini ou d\'enombrable) de courbes
rationnelles $C_i$ dans $X$ de sorte que~:

\vspace{+1mm}
(i) pour tout $i$, on a $\displaystyle{0 < -K_X \cdot C_i \leq \dim X +1}$,

\vspace{+2mm}
(ii) $\displaystyle{
\overline{\NE }(X) = \overline{\NE }(X)_{K_X \geq 0} + \sum_i \RR_+ \ [C_i]}$,

\noindent o\`u
$\displaystyle{
\overline{\NE }(X)_{K_X \geq 0} := \{ [C] \in N_1(X,\RR) \ |
\ K_X \cdot C \geq 0 \} }$.
}

\medskip

Les courbes rationnelles $C_i$ sont appel\'es
{\bf courbes rationnelles extr\^emales} et les $\RR_+ \ [C_i]$
sont appel\'ees {\bf rayons extr\^emaux}.

On ne connait pas de version analogue de ce th\'eor\`eme
pour les vari\'et\'es de Moishezon. \`A notre connaissance,
la question suivante est ouverte~:

\medskip

\noindent {\bf Question } {\em Soit $X$ une
vari\'et\'e de Moishezon, dont le fibr\'e canonique
n'est pas nef. Existe-t-il alors une courbe
rationnelle $C$ telle que $K_X \cdot C < 0$ \ ?
}

\subsection{Contraction de Mori}

Les rayons extr\^emaux jouent le m\^eme r\^ole
dans la th\'eorie de Mori que les courbes rationnelles
lisses d'auto-intersection n\'egative dans la th\'eorie
des surfaces complexes~: ils peuvent \^etre contract\'es.
Le th\'eor\`eme suivant d\'ecrit les diff\'erents types
de contraction obtenus en contractant un rayon extr\^emal.

\bigskip

\noindent {\bf Th\'eor\`eme de contraction \cite{CKM88} }
{\em Soit $X$ une vari\'et\'e projective
dont le fibr\'e canonique n'est pas nef. Soient
$C$ une courbe rationnelle extr\^emale et $R := \RR_+ \ [C]$
le rayon extr\^emal engendr\'e par $C$.

Alors, il existe une vari\'et\'e (\'eventuellement singuli\`ere)
projective,
normale, et une application $$f : X \to Y ,$$ not\'ee
aussi $\con_R$, de sorte que~:

(i) une courbe de $X$ est contract\'ee
par $f$ si et seulement si sa classe d'homologie
appartient au rayon $R$,

(ii) le fibr\'e $-K_X$ est $f$-ample (i.e la restriction
de $-K_X$ \`a toute fibre de $f$ est ample).

\noindent De plus, on distingue trois types de contractions~:

(a) $\dim X > \dim Y$ et $f$ est une fibration
Fano (i.e la fibre g\'en\'erique de $f$ est une
vari\'et\'e lisse dont le fibr\'e anti-canonique est ample),

(b) $\dim X = \dim Y$ et $f$ est une contraction divisorielle
(i.e $f$ est birationnelle et contracte un diviseur),

(c) $\dim X = \dim Y$ et $f$ est une petite contraction
(i.e $f$ est birationnelle et contracte un sous-ensemble alg\'ebrique
de codimension sup\'erieure ou \'egale \`a $2$).}

\medskip

Les cas (a) et (b) sont les ``bons" cas~: dans le cas (a),
on r\'eduit la compr\'ehension de la vari\'et\'e $X$
\`a celle d'une vari\'et\'e de dimension plus petite
et \`a la structure des fibres qui sont des vari\'et\'es
de Fano. Dans le cas (b), la vari\'et\'e singuli\`ere $Y$ est
$\QQ$-factorielle, \`a singularit\'es
terminales (ce sont les singularit\'es qui permettent de donner
encore un sens \`a l'expression ``$K_Y$ est ou n'est
pas nef") et le nombre de Picard de $Y$ est strictement
plus petit que celui de $X$. Le cas (c) est le
``mauvais" cas~: les singularit\'es de $Y$ sont telles
que $Y$ ne poss\`ede pas de fibr\'e canonique et il n'est pas clair
que $Y$ soit ``plus simple" que $X$.

\subsection{Contractions divisorielles}

Nous utilisons dans la suite un r\'esultat plus
pr\'ecis que le th\'eor\`eme de contraction dans
le cas d'une contraction divisorielle. Sous cette
forme, il appara\^{\i}t dans les travaux de T.\ Ando \cite{And85}
et M.\ Beltrametti \cite{Bel86}~:

\bigskip

\noindent {\bf Th\'eor\`eme (T.\ Ando,  M.\ Beltrametti, 1985) }
{\em Soient $X$ une vari\'et\'e projective et
$f : X \to Y$ une contraction divisorielle d'un rayon
extr\^emal. Soit enfin $F$ une fibre g\'en\'erale
de $f_E : E \to f(E)$ o\`u $E$ est le diviseur exceptionnel
de $f$.

\noindent Alors il existe un fibr\'e en droites $L$ sur $X$ tel que~:

(i) $\ima (\Pic (X) \to \Pic (F)) = \ZZ \cdot L_{|F}$ o\`u $L_{|F}$
est ample sur $F$,

(ii) ${\cal O}_F(-K_X) \simeq {\cal O}_F(pL)$ et
${\cal O}_F(-E) \simeq {\cal O}_F(qL)$ o\`u $p$ et $q$ sont
deux entiers positifs.

\noindent Enfin, si $F$ est de dimension
$2$, alors $F$ est isomorphe \`a $\PP ^2$ ou \`a la
quadrique ${\cal Q}_2$.
}

\medskip

Nous utiliserons ce r\'esultat lorsque $X$ est de dimension $4$.

\subsection{L'in\'egalit\'e de Wi\'sniewski}

Avant d'\'enoncer cette in\'egalit\'e, nous avons besoin de deux notations~:
soient $X$ une vari\'et\'e projective, $R$ un rayon extr\^emal
de $\overline{\NE } (X)$ et $f$ la contraction de Mori associ\'ee.

On note
$\displaystyle{ l(R) = \min \{ -K_X \cdot C \ | \
C \mbox{ est une courbe rationnelle telle que } [C] \in R \} }.$
Le nombre $l(R)$ est la {\bf longueur} du rayon $R$.

On note aussi $A(R)$ le lieu de $X$ couvert par les courbes
dont la classe appartient \`a $R$.

Dans \cite{Wis91}, J.\ Wi\'sniewski d\'emontre l'in\'egalit\'e
fondamentale suivante~:

\bigskip

\noindent {\bf Th\'eor\`eme (J.\ Wi\'sniewski, 1991) }
{\em Pour toute fibre non triviale $F$ de $f$, on a~:
$$ \dim F + \dim A(R) \geq \dim X + l(R) - 1.$$}

Cette in\'egalit\'e, dont une forme faible
est due \`a P.\ Ionescu \cite{Ion86}, a de nombreuses cons\'equences
dans la classification des vari\'et\'es projectives
de dimension sup\'erieure ou \'egale \`a $3$.

\section{Un premier r\'esultat}

Rappelons pour commencer ce paragraphe que les exemples de J.\ Koll\'ar
et K.\ O\-gui\-so d\'ej\`a cit\'es v\'erifient la propri\'et\'e
suppl\'ementaire
suivante~: les vari\'et\'es $X$ cons\-trui\-tes ne sont, bien s\^ur,
pas projectives, mais le deviennent apr\`es exactement un \'eclatement
le long d'une sous-vari\'et\'e $Y \subset X$.

Dans toute cette partie, nous nous pla\c cons dans la situation
analogue suivante~:

\medskip

\noindent {\bf Hypoth\`ese }
{\em $X$ est une vari\'et\'e de Moishezon non projective de dimension
$n$ dont le groupe de Picard
est \'egal \`a $\ZZ$, dont
le fibr\'e canonique est gros et telle qu'il existe
une sous-vari\'et\'e $Y \subset X$ de sorte que l'\'eclatement
$\pi : \tilde{X} \to X$ de $X$ le
long de $Y$ d\'efinisse une vari\'et\'e projective $\tilde{X}$.
On note
$E$ le diviseur exceptionnel de l'\'eclatement.
}

\medskip

La remarque suivante est tr\`es importante~:

\medskip

\noindent {\bf Remarque }
D'apr\`es le corollaire 3.1.2, on sait
alors que $K_X$ (respectivement $K_{\tilde{X}}$) est
strictement positif sur les courbes
non incluses dans $Y$ (respectivement dans $E$). La cons\'equence suivante
sera utilis\'ee dans la suite~: si $C$ est une courbe
de $Y$ sur laquelle $K_X$ est strictement n\'egatif, cette courbe
ne peut pas se d\'eformer (dans $X$) hors de $Y$.

\medskip

La m\'ethode que nous adoptons pour \'etudier
$X$ et $Y$ est d'appliquer la th\'eorie de Mori
\`a la vari\'et\'e projective $\tilde{X}$.

\subsection{C\^one des courbes sur $\tilde{X}$}

Dans le cadre de notre \'etude, comme $\Pic (X) = \ZZ$,
l'espace vectoriel
$N_1(\tilde{X},\RR)$ est isomorphe \`a $\RR ^2$.

Dans la suite, nous repr\'esentons dans $N_1(\tilde{X},\RR) \simeq \RR ^2$
le c\^one ferm\'e $\overline{\NE }(\tilde{X})$~;
dans les figures ci-dessous,
ce dernier correspond \`a la
partie hachur\'ee. Si $D$ est un \'el\'ement de $\Pic (\tilde{X})$,
nous notons $D > 0$ (respectivement $D = 0$, respectivement $D < 0$)
les ensembles
$$\{ [C] \in \overline{\NE }(\tilde{X}) \ | \ D \cdot C > 0 \}$$
(respectivement $D \cdot C = 0$, respectivement $D \cdot C < 0$).

Deux cas se pr\'esentent suivant que $K_X$ est nef ou non.
Ces deux cas se distinguent naturellement~; ils correspondent,
comme nous le verrons plus loin,
au fait que $X$ admet ou non un morphisme vers une vari\'et\'e
(\'eventuellement singuli\`ere) projective de m\^eme dimension,

\medskip

(i) soit $K_X$ est nef et le dessin est le suivant~:

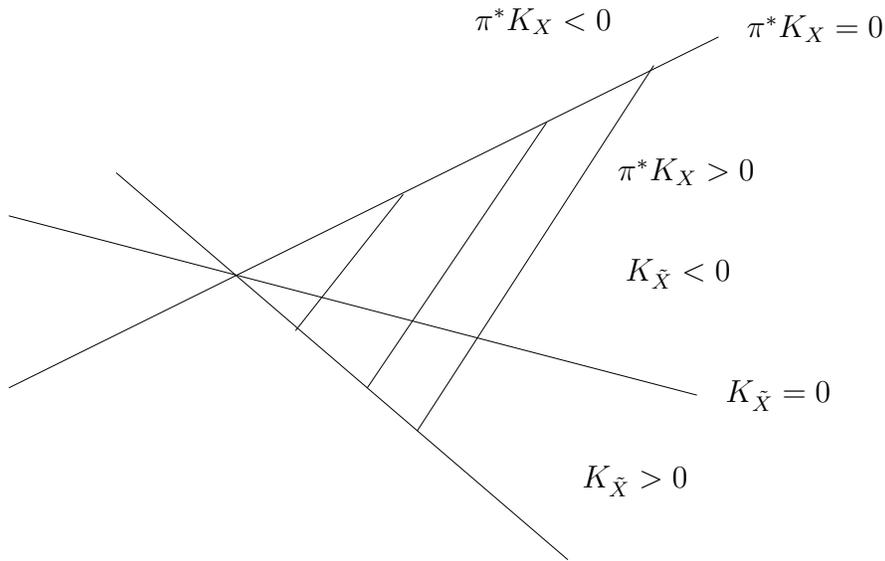
\setlength{\unitlength}{0.0125in}
\begin{picture}(404,250)(0,-10)
\path(0,144)(288,69)
\path(165,153)(120,96)
\path(225,183)(150,72)
\path(45,162)(234,0)
\path(171,54)(270,207)
\path(0,72)(297,219)
\put(309,219){\makebox(0,0)[lb]{\raisebox{0pt}[0pt][0pt]{\shortstack[l]{{\twlrm
$\pi^{*}K_{X}=0$}}}}}
\put(258,117){\makebox(0,0)[lb]{\raisebox{0pt}[0pt][0pt]{\shortstack[l]{{\twlrm
$K_{\tilde{X}}<0$}}}}}
\put(300,66){\makebox(0,0)[lb]{\raisebox{0pt}[0pt][0pt]{\shortstack[l]{{\twlrm
$K_{\tilde{X}}=0$}}}}}
\put(195,222){\makebox(0,0)[lb]{\raisebox{0pt}[0pt][0pt]{\shortstack[l]{{\twlrm
$\pi^{*}K_{X}<0$}}}}}
\put(255,159){\makebox(0,0)[lb]{\raisebox{0pt}[0pt][0pt]{\shortstack[l]{{\twlrm
$\pi^{*}K_{X}>0$}}}}}
\put(240,30){\makebox(0,0)[lb]{\raisebox{0pt}[0pt][0pt]{\shortstack[l]{{\twlrm
$K_{\tilde{X}}>0$}}}}}
\end{picture}


(ii) soit $K_X$ n'est pas nef et le dessin est le suivant~:

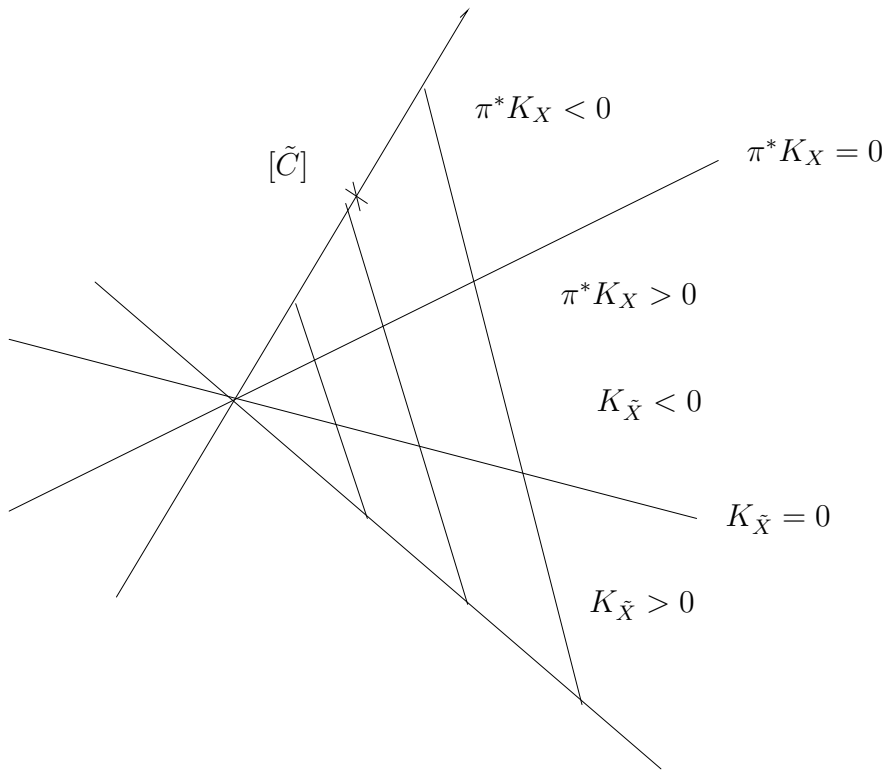
\setlength{\unitlength}{0.0125in}
\begin{picture}(404,333)(0,-10)
\path(0,180)(288,105)
\path(45,72)(192,318)(189,315)
\path(141,243)(150,237)
\path(120,195)(150,105)
\path(141,237)(192,69)
\path(174,285)(240,27)
\path(144,246)(147,234)
\path(0,108)(297,255)
\path(36,204)(273,0)
\put(108,249){\makebox(0,0)[lb]{\raisebox{0pt}[0pt][0pt]{\shortstack[l]{{\twlrm
$[\tilde{C}]$}}}}}
\put(195,273){\makebox(0,0)[lb]{\raisebox{0pt}[0pt][0pt]{\shortstack[l]{{\twlrm
$\pi^{*}K_{X}<0$}}}}}
\put(309,255){\makebox(0,0)[lb]{\raisebox{0pt}[0pt][0pt]{\shortstack[l]{{\twlrm
$\pi^{*}K_{X}=0$}}}}}
\put(231,195){\makebox(0,0)[lb]{\raisebox{0pt}[0pt][0pt]{\shortstack[l]{{\twlrm
$\pi^{*}K_{X}>0$}}}}}
\put(300,102){\makebox(0,0)[lb]{\raisebox{0pt}[0pt][0pt]{\shortstack[l]{{\twlrm
$K_{\tilde{X}}=0$}}}}}
\put(243,66){\makebox(0,0)[lb]{\raisebox{0pt}[0pt][0pt]{\shortstack[l]{{\twlrm
$K_{\tilde{X}}>0$}}}}}
\put(246,150){\makebox(0,0)[lb]{\raisebox{0pt}[0pt][0pt]{\shortstack[l]{{\twlrm
$K_{\tilde{X}}<0$}}}}}
\end{picture}


\noindent {\bf Quelques commentaires sur ces diagrammes}

- dans les deux cas, le fait que la droite $\{ K_{\tilde{X}} = 0 \}$
coupe le c\^one effectif vient du fait qu'il y a \`a la fois des
courbes sur lesquelles $ K_{\tilde{X}}$ est strictement positif
(celles non contenues dans le diviseur exceptionnel) et
des
courbes sur lesquelles $ K_{\tilde{X}}$ est strictement
n\'egatif (toute courbe incluse dans les fibres de l'\'eclatement),

- dans le deuxi\`eme cas, la position relative de
$\{ \pi ^{\ast} K_{X} = 0 \}$ est justifi\'ee par le fait qu'il
y a des courbes sur lesquelles $\pi ^{\ast} K_{X}$ et $K_{\tilde{X}}$
sont strictement positifs (celles non contenues dans le diviseur
exceptionnel d'apr\`es 3.1.2) et que $\pi ^{\ast} K_{X}$ est nul sur
toute courbe incluse dans les fibres de l'\'eclatement.

\subsub{Quelques cons\'equences de ces diagrammes}

Nous regroupons ici les renseignements provenant
directement
de la description de $\overline{\NE }(\tilde{X})$.
Pour cela, appliquons le th\'eor\`eme du c\^one
\`a la vari\'et\'e projective $\tilde{X}$. On en
d\'eduit que le rayon extr\^emal du c\^ot\'e $K_{\tilde{X}} < 0$ est
engendr\'e par la classe d'une courbe rationnelle $\tilde{C}$
dans $\tilde{X}$. Alors,
le th\'eor\`eme de contraction assure
l'existence d'une
vari\'et\'e (en g\'en\'eral singuli\`ere)
projective $Z$ et d'un morphisme $f$ associ\'es
\`a la courbe extr\^emale rationnelle
$\tilde{C}$ de sorte
que la situation suivante ait lieu~:

\begin{center}
\setlength{\unitlength}{0.0125in}
\begin{picture}(216,110)(0,-10)
\path(147,69)(147,18)
\path(145.000,26.000)(147.000,18.000)(149.000,26.000)
\path(99,81)(24,81)
\path(32.000,83.000)(24.000,81.000)(32.000,79.000)
\put(60,84){\makebox(0,0)[lb]{\raisebox{0pt}[0pt][0pt]{\shortstack[l]{{\twlrm
$f$}}}}}
\path(117,72)(117,15)
\path(115.000,23.000)(117.000,15.000)(119.000,23.000)
\put(108,75){\makebox(0,0)[lb]{\raisebox{0pt}[0pt][0pt]{\shortstack[l]{{\twlrm
$\tilde{X} \supset E$}}}}}
\put(0,75){\makebox(0,0)[lb]{\raisebox{0pt}[0pt][0pt]{\shortstack[l]{{\twlrm
$Z$}}}}}
\put(126,39){\makebox(0,0)[lb]{\raisebox{0pt}[0pt][0pt]{\shortstack[l]{{\twlrm
$\pi$}}}}}
\put(111,0){\makebox(0,0)[lb]{\raisebox{0pt}[0pt][0pt]{\shortstack[l]{{\twlrm
$X \supset Y$}}}}}
\end{picture}
\end{center}

\medskip

\noindent (i) Cas o\`u $K_X$ est nef.

Alors le rayon extr\^emal est engendr\'e
par la classe d'une courbe rationnelle incluse
dans une fibre non triviale de $\pi$. Toutes les
fibres de $\pi$ sont donc contract\'ees par $f$ si bien
que $f$ se factorise en une application $g : X \to Z$
$$ \tilde{X} \stackrel{\pi}{\to} X \stackrel{g}{\to} Z  \
\mbox{et} \ \ f = g \circ \pi.$$

\medskip

\noindent (ii) Cas o\`u $K_X$ n'est pas nef.

Alors les fibres de $f$ et les fibres de $\pi$
ne se coupent que sur un nombre fini de points~: en effet,
il n'existe pas de courbes simultan\'ement contract\'ees
par $\pi$ et $f$ car les rayons engendr\'es par $[\tilde{C}]$
et la classe d'une courbe rationnelle incluse
dans une fibre non triviale de $\pi$ sont distincts.

\subsub{Une application imm\'ediate}

Dans le cas o\`u $K_X$ n'est pas nef, la courbe
rationnelle $\tilde{C}$
n'\'etant pas contract\'ee par
$\pi$, la courbe rationnelle
$C = \pi (\tilde{C})$ v\'erifie $K_X \cdot C < 0$. Le r\'esultat suivant
en d\'ecoule~:

\medskip

\noindent {\bf Proposition }{\em
Sous les hypoth\`eses pr\'ec\'edentes, si $K_X$ n'est pas nef,
il existe une courbe rationnelle $C \subset Y $
sur laquelle $K_X$ est strictement n\'egatif.
}

\subsection{Contraction de Mori de $\tilde{X}$}

Nous \'etudions ici plus en d\'etail la contraction de Mori $f$
associ\'ee \`a la courbe extr\^emale rationnelle
$\tilde{C} \subset
\tilde{X}$ obtenue pr\'ec\'edemment pour en d\'eduire
une estimation de la dimension de $Y$ en toute dimension.

\subsub{\'Enonc\'e du r\'esultat}

Nous avons vu qu'il y a trois types de contractions
extr\^emales. Le th\'eor\`eme suivant restreint
les possibilit\'es dans notre situation~:

\bigskip

\noindent {\bf Th\'eor\`eme F } {\em
Soit $X$ de Moishezon avec $\textstyle{\Pic (X)} = \ZZ$ et $K_X$
gros. Si $X$ est rendue projective apr\`es \'eclatement $\pi : \tilde{X} \to X$
le
long de $Y$, alors~:

(i) on a $\dim \tilde{X} = \dim Z$, autrement dit
$f$ est une application birationnelle,

(ii) si $f$ est une contraction
divisorielle, son diviseur
exceptionnel est \'egal \`a celui de $\pi$
(not\'e $E$ pr\'ec\'edemment)~; ce cas est le seul possible
lorsque $K_X$ est nef,

(iii) si $f$ est une contraction
divisorielle et si
$K_X$ n'est pas nef, les in\'egalit\'es
suivantes sont satisfaites~:
$$ \codim Y -1 \leq \dim f(E) < \dim Y \ \ \mbox{\rm et} \ \
\dim Y > \frac{n-1}{2}, $$

(iv) si $f$ est une petite contraction et si
$K_X$ n'est pas nef, l'in\'egalit\'e
suivante est satisfaite~:
$$ \dim Y \geq \frac{n+1}{2}.$$
}

\smallskip

On peut reformuler ce r\'esultat sans faire intervenir
la contraction de Mori~:

\medskip

\noindent {\bf Th\'eor\`eme F' } {\em
Soit $X$ une vari\'et\'e de Moishezon de dimension
$n$ avec $\textstyle{\Pic (X)} = \ZZ$
et $K_X$ gros. Supposons que $X$ est rendue projective
apr\`es \'eclatement le long d'une sous-vari\'et\'e lisse $Y$.

\noindent Alors, si $K_X$ n'est pas nef, on a
$\displaystyle{ \dim Y > \frac{n-1}{2}}$.}

\bigskip

\noindent {\bf Remarque }
Le fait que $K_X$ soit gros est ici essentiel.
En effet, les constructions de J.\ Koll\'ar
et K.\ Oguiso montrent que les in\'egalit\'es
du point (iii) ne sont pas vraies en g\'en\'eral.
La construction de J.\ Koll\'ar donne aussi un exemple
o\`u les diviseurs exceptionnels de $\pi$ et $f$ ne
sont pas \'egaux.

\subsub{D\'emonstration du th\'eor\`eme F-F'}

Les points (i) et (ii) du th\'eor\`eme~F sont faciles~:
par hypoth\`ese, $f$ ne contracte que des courbes sur lesquelles
$K_{\tilde{X}}$ est strictement n\'egatif, donc incluses dans $E$
d'apr\`es le corollaire 3.2.1~; en particulier le point (i) est
d\'emontr\'e.

\noindent Ceci montre aussi que si $f$ est une contraction
divisorielle, son diviseur exceptionnel \'etant
inclus dans $E$ est donc \'egal \`a $E$. R\'eciproquement,
si $K_X$ est nef, $f$ se factorise \`a travers $\pi$ et donc
est une contraction divisorielle. Le point (ii) est d\'emontr\'e.

\bigskip

Montrons le point (iii) du th\'eor\`eme F~:
$f$ est une contraction divisorielle
et $K_X$ n'est pas nef. La situation
est r\'esum\'ee par le diagramme suivant~:

\begin{center}
\setlength{\unitlength}{0.0125in}
\begin{picture}(207,130)(0,-10)
\path(138,60)(138,15)
\path(136.000,23.000)(138.000,15.000)(140.000,23.000)
\path(168,90)(168,15)
\path(166.000,23.000)(168.000,15.000)(170.000,23.000)
\path(44.000,77.000)(36.000,75.000)(44.000,73.000)
\path(36,75)(126,75)
\path(147,87)	(151.921,89.353)
	(155.375,90.591)
	(159.000,90.000)

\path(159,90)	(158.773,86.022)
	(156.650,82.657)
	(153.000,78.000)

\put(147,45){\makebox(0,0)[lb]{\raisebox{0pt}[0pt][0pt]{\shortstack[l]{{\twlrm
$\pi$}}}}}
\path(41.000,107.000)(33.000,105.000)(41.000,103.000)
\path(33,105)(153,105)
\put(132,0){\makebox(0,0)[lb]{\raisebox{0pt}[0pt][0pt]{\shortstack[l]{{\twlrm
$X \supset Y $}}}}}
\put(135,69){\makebox(0,0)[lb]{\raisebox{0pt}[0pt][0pt]{\shortstack[l]{{\twlrm
$\tilde{X}$}}}}}
\put(0,102){\makebox(0,0)[lb]{\raisebox{0pt}[0pt][0pt]{\shortstack[l]{{\twlrm
$f(E)$}}}}}
\put(12,87){\makebox(0,0)[lb]{\raisebox{0pt}[0pt][0pt]{\shortstack[l]{{\twlrm
$\bigcap$}}}}}
\put(165,99){\makebox(0,0)[lb]{\raisebox{0pt}[0pt][0pt]{\shortstack[l]{{\twlrm
$E$}}}}}
\put(9,69){\makebox(0,0)[lb]{\raisebox{0pt}[0pt][0pt]{\shortstack[l]{{\twlrm
$Z$}}}}}
\end{picture}
\end{center}

\noindent Dans cette situation, on \'ecrit
$$ K_{\tilde{X}}=f^{\ast}K_Z + aE =
\pi^{\ast}K_X + (r-1)E $$
o\`u $r= \codim  Y$ et o\`u $a$ est un nombre rationnel.
Cette \'egalit\'e est une \'egalit\'e
de $\QQ$-diviseurs de Cartier~: dans le cas
d'une contraction divisorielle, le diviseur
canonique de $Z$ n'est pas de Cartier en g\'en\'eral
mais l'un de ses multiples entiers l'est.
Rappelons que toutes les notions de positivit\'e
(telle que par exemple \^etre gros, nef ou ample)
s'\'etendent naturellement aux $\QQ$-diviseurs
de Cartier.

\noindent Les in\'egalit\'es cherch\'ees d\'ecoulent
imm\'ediatement
du lemme suivant~:

\medskip

\noindent {\bf Lemme } {\em
Si $K_X$ n'est pas nef, les nombres $a$ et $r$ v\'erifient
les deux in\'egalit\'es suivantes~:

(i) $a > r-1$,

(ii) $ \codim f(E) + r \leq n + 1$.

\noindent L'in\'egalit\'e suivante est vraie en toute
g\'en\'eralit\'e
pour une contraction divisorielle~:

(iii) $a \leq \codim  f(E) - 1$.
}

\bigskip

\noindent {\bf D\'emonstration du lemme}

\medskip

\noindent {\bf In\'egalit\'e (i)~:}

Comme $Z$ est projective avec $\Pic (Z) = \ZZ$ et $K_Z$ gros,
on en d\'eduit que $K_Z$ est ample et donc
que $f^{\ast}K_Z$ est nef, et strictement
positif sur les courbes de $\tilde{X}$ non
contract\'ees par $f$. Choisissons
alors une courbe rationnelle $R$ incluse dans une fibre
non triviale
de $\pi$ (ces derni\`eres sont des $\PP ^{r-1}$, on
prend pour $R$ une droite $\PP ^{1}$).

\noindent L'\'egalit\'e
$$f^{\ast}K_Z \cdot R + a E \cdot R =
\pi^{\ast}K_X \cdot R + (r-1)E \cdot R$$
donne alors~:
$$a-(r-1) = f^{\ast}K_Z \cdot R > 0$$
car, $K_X$ n'etant pas nef,
$R$ n'est pas contract\'ee par $f$. \finpreuve

\medskip

\noindent {\bf In\'egalit\'e (ii)~:}

Cette in\'egalit\'e d\'ecoule de suite du fait
que les fibres de $f$ et $\pi$ dans $E$ ne peuvent se couper
qu'en un nombre fini de points. De l\`a ~:
$$ (n-1 - \dim f(E) ) + (r-1) \leq n-1. \finpreuve$$

\medskip

\noindent {\bf In\'egalit\'e (iii)~:}

Soit $F$ une fibre g\'en\'erique de la restriction de $f$
\`a $E$, et soit $\tilde{C}$ une courbe dans $F$.

\noindent Alors, on a
$$aE \cdot \tilde{C} = K_{\tilde{X}} \cdot \tilde{C}$$
et par la formule d'adjonction~:
$$ K_{\tilde{X} | E} = K_E - E_{| E}.$$
\noindent De l\`a, on en
d\'eduit~:
$$a+1 = \frac{K_E \cdot \tilde{C}}{E \cdot \tilde{C}}.$$

\noindent Comme le fibr\'e canonique $K_F$ est simplement la restriction
de $K_{E}$ \`a $F$, on obtient~:
$$a = \frac{K_F \cdot \tilde{C}}{E \cdot \tilde{C}} - 1.$$

\noindent Or, la vari\'et\'e $F$ est Fano, et par le th\'eor\`eme
du c\^one appliqu\'e \`a $F$, on peut supposer
que $\tilde{C}$ est une courbe (rationnelle) satisfaisant~:

$$ 0 < -K_F \cdot \tilde{C} \leq \dim F + 1 = n-1 -\dim f(E) +1 = \codim f(E)
.$$

\noindent De l\`a, comme $E \cdot \tilde{C}$ est un entier strictement
n\'egatif, il vient $a  \leq \codim f(E) - 1$. Ceci termine la
preuve du lemme.\finpreuve

\medskip

\noindent {\bf Remarque  } L'in\'egalit\'e (iii) peut aussi se
d\'eduire de l'in\'egalit\'e de Wi\'sniewski~: en effet $f$
\'etant divisorielle, on a $\dim A(R) = n-1$ et la dimension
de la fibre g\'en\'erique non triviale est $n-1-\dim f(E)$.
De l\`a, l'in\'egalit\'e de Wi\'sniewski donne $$\codim f(E) -1 \geq l(R).$$
Or, comme
$K_{\tilde{X}} = f^*K_Z + a E$, on a $-K_{\tilde{X}}\cdot C \geq a$
pour toute courbe contract\'ee, d'o\`u $l(R) \geq a$ comme souhait\'e.

\bigskip

Montrons maintenant le point (iv) du th\'eor\`eme~F.
Pour cela, on applique l'in\'egalit\'e de Wi\'sniewski~:
comme $f$ est une petite contraction, on a
\'evidemment $$\dim A(R) \leq n-2,$$ et si $F$ est une fibre
non triviale de $f$ cette derni\`ere est incluse dans $E$
et ne coupe les fibres de $\pi$ que sur un ensemble
fini. On en d\'eduit que $\dim F \leq \dim Y$ d'o\`u~:
$$ n-2 + \dim Y \geq n + l(R) -1,$$
soit $$ \dim Y \geq l(R) +1.$$
Il suffit alors d'estimer $l(R)$. Or, on a
$K_{\tilde{X}} = \pi^*K_X + (r-1) E$ et comme
$K_X$ n'est pas nef, $ \pi^*K_X$ est strictement n\'egatif sur les
courbes contract\'ees par $f$. On en d\'eduit que
$$ -K_{\tilde{X}} \cdot C \geq r$$
pour toute courbe contract\'ee par $f$. De l\`a,
$ l(R) \geq r$
et en reportant
$$ 2 \dim Y \geq n+1$$
qui est l'in\'egalit\'e souhait\'ee.\finpreuve

\subsection{Application \`a la dimension $3$}

On d\'eduit du th\'eor\`eme F un r\'esultat pr\'ecisant celui
de J.\ Koll\'ar dans notre situation~:

\medskip

\noindent {\bf Corollaire } {\em
Soit $X$ une vari\'et\'e de Moishezon non projective
de dimension $3$,
avec $\Pic(X) = \ZZ$ et $K_X$ gros.

\noindent Si $X$ peut \^etre rendue projective apr\`es
un \'eclatement seulement,
alors $X$ est une petite modification d'une vari\'et\'e
singuli\`ere projective ayant une unique singularit\'e
nodale ordinaire (dont le mod\`ele local est $xy-zt = 0$ dans
$(\CC ^4,0)$). En particulier, le fibr\'e canonique $K_X$
est nef.}

\medskip

\noindent {\bf D\'emonstration du corollaire}

On note toujours $\pi$ l'\'eclatement
rendant $X$ projective et $f$ la contraction
de Mori d\'efinie sur la vari\'et\'e projective
$\tilde{X}$. D'apr\`es le th\'eor\`eme F, $f$ est
birationnelle, et comme il n'y a pas de petites
contractions en dimension $3$ d'une vari\'et\'e
non singuli\`ere, c'est que $f$ est une contraction
divisorielle.
De plus, les in\'egalit\'es (iii) du th\'eor\`eme F
ne peuvent \^etre v\'erifi\'ees ici car elles
impliquent $\codim  Y = 1$.
C'est donc que $K_X$ est nef (on retrouve ainsi le r\'esultat
de J.\ Koll\'ar), et que le rayon extr\^emal
du c\^ot\'e $K_{\tilde{X}} < 0$ est engendr\'e par la classe d'homologie
des fibres de $\pi$. Il y a alors exactement deux possibilit\'es~:

- les fibres de la contraction de Mori (restreinte au diviseur
exceptionnel $E$) sont de dimension $1$ et alors cette derni\`ere
co\"{\i}ncide avec $\pi$. Dans ce cas, $X$ est projective, ce que
l'on a exclu,

- la contraction de Mori contracte le diviseur exceptionnel $E$
sur un point. Dans ce cas, nous appliquons le r\'esultat
de S.\ Mori \cite{Mor82} qui donne la liste de toutes les contractions
extr\^emales d'une vari\'et\'e non singuli\`ere de dimension $3$.
On en d\'eduit que le diviseur exceptionnel de $\pi$ (\'egal \`a celui
de $f$)
est isomorphe \`a
$\PP ^1 \times \PP ^1$ et que
$\displaystyle{{\cal O}_{E} (E) = N_{E/\tilde{X}}}$ est de type
$(-1,-1)$. La situation est alors la suivante~:
$$ \tilde{X} \stackrel{\pi}{\to} X \stackrel{g}{\to} Z  \
\mbox{et} \ \ f = g \circ \pi,$$
o\`u $Z$ est une vari\'et\'e
singuli\`ere projective ayant une unique singularit\'e
nodale ordinaire (dont le mod\`ele local est $xy-zt = 0$ dans
$(\CC ^4,0)$).
Dans ce cas, la contraction de Mori est alors $g \circ \pi$
et correspond \`a l'\'eclatement du point singulier~: le centre $Y$
de l'\'eclatement $\pi$ est une courbe rationnelle lisse.\finpreuve

\bigskip

\noindent {\bf Exemple }
La situation pr\'ec\'edente peut effectivement \^etre
r\'ealis\'ee~: soit $Z$ une hypersurface de $\PP ^4$
d'\'equation
$$h_0x_0^2 + h_1x_1^2 + h_2x_2^2 + h_3x_3^2 = 0,$$
o\`u $[x_0 : \cdots : x_4]$ sont les coordonn\'ees homog\`enes
dans $\PP^4$ et o\`u les $h_i$ sont quatre polyn\^omes
homog\`enes de degr\'e $d$ sup\'erieur ou
\'egal \`a $4$ ne s'annulant pas en $[0:0:0:0:1]$ et g\'en\'eriques
parmi les polyn\^omes ayant ces propri\'et\'es.
Alors l'hypersurface $Z$ est lisse except\'e au point
$[0:0:0:0:1]$ o\`u elle poss\`ede une singularit\'e
nodale ordinaire. On obtient $X$ comme d\'ecrit pr\'ec\'edemment
en r\'esolvant la singularit\'e puis en contractant dans une direction
de la quadrique exceptionnelle.

\medskip

Nous donnons dans la suite d'autres applications du th\'eor\`eme F
mais nous commen\c cons par montrer dans le paragraphe suivant
que le r\'esultat de J.\ Koll\'ar ne s'\'etend pas en dimension
sup\'erieure ou \'egale \`a $4$.

\section{Une famille de vari\'et\'es de Moishezon}

Le but de cette partie est de montrer le r\'esultat suivant~:

\medskip

\noindent {\bf Th\'eor\`eme G} {\em
Pour tout entier $n$ sup\'erieur
ou \'egal \`a $4$, il existe des vari\'et\'es de Moishezon $X$
de dimension $n$ v\'erifiant~:

(i) $\Pic(X) = \ZZ$, (ii) $K_X$ est gros, (iii) $K_X$ n'est pas nef.
}

\medskip

Ainsi, le r\'esultat de J.\ Koll\'ar est propre \`a la dimension
$3$.

\medskip

\noindent {\bf Remarque } Les vari\'et\'es obtenues dans la
construction qui suit
rel\`event toutes du cas ``contraction divisorielle"
\'evoqu\'e dans le paragraphe pr\'ec\'edent.
Il serait bien s\^ur int\'eressant de cons\-truire de telles
vari\'et\'es relevant du cas ``petite contraction".
Cependant, nous verrons plus loin que ce cas ne peut pas
se produire en dimension~$4$.

\subsection{Un r\'esultat interm\'ediaire}

La d\'emonstration du th\'eor\`eme G repose
sur la proposition suivante, que nous prouvons plus loin.
Mentionnons qu'il nous a \'et\'e signal\'e par un rapporteur
anonyme que cette proposition se trouve dans \cite{BVV78}.

Pour $n$ entier, nous notons $[x_0: \cdots :x_{n+1}]$
les coordonn\'ees homog\`enes dans $\PP ^{n+1}$. On d\'esigne
par $\PP _{x_n} ^{1}$ la droite
$\{ x_0 = \dots = x_{n-1} = 0 \}$.
Choisissons alors $n$ polyn\^omes
homog\`enes $h_0,\ldots \!,h_{n-1}$ de degr\'e $2n-2$ et
consid\'erons
l'hypersurface $Z$ de degr\'e $2n-1$ dans $\PP ^{n+1}$
et d'\'equation
$$\sigma = x_0h_0 +\cdots+ x_{n-1}h_{n-1} = 0.$$
Cette hypersurface contient $\PP _{x_n} ^{1}$ et peut \^etre
singuli\`ere. On a cependant le r\'esultat suivant~:

\medskip

\noindent {\bf Proposition }{\em
Si $n$ est sup\'erieur ou \'egal \`a $3$ et si les $h_{i}$
sont choisis g\'en\'eriquement, alors~:

\smallskip

\!(i) l'hypersurface $Z$ est non singuli\`ere,

\smallskip

\!(ii) le fibr\'e normal $N_{\PP ^{1}/ Z}$ est \'egal \`a
${\cal O}_{\PP ^{1}}(-1)^{\oplus n-1}$,

\smallskip

\!(iii) $K_Z$ est \'egal \`a ${\cal O}_{\PP ^{n+1}}(n-3)_{| Z}$,

\smallskip

\!(iv) $\Pic (Z) = \ZZ$.
}

\subsection{D\'emonstration du th\'eor\`eme G }

La construction qui suit nous a \'evidemment
\'et\'e inspir\'ee par l'analyse du paragraphe pr\'ec\'edent, dans
le cas o\`u la contraction de Mori est une contraction
divisorielle~: si une vari\'et\'e de dimension $4$
de Moishezon satisfait le point (iii) du th\'eor\`eme~F,
c'est en \'eclatant une surface, puis en contractant sur une
courbe rationnelle que l'on obtient un mod\`ele projectif. Nous
donnons cependant la construction g\'en\'erale en toute dimension.

\medskip

\noindent {\bf Construction explicite : }

On se fixe dor\'enavant une hypersurface $Z$ donn\'ee par la proposition
pr\'ec\'edente. La vari\'et\'e $X$ cherch\'ee va \^etre obtenue
en effectuant un ``flip" (plus exactement l'inverse
d'un flip) \`a partir de $Z$.

\medskip

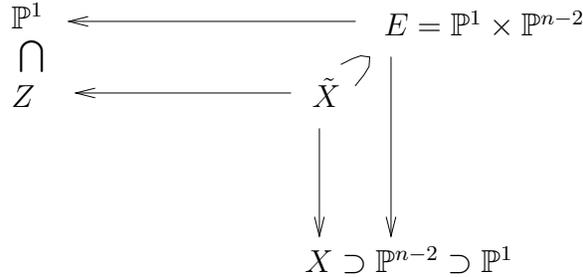
\begin{center}
\setlength{\unitlength}{0.0125in}
\begin{picture}(322,130)(0,-10)
\path(129,60)(129,15)
\path(127.000,23.000)(129.000,15.000)(131.000,23.000)
\path(159,90)(159,15)
\path(157.000,23.000)(159.000,15.000)(161.000,23.000)
\path(35.000,77.000)(27.000,75.000)(35.000,73.000)
\path(27,75)(117,75)
\path(138,87)	(142.921,89.353)
	(146.375,90.591)
	(150.000,90.000)

\path(150,90)	(149.773,86.022)
	(147.650,82.657)
	(144.000,78.000)

\path(32.000,107.000)(24.000,105.000)(32.000,103.000)
\path(24,105)(144,105)
\put(3,87){\makebox(0,0)[lb]{\raisebox{0pt}[0pt][0pt]{\shortstack[l]{{\twlrm
$\bigcap$}}}}}
\put(123,0){\makebox(0,0)[lb]{\raisebox{0pt}[0pt][0pt]{\shortstack[l]{{\twlrm
$X \supset \PP^{n-2}  \supset \PP ^1$}}}}}
\put(0,69){\makebox(0,0)[lb]{\raisebox{0pt}[0pt][0pt]{\shortstack[l]{{\twlrm
$Z$}}}}}
\put(126,69){\makebox(0,0)[lb]{\raisebox{0pt}[0pt][0pt]{\shortstack[l]{{\twlrm
$\tilde{X}$}}}}}
\put(0,102){\makebox(0,0)[lb]{\raisebox{0pt}[0pt][0pt]{\shortstack[l]{{\twlrm
$\PP ^{1}$}}}}}
\put(156,99){\makebox(0,0)[lb]{\raisebox{0pt}[0pt][0pt]{\shortstack[l]{{\twlrm
$E=\PP^{1} \times \PP ^{n-2} $}}}}}
\end{picture}
\end{center}

Suivant la figure ci-dessus, notons $\tilde{X}$ la vari\'et\'e
projective obtenue en \'eclatant $Z$ le long de $\PP ^{1}$.
Le diviseur exceptionnel de l'\'eclatement est alors
$E = \PP ^{1} \times \PP ^{n-2}$, et pour pouvoir
contracter dans l'autre direction, il s'agit de montrer,
d'apr\`es le crit\`ere de contraction de Fujiki-Nakano,
que $\displaystyle{ {\cal O}(E)_{| \PP ^1} = {\cal O}_{\PP ^1}(-1) }$.
Pour cela, les deux suites exactes suivantes~:
$$ 0 \to N_{\PP^1 /E} = {\cal O}_{\PP ^1}^{\oplus n-2}
\to N_{\PP^1/\tilde{X}} \to N_{E/ \tilde{X} | \PP ^1} =
{\cal O}(E)_{| \PP ^1} \to 0 ,$$
$$ 0 \to T\PP ^1 \to T\tilde{X}_{| \PP ^1}
\to N_{\PP^1/\tilde{X}} \to 0 $$
\noindent donnent successivement~:
$$\deg (N_{\PP^1/\tilde{X}})=
\deg ({\cal O}(E)_{| \PP ^1}) \ ,\
\deg (K_{\tilde{X} | \PP ^1}) =
-2 -\deg (N_{\PP^1/\tilde{X}}).$$
\noindent Comme $K_{\tilde{X}} = f^{\ast}K_Z + (n-2){\cal O}(E)$
et $K_Z = {\cal O}_{\PP ^{n+1}}(n-3)_{| Z}$, on en
d\'eduit bien que ${\cal O}(E)_{| \PP ^1} = {\cal O}_{\PP ^1}(-1)$.

\noindent La contraction de $\PP ^{1}$ d\'efinit donc une
vari\'et\'e de Moishezon,
contenant un $\PP ^{n-2}$ et telle
que $\Pic (X) = \ZZ$.

Montrons maintenant que $N_{\PP ^{n-2} / \tilde{X}} =
{\cal O}_{\PP ^{n-2}}(-1) \oplus {\cal O}_{\PP ^{n-2}}(-1)$.
Comme $$E =
\PP ^{1} \times \PP ^{n-2} = \PP(N^{\ast}_{\PP ^{n-2} / \tilde{X}}),$$
le fibr\'e normal $N_{\PP ^{n-2} / \tilde{X}}$ est de la forme
${\cal O}_{\PP ^{n-2}}(a) \oplus {\cal O}_{\PP ^{n-2}}(a)$.
Comme pr\'ec\'edemment, la suite exacte~:
$$ 0 \to T\PP ^{n-2} \to TX_{| \PP ^{n-2}}
\to N_{\PP^{n-2}/X} \to 0 $$
\noindent donne $2a = -\deg (K_{X | \PP ^{n-2}}) -n+1$,
puis
$$\deg (K_{X | \PP ^{n-2}}) = \deg (K_{\tilde{X} | \PP ^{n-2}}
- {\cal O}(E)_{| \PP ^{n-2}}) = -(n-2) + 1,$$
\noindent d'o\`u finalement $a=-1$.

\noindent Par ailleurs, nous venons de montrer que
$$ K_{X | \PP ^{n-2}} = {\cal O}_{\PP ^{n-2}}(3-n).$$

\noindent Ainsi, si $n$ est sup\'erieur ou \'egal \`a $4$, $-K_X$
est ample sur $\PP ^{n-2}$. Finalement, le fibr\'e $K_X$ bien
que gros n'est pas nef et le th\'eor\`eme est d\'emontr\'e.
\finpreuve

\bigskip

\noindent {\bf Remarque  }
La construction pr\'ec\'edente, en dimension $3$, donne un nouvel
exemple de vari\'et\'e de Moishezon, de ``Calabi-Yau"
satisfaisant
$\displaystyle{\Pic (X) = \ZZ \cdot {\cal O}_X (1) ,}$
(voir aussi \cite{Ogu94}).

\subsection{D\'emonstration de la proposition}

On d\'emontre (i) et (ii) simultan\'ement.

Les points singuliers de $Z$ sont des z\'eros communs des
\'equations
$$ x_0h_0 +\cdots+ x_{n-1}h_{n-1} = 0$$ et
$$x_0 \frac{\partial h_0}{\partial x_i} +\cdots+
x_{n-1} \frac{\partial h_{n-1}}{\partial x_i} + h_i = 0
 \ , \ i=0,\ldots \!,n-1.$$
En particulier, $Z$ est lisse au voisinage de
$\PP ^{1} = \{ x_0 = \dots = x_{n-1} = 0 \}$ d\`es que les
$h_i$ ne s'annulent pas simultan\'ement sur $\PP ^{1}$. Ceci
est vrai pour un choix g\'en\'erique des $h_i$ d\`es que
$n$ est sup\'erieur ou \'egal \`a $2$.

On d\'eduit alors
du th\'eor\`eme de Bertini \cite{G-H78}
que si les $h_i$ sont \`a
nouveau g\'en\'eriques, l'hypersurface $Z$ est non singuli\`ere
partout~: en effet, de fa\c{c}on g\'en\'erale, une relation
$$\displaystyle{ \sum_i s_i f_i = 0 }$$
d\'efinit une vari\'et\'e
non singuli\`ere en dehors des z\'eros communs des $s_i$
d\`es que les $f_i$ sont g\'en\'eriques dans l'espace des
sections holomorphes d'un fibr\'e engendrant en tout point
les jets
d'ordre inf\'erieur ou \'egal \`a $1$.

\medskip

D\'eterminons ensuite le fibr\'e normal $N_{\PP ^{1}/ Z}$, et pour cela,
consid\'erons la suite exacte des fibr\'es normaux :
$$ 0 \to N_{\PP ^{1}/ Z} \to N_{\PP ^{1}/ \PP ^{n+1}} =
{\cal O}_{\PP ^{1}}(1)^{\oplus n} \stackrel{d \sigma}{\to}
{\cal O}_{\PP ^{n+1}}(2n-1)_{| \PP ^1} \to 0 .$$
Il est alors clair
que $N_{\PP ^{1}/ Z}$ est de degr\'e $-(n-1)$.
Pour montrer qu'il est \'egal \`a ${\cal O}_{\PP ^{1}}(-1)^{\oplus n-1}$,
il suffit donc de montrer qu'il n'a pas de sections (rappelons
en effet qu'un th\'eor\`eme d'A.\ Grothendieck affirme que
tout fibr\'e vectoriel sur $\PP ^1$ est scind\'e).

Par la suite exacte pr\'ec\'edente, une section de $N_{\PP ^{1}/ Z}$
peut \^etre vue
comme une section de ${\cal O}_{\PP ^{1}}(1)^{\oplus n}$,
annul\'ee par $d \sigma$. Une telle section
correspond \`a la donn\'ee d'un $n$-uplet
$(s_0,\ldots \!,s_{n-1})$ o\`u les $s_i$ sont des polyn\^omes
homog\`enes de degr\'e $1$ en les variables $x_n , x_{n+1}$,
que l'on \'ecrit $s_i (x) = s_{i,n}x_n + s_{i,n+1}x_{n+1}$.
Dans
$N_{\PP ^{1}/ \PP ^{n+1}}$, on a alors~:
$$\displaystyle{s = \sum_{i=0}^{n-1} s_i \frac{\partial}{\partial x_i}}.$$
De m\^eme, notons
$$\displaystyle{h_i(x) = \sum_{p=o}^{2n-2} h_{i,p}x_{n}^{p}x_{n+1}^{2n-2-p}}$$
la restriction de $h_i$ \`a $\PP ^1$.
La relation $d \sigma (s) = 0$ donne ici~:
$$\displaystyle{\sum_{i=0}^{n-1} s_i h_i = 0}.$$
Comme $$\displaystyle{ d \sigma = \sum_{i=0}^{n-1} h_idx_i }$$
le long de $\PP ^1$, cette relation
se traduit par un syst\`eme lin\'eaire \`a $2n$ \'equations en les $2n$
inconnues
$s_{i,n},s_{i,n+1}$. Il s'agit de montrer que pour un choix
g\'en\'erique des $h_i$, le d\'eterminant de la matrice suivante~:

$$
\left(
\begin{matrix}
h_{0,0}&h_{1,0}&\dots&h_{n-1,0}&0&0&\dots&0\\
h_{0,1}&h_{1,1}&\dots&h_{n-1,1}&h_{0,0}&h_{1,0}&\dots&h_{n-1,0}\\
\vdots&\vdots&&\vdots&\vdots&\vdots&&\vdots\\
h_{0,2n-2}&h_{1,2n-2}&\dots&h_{n-1,2n-2}&h_{0,2n-3}&h_{1,2n-3}&\dots&h_{n-1,2n-3}\\
0&0&\dots&0&h_{0,2n-2}&h_{1,2n-2}&\dots&h_{n-1,2n-2}
\end{matrix}
\right)
$$
n'est pas nul, ce qui est clair
en prenant par exemple $$h_{0,0}= \lambda _0, \ldots \!,
h_{n-1,n-1}= \lambda _{n-1},
h_{0,n-1}= \mu _{0}, \ldots \!, h_{n-1,2n-2}= \mu _{n-1}$$
avec $\lambda _i \neq 0$, $\mu_i \neq 0$.

\noindent Ainsi, il existe un choix des $h_i$ de sorte que
l'hypersurface $Z$ (\'eventuellement singuli\`ere) est
lisse au voisinage de $\PP ^1$, avec $N_{\PP ^{1}/ Z} =
{\cal O}_{\PP ^{1}}(-1)^{\oplus n-1}$.

\medskip

Si maintenant les $h_i$ sont choisis de sorte que (i) et (ii) soient
satisfaits, alors (iii) est imm\'ediat par adjonction et
(iv) d\'ecoule du th\'eor\`eme de Lefschetz~:
l'hypersurface $Z$ est le diviseur d'une section
d'un fibr\'e ample et les fibr\'es en droites
sur $Z$ sont restriction de fibr\'es en droites
sur $\PP ^{n+1}$.\finpreuve

\bigskip

\noindent {\bf Remarque }
Le cas $n=3$ de la proposition pr\'ec\'edente correspond \`a celui des
quintiques dans $\PP ^4$. Il a \'et\'e consid\'er\'e par S.\ Katz
dans \cite{Kat86}. Dans ce cadre, S.\ Katz d\'etermine le fibr\'e
normal $N_{\PP ^{1}/ Z}$ dans le cas o\`u les $h_i$ sont
g\'en\'eriques, mais analyse aussi la situation non
g\'en\'erique. Signalons aussi \cite{Cle83}, o\`u H.\ Clemens
consid\`ere des questions analogues, toujours en dimension $3$.

\medskip

\noindent {\bf Remarque }
On peut reprendre plus g\'en\'eralement la construction
pr\'ec\'edente
pour les hypersurfaces de $\PP ^{n+1}$ de degr\'e $2n-2k+1$
passant par un $\PP ^{k}$ lin\'eaire. On peut alors
\`a nouveau montrer que (g\'en\'eriquement)
$N_{\PP ^{k}/ Z}$ (qui est de degr\'e $k-n$) n'a pas de sections.
Cependant, le fibr\'e $N_{\PP ^{k}/ Z}$ n'est pas
scind\'e si $k \geq 2$.

\medskip

D\'emontrons ce dernier point en consid\'erant \`a nouveau
la suite exacte des fibr\'es
normaux~:

$$ 0 \to N_{\PP ^{k}/ Z} \to N_{\PP ^{k}/ \PP ^{n+1}} =
{\cal O}_{\PP ^{k}}(1)^{\oplus n+1-k} \stackrel{d \sigma}{\to}
{\cal O}_{\PP ^{n+1}}(2n-2k+1)_{| \PP ^k} \to 0 .$$
En dualisant cette suite, il vient~:

$$ 0 \to {\cal O}_{\PP ^{n+1}}(-2n+2k-1)_{| \PP ^k} \to
{\cal O}_{\PP ^{k}}(-1)^{\oplus n+1-k} \to N_{\PP ^{k}/ Z}^{\ast}
\to 0 .$$
Si $k \geq 2$, comme le fibr\'e
${\cal O}_{\PP ^{n+1}}(-2n+2k-1)_{| \PP ^k}$ est n\'egatif,
la suite exacte longue de cohomologie donne
$H^0 (\PP ^k, N_{\PP ^{k}/ Z}^{\ast}) = 0$. Ceci exclut de suite
le fait que $N_{\PP ^{k}/ Z}$ soit scind\'e car il serait
alors \'egal \`a ${\cal O}_{\PP ^{k}}(-1)^{\oplus n-k}$.~\finpreuve

\section{Une classification en dimension $4$}

Dans ce paragraphe, nous consid\'erons des vari\'et\'es
de Moishezon (non projectives) $X$ de dimension $4$.
Comme pr\'ec\'edemment, nous
supposons que $\Pic (X) = \ZZ$, que $K_X$ est gros
et que $X$ est rendue projective apr\`es \'eclatement
le long d'une sous-vari\'et\'e lisse $Y$.

\subsection{\'Enonc\'e des r\'esultats}

Nous montrons les deux r\'esultats suivants~:

\medskip

\noindent {\bf Th\'eor\`eme H }{\em
Sous les hypoth\`eses pr\'ec\'edentes, $Y$ est n\'ecessairement
une surface.

\noindent Autrement dit, et dans cette situation particuli\`ere,
il ne suffit pas d'\'eclater une courbe pour rentrer dans
le monde projectif.
}

\medskip

Nous avons vu pr\'ec\'edemment que $K_X$ n'est pas n\'ecessairement
nef \`a partir de la dimension $4$. Le r\'esultat suivant montre
que l'exemple construit dans le paragraphe pr\'ec\'edent
est le ``seul possible" dans le cas
o\`u $K_X$ n'est pas nef.
Nous reprenons les notations des paragraphes
pr\'ec\'edents,
$\displaystyle{ \pi : \tilde{X} \to X}$
d\'esigne
l'\'eclatement de $X$ le
long de $Y$ et
$\displaystyle{ f : \tilde{X} \to Z}$
d\'esigne la contraction extr\^emale de Mori sur $\tilde{X}$.

\medskip

\noindent {\bf Th\'eor\`eme I }{\em
Sous les hypoth\`eses pr\'ec\'edentes et si $K_{X}$ n'est
pas nef, alors~:

(i) le couple $(Y,N_{Y/X})$ est \'egal \`a
$(\PP ^2, {\cal O}_{\PP ^{2}}(-1)^{\oplus 2})$,

(ii) $f$ contracte le
diviseur exceptionnel de $\pi$ sur une courbe rationnelle
lisse \`a fibr\'e normal ${\cal O}_{\PP ^{1}}(-1)^{\oplus 3}$
dans une vari\'et\'e projective lisse $Z$.
En particulier, $f$ est une contraction divisorielle.}

\bigskip

Ces r\'esultats sont accessibles en dimension $4$ car les
contractions de Mori sont ``bien comprises"
gr\^ace aux r\'esultats, rappel\'es pr\'ec\'edemment
de T.\ Ando \cite{And85}
et M.\ Beltrametti \cite{Bel86} pour les contractions divisorielles.

\bigskip

\subsection{D\'emonstration du th\'eor\`eme H}

\medskip

Nous traitons s\'epar\'ement les cas $K_X$ nef
et $K_X$ non nef.

\noindent (i) Le cas $K_X$ non nef~.

Il d\'ecoule directement du th\'eor\`eme~F~: en effet,
si la contraction de Mori $f$
est une contraction divisorielle, le point (iii)
assure
que $Y$ est une surface et que $f(E)$ est une courbe.
Par ailleurs, le point (iv)
exclut la possibilit\'e que $f$ soit une petite contraction
(car sinon $\dim Y \geq 3$ !).
Mentionnons qu'une premi\`ere version de ce travail
\cite{Bo95b} excluait ce cas en utilisant le
difficile th\'eor\`eme de structure des petites contractions de
Kawamata \cite{Kaw89}. Ceci ach\`eve le cas $K_X$ non nef.

\bigskip

\noindent (ii) Le cas $K_X$ nef.

Dans ce cas, rappelons que la contraction de Mori se factorise
par $\pi$. Notons $\displaystyle{ g : X \to Z}$ de sorte
que $f = g \circ \pi$.

Raisonnons alors par l'absurde en supposant que $Y$ est une courbe.
Dans ce cas, il est clair que $f(E)$ (\'egal \`a $g(Y)$) est un point.
En effet, dans le cas contraire, $f(E)$ est une courbe et l'application
$g$ est finie. Comme $Z$ est projective, on en d\'eduit que $X$ est
projective, ce que l'on a rejet\'e.
Ainsi $f(E)$ est un point et le diviseur $E$ est une vari\'et\'e
de Fano.

\noindent Nous allons montrer que $E$ est en fait isomorphe \`a
la quadrique de dimension $3$, ce qui fournira la contradiction~;
une quadrique de dimension $3$, dont le nombre de Picard est $1$,
ne pouvant \^etre \'egale au
projectivis\'e d'un fibr\'e de rang $3$ sur une courbe,
pour lequel le nombre de Picard est $2$ !

\noindent Pour cela, remarquons que $Z$ \'etant
$\QQ$-factorielle \`a singularit\'es
terminales, il existe un entier $m$ non nul tel que $mK_Z$ est de Cartier.
Alors~:
$$mK_X = g^{\ast}(mK_Z).$$
En particulier, la restriction de
$K_X$ \`a $Y$ est triviale. Il en d\'ecoule que $K_Y = \det N_{Y/X}$,
et par cons\'equent,
$$K_E = \pi ^{\ast}(K_Y - \det N_{Y/X} ) +
3{\cal O}_{E}(-1) = 3{\cal O}_{E}(-1).$$
On en d\'eduit que
${\cal O}_{E}(1)$ est ample et que $E$ est une vari\'et\'e
(de dimension $3$) d'indice $3$~; rappelons que l'indice
d'une vari\'et\'e de Fano $V$ est le plus grand entier
$r >0$ tel qu'il existe un fibr\'e en droites
$L$ avec $-K_V = rL$.
Or, le th\'eor\`eme de Kobayashi-Ochiai \cite{KoO73} affirme
qu'{\em une vari\'et\'e de Fano de dimension $n$ et
d'indice $n$ est isomorphe \`a la quadrique ${\cal Q}_n$}. On en
d\'eduit ici que $E$ est
la quadrique ${\cal Q}_3$ comme annonc\'e.\finpreuve

\bigskip

Dans la situation du th\'eor\`eme H et lorsque
$K_X$ est nef, nous avons
vu que la contraction de Mori $f$ sur $\tilde{X}$
se factorise en une application birationnelle
$g : X \to Z$ qui contracte la surface $Y$.
La proposition suivante pr\'ecise le cas o\`u
$g(Y)$ est r\'eduit \`a un point~:

\bigskip

\noindent {\bf Proposition } {\em Si $K_X$ est nef et si
$f(E)$ (\'egal \`a $g(Y)$) est un point, alors
le couple
$(Y,N_{Y/X})$ est \'egal \`a $(\PP ^2, T^{\ast}\PP ^2)$,
$(\PP ^2, {\cal O}_{\PP ^2}(-1) \oplus {\cal O}_{\PP ^2}(-2) )$
ou $({\cal Q}_2, {\cal O}_{{\cal Q}_2}(-1,-1)^{\oplus 2} )$.}

\bigskip

Nous ne connaissons pas d'exemples explicites o\`u
ces possibilit\'es sont effectivement r\'ealis\'ees, mais
nous pouvons remarquer qu'aucune n'est exclue {\em a priori}
par les r\'esultats de T.\ Ando et M.\ Beltrametti.

\bigskip

La d\'emonstration de la proposition d\'ecoule directement
du th\'eor\`eme suivant de T.\ Peternell \cite{Pet91}~:

\medskip

\noindent {\bf Th\'eor\`eme (T.\ Peternell, 1991) } {\em
Soit $V$ une vari\'et\'e projective de dimension $n$
et soit $E$ un fibr\'e vectoriel de rang $n$ sur $V$
de sorte que $c_1(E) = c_1(X)$. Alors, le couple
$(V,E)$ est \'egal \`a
$(\PP ^n, {\cal O}_{\PP ^n}(2) \oplus {\cal O}_{\PP ^n}(1)^{\oplus n+1})$,
$(\PP ^n, T \PP ^n)$ ou $({\cal Q}_n, {\cal O}_{{\cal Q}_n}(1)^{\oplus n})$.}

\bigskip

\noindent {\bf D\'emonstration de la proposition}

La d\'emonstration du th\'eor\`eme H montre que
$$K_Y = \det N_{Y/X}$$
et que
$${\cal O}_E(1) = {\cal O}_{\PP (N_{Y/X}^{\ast})}(1)$$
est ample,
donc que $N_{Y/X}^{\ast}$ est aussi ample.
Le r\'esultat d\'ecoule
du th\'eor\`eme de T.\ Peternell appliqu\'e au couple
$(Y, N_{Y/X}^{\ast})$. \finpreuve

\subsection{D\'emonstration du th\'eor\`eme I}

Notons $F$ la fibre g\'en\'erale de $f$ restreinte
au diviseur exceptionnel $E$. Comme $f(E)$ est
une courbe, $F$ est de dimension $2$. D'apr\`es le
th\'eor\`eme de T.\ Ando et M.\ Beltrametti,
$F$ est \'egal \`a $\PP ^2$ ou \`a
la quadrique ${\cal Q}_2$. De plus, il a \'et\'e vu
pr\'ec\'edemment que $F$ coupe les fibres de $\pi$
sur des points. On en d\'eduit que $\pi _{| F} : F \to Y$
est une application surjective finie. Deux cas sont \`a distinguer~:

\medskip

- $F$ est \'egal \`a $\PP ^2$.

\noindent Dans ce cas, $Y$ est aussi \'egal
\`a
$\PP ^2$. En effet, un r\'esultat de R.\ Lazarsfeld \cite{Laz84}
affirme que si {\em $h : \PP ^n \to V$ est une application holomorphe
surjective finie sur une vari\'et\'e de dimension $n$,
alors $V$ est isomorphe \`a $\PP ^n$}~; en dimension $2$,
on peut trouver une d\'emonstration \'el\'ementaire
dans \cite{BPV84}.

\noindent Montrons alors que
$$\displaystyle{ \pi _{| F} : F \simeq \PP ^2 \to Y \simeq \PP ^2 }$$
est un isomorphisme. Pour cela, il suffit de montrer que
$\pi _{| F}$ est un isomorphisme local, car alors $\pi _{| F}$
est un rev\^etement donc le rev\^etement trivial.
Soient donc $x$ dans $F$ et $L$ un $\PP ^1$ quelconque passant par $\pi (x)$.
Sa pr\'e-image $\pi ^{-1}(L)$ est une surface d'Hirzebruch bi-r\'egl\'ee
donc $\PP ^1 \times \PP ^1$. L'intersection
$F \cap \pi ^{-1}(L)$ est alors une r\'eunion
de $\PP ^1$ ``horizontaux". La restriction de $\pi$
au $\PP ^1$ horizontal passant par $x$ est donc un
isomorphisme sur son image. Ceci \'etant vrai pour tout
$\PP ^1$ passant par $\pi (x)$, ceci montre bien
que $d \pi_{| F}(x)$ est surjective, donc inversible,
et que $\pi _{| F}$ est un isomorphisme
local.
Ainsi,
$$\displaystyle{ \pi _{| F} : F \simeq \PP ^2 \to Y \simeq \PP ^2 }$$
est un isomorphisme.
On en d\'eduit que le fibr\'e normal $N_{Y/X}$
est scind\'e~; on d\'efinit alors $a$ et $b$ en posant~:
$$N_{Y/X} = {\cal O}_{\PP ^2}(a) \oplus {\cal O}_{\PP ^2}(b).$$
Le fait que $\pi ^{-1}(L) \simeq \PP ^1 \times \PP ^1$
montre m\^eme que $a=b$.

\noindent Comme $K_X$
n'est pas nef, $K_X$ est n\'egatif sur $Y$. Il vient alors~:
$$\displaystyle{\deg (K_{X | Y}) = -3-2a < 0}$$ d'o\`u
$a \geq -1$. L'affirmation suivante permet de conclure~:

\medskip

\noindent {\bf Affirmation }
{\em L'entier $a$ est strictement n\'egatif.}

\medskip

\noindent {\bf D\'emonstration }

Par l'absurde, supposons que $a \geq 0$. Alors, si $C$
d\'esigne un $\PP ^1$ de $Y = \PP ^2$, on a~:
$$ H^1(C,N_{C/X}) =
H^1(\PP^1, {\cal O}_{\PP ^1}(1) \oplus {\cal O}_{\PP ^1}(a)^{\oplus 2}) = 0,$$
d'o\`u~:
$$ \dim_{[C]} \Hilb (X) =
\dim H^0(\PP^1, {\cal O}_{\PP ^1}(1) \oplus {\cal O}_{\PP ^1}(a)^{\oplus 2})
= 2a + 4.$$
Or, $$\dim_{[C]} \Hilb (Y) = \dim_{[\PP ^1]} \Hilb (\PP ^2) = 2.$$
Comme $a \geq 0$, on en d\'eduit que~:
$$ \dim_{[C]} \Hilb (X) > \dim_{[C]} \Hilb (Y)$$
si bien que $C$ se d\'eforme dans $X$ hors de $Y$. Ceci n'est
pas possible comme nous l'avons d\'ej\`a rencontr\'e car
$K_X$ est positif sur les courbes non incluses dans $Y$.
Ici, $K_X$, n'\'etant pas nef, est n\'egatif sur $C$. Contradiction !
\finpreuve

\medskip

Ainsi, $a=-1$ et
$f$ contracte $E$ sur une courbe
rationnelle lisse \`a fibr\'e normal ${\cal O}_{\PP ^1}(-1)^{\oplus 3}$
dans la vari\'et\'e projective lisse $Z$.

\medskip

- $F$ est \'egal \`a la quadrique ${\cal Q}_2$.

\noindent Nous montrons
que ce cas ne peut pas arriver.
En effet, $Y$ est alors isomorphe \`a $\PP ^2$ ou ${\cal Q}_2$. Le cas
$Y \simeq \PP ^2$ s'exclut exactement comme pr\'ec\'edemment~: $\pi _{| F}$
r\'ealise un isomorphisme entre la quadrique et $\PP ^2$ !

\noindent Si $Y \simeq {\cal Q}_2$, le raisonnement est plus simple et il
est inutile de montrer que
$$\pi _{| F} : F \simeq {\cal Q}_2 \to Y \simeq {\cal Q}_2$$
est un isomorphisme. Choisissons en effet un $\PP ^1$ dans $Y$,
\`a savoir un des
g\'en\'erateurs de $H_2({\cal Q}_2,\ZZ)$, sur lequel
$K_X$ est strictement n\'egatif (il en existe car $K_X$
n'est pas nef). Alors $N_{Y/X}$ restreint \`a
$\PP ^1$ est de la forme $${\cal O}_{\PP ^1}(a) \oplus {\cal O}_{\PP ^1}(a)$$
(ceci comme pr\'ec\'edemment car
$\pi ^{-1}(\PP ^1) \simeq \PP ^1 \times \PP ^1$).
La suite exacte~:
$$ 0 \to T \PP ^1 \to TX_{| \PP ^1} \to N_{\PP ^1/X} \to 0,$$
et le fait que $N_{\PP ^1/{\cal Q}_2}$ est trivial entrainent
que $$\deg (-K_{X | \PP ^1}) = 2 + 2a > 0$$ et donc que
$a \geq 0$. Ceci est, comme dans
le cas pr\'ec\'edent, absurde
car ce $\PP ^1$ se d\'eformerait
alors dans $X$ hors de $Y$ ! \finpreuve

\subsection{Quelques commentaires}
Comme nous venons de le voir, la situation en dimension
$4$ est tr\`es satisfaisante lorsque $K_X$ n'est pas nef.
Dans le cas o\`u $K_X$ est nef, nous avons obtenu
une restriction sur le centre de l'\'eclatement seulement
lorsque le diviseur exceptionnel $E$ est contract\'e sur un point.
Au moment o\`u nous finissions la r\'edaction de cette th\`ese,
nous avons appris que M.\ Andreatta et J.A.\ Wi\'sniewski
terminent la r\'edaction d'un travail consistant \`a classifier
les contractions extr\^emales divisorielles en dimension $4$
sur une vari\'et\'e non-singuli\`ere,
\'etendant ainsi les r\'esultats de M.\ Beltrametti au cas
o\`u le diviseur est contract\'e sur une courbe
ou sur une surface. Nous sommes en mesure d'appliquer leurs
r\'esultats dans notre situation pour obtenir la
proposition suivante. Pr\'ecisons cependant que nous n'avons
pas encore une version \'ecrite du travail en question
mais que notre seule r\'ef\'erence est une s\'erie
de discussions informelles avec M.\ Andreatta, M.\ Mella et J.A.\ Wi\'sniewski.

\medskip

\noindent {\bf Proposition  }
{\em Soit $X$ comme dans le th\'eor\`eme H. On suppose
que $K_X$ est nef. Si $f$ est la contraction de Mori
d\'efinie sur $\tilde{X}$, alors~:

(i) le diviseur exceptionnel $E$ est contract\'e
sur une courbe ou un point. Autrement dit, $f(E)$
n'est pas une surface,

(ii) si $f(E)$ est une courbe, cette derni\`ere est une
courbe lisse de singularit\'es nodales ordinaires
$3$-dimensionelles et le centre $Y$
de l'\'eclatement $\pi$ est une surface r\'egl\'ee
dont les fibres $\PP ^1$ ont pour fibr\'e
normal ${\cal O}_{\PP ^1}\oplus {\cal O}_{\PP ^1}(-1)^{\oplus 2}$.
Autrement dit, la situation est localement le
produit d'une courbe par le mod\`ele analogue en dimension~$3$.}

\medskip

Cette proposition termine la description des
situations possibles~; cependant nous ne connaissons
pas \`a l'heure actuelle d'exemple explicite o\`u le
point (ii) est r\'ealis\'e.

\medskip

\noindent {\bf ``D\'emonstration"}

Tout d'abord, mentionnons que la contraction divisorielle que
nous \'etudions est tr\`es
particuli\`ere car nous savons {\em a priori} que le diviseur
exceptionnel a une structure de fibration en espaces
projectifs sur une base lisse.

Pour le point (i), supposons par l'absurde que $f(E)$ est une
surface. Dans ce cas, la fibre g\'en\'erale est un $\PP ^1$
et M.\ Andreatta et J.A.\ Wi\'sniewski montrent qu'une
\'eventuelle fibre particuli\`ere est soit $\PP ^2$, soit
la quadrique ${\cal Q}_2$, soit la quadrique
singuli\`ere ${\cal Q}_2^0$. Dans notre situation,
une \'eventuelle fibre particuli\`ere est donc ${\cal Q}_2$
et l'image $\pi ({\cal Q}_2)$ dans $Y$ est une courbe
rationnelle $C$ d'auto-intersection $-1$.

\noindent Montrons que ceci n'est pas possible, \`a
nouveau par un argument de d\'eformation. En effet,
$K_X$ est trivial sur $C$, donc
$$ N_{C/X} =  {\cal O}_{\PP ^1}(-1) \oplus {\cal O}_{\PP ^1}(a) \oplus
{\cal O}_{\PP ^1}(b)$$
o\`u $a$ et $b$ sont deux entiers satisfaisant la relation $a + b = -1$.
De l\`a
$$ \dim \Hilb _{[C]}(X) \geq \dim H^0(C,N_{C/X}) - \dim H^1(C,N_{C/X}) =
a+b+2 = 1 > 0$$
d'o\`u l'on d\'eduit que $C$ se d\'eforme dans $X$ et ce hors de $Y$.

Pour le point (ii), nous sommes dans la situation ``facile"
du travail de M.\ Andreatta et J.A.\ Wi\'sniewski car les fibres
de $f$ restreinte \`a $E$ sont \'equi-dimensionnelles.
Dans notre situation, la fibre g\'en\'erale est une quadrique
${\cal Q}_2$ et il n'y a pas de fibres particuli\`eres~: la
situtation est, transversalement \`a $f(E)$, la r\'esolution
d'une singularit\'e nodale $3$-dimensionnelle.\finpreuve


\end{document}